\UseRawInputEncoding
\documentclass[11pt]{article}
\usepackage{amsfonts,mathrsfs}
\usepackage{amssymb}
\usepackage{amsmath}
\usepackage{color}
=0pt
\def\sqr#1#2{{\vcenter{\vbox{\hrule height.#2pt
              \hbox{\vrule width.#2pt height#1pt \kern#1pt \vrule
width.#2pt}
              \hrule height.#2pt}}}}
\def\signed #1{{\unskip\nobreak\hfil\penalty50
              \hskip2em\hbox{}\nobreak\hfil#1
              \parfillskip=0pt \finalhyphendemerits=0 \par}}
\def\endpf{\signed {$\sqr69$}}
\def\dbE{{\mathbb{E}}}
\def\dbF{{\mathbb{F}}}

\def\dbN{{\mathbb{N}}}

\def\dbP{{\mathbb{P}}}

\def\dbR{{\mathbb{R}}}
\def\dbS{{\mathbb{S}}}

\def\a{\alpha}

\def\d{\delta}
\def\e{\varepsilon}

\def\l{\lambda}

\def\si{\sigma}

\def\f{\varphi}

\def\om{\omega}

\def\3n{\negthinspace \negthinspace \negthinspace }
\def\2n{\negthinspace \negthinspace }
\def\1n{\negthinspace }
\def\ns{\noalign{\smallskip} }

\def\ds{\displaystyle}
\def\G{\Gamma}
\def\D{\Delta}
\def\Th{\Theta}
\def\L{\Lambda}

\def\Om{\Omega}
\def\om{\omega}
\def\cA{{\cal A}}

\def\cC{{\cal C}}
\def\cD{{\cal D}}

\def\cF{{\cal F}}

\def\cJ{{\cal J}}

\def\cL{{\cal L}}

\def\cO{{\cal O}}

\def\cS{{\cal S}}
\def\cT{{\cal T}}

\def\cX{{\cal X}}
\def\cY{{\cal Y}}

\def\cl{{\cal l}}
\def\mE{{\mathbb{E}}}

\def\no{\noindent}

\def\ss{\smallskip}
\def\ms{\medskip}
\def\bs{\bigskip}
\def\q{\quad}
\def\qq{\qquad}
\def\hb{\hbox}

\def\limsup{\mathop{\overline{\rm lim}}}

\def\lan{\big\langle}
\def\ran{\big\rangle}

\def\pa{\partial}

\def\wt{\widetilde}
\def\cd{\cdot}

\def\ae{\hbox{\rm a.e.{ }}}

\def\cl{\overline}

\def\deq{\mathop{\buildrel\D\over=}}

\def\({\Big (}
\def\){\Big )}
\def\[{\Big[}
\def\]{\Big]}

\def\={\buildrel \triangle \over =}

\def\resp{{\it resp. }}

\def\be{\begin{equation}}
\def\bel{\begin{equation}\label}
\def\ee{\end{equation}}
\def\bea{\begin{eqnarray}}
\def\eea{\end{eqnarray}}
\def\bt{\begin{theorem}}
\def\et{\end{theorem}}
\def\bc{\begin{corollary}}
\def\ec{\end{corollary}}
\def\bl{\begin{lemma}}
\def\el{\end{lemma}}
\def\bp{\begin{proposition}}
\def\ep{\end{proposition}}
\def\br{\begin{remark}}
\def\er{\end{remark}}
\def\ba{\begin{array}}
\def\ea{\end{array}}
\def\bde{\begin{definition}}
\def\ede{\end{definition}}

\newtheorem{lemma}{Lemma}[section]
\newtheorem{remark}{Remark}[section]

\newtheorem{theorem}{Theorem}[section]
\newtheorem{corollary}{Corollary}[section]

\newtheorem{definition}{Definition}[section]
\newtheorem{proposition}{Proposition}[section]

\allowdisplaybreaks
 \makeatletter
   
   \@addtoreset{equation}{section}
\makeatother

\begin{document}

\title{\bf Optimal Feedback Controls of Stochastic Linear Quadratic Control Problems in Infinite Dimensions with Random Coefficients
}

\author{Qi L\"{u}\footnote{School of
Mathematics, Sichuan University, Chengdu, P. R. China. Email: lu@scu.edu.cn. Qi L\"u is supported by the NSF of China under grants 12025105 and 11971334}  ~~~ and ~~~ Tianxiao Wang \footnote{School of
Mathematics, Sichuan University, Chengdu, P. R. China. Email: wtxiao2014@scu.edu.cn. Tianxiao Wang is supported by the NSF of China under grants 11971332 and 11931011.}~\footnote{Corresponding author. } }

\date{}

\maketitle

\begin{abstract}

It is a longstanding unsolved problem to
characterize the optimal feedback controls for
general linear quadratic optimal control problem
of stochastic evolution equation with random
coefficients. A solution  to this problem is
given in \cite{Lu-Zhang-arxiv-2019} under some
assumptions which can be verified for
interesting concrete models, such as controlled
stochastic wave equations, controlled stochastic
Schr\"odinger equations, etc. More precisely,
the authors establish the equivalence between
the existence of optimal feedback operator and
the solvability of the corresponding
operator-valued, backward stochastic Riccati
equations.  However, their result cannot cover
some important stochastic partial differential
equations, such as stochastic heat equations,
stochastic stokes equations, etc. A key
contribution of the current work is to relax the
$C_0$-group assumption of unbounded linear
operator $A$ in \cite{Lu-Zhang-arxiv-2019} and
using contraction semigroup assumption instead.
Therefore, our result can be well applicable in
the linear quadratic problem of stochastic
parabolic partial differential equations. To
this end, we introduce a suitable notion to the
aforementioned Riccati equation, and some
delicate techniques which are even new in the
finite dimensional case.

\end{abstract}

\bs

\no{\bf 2010 Mathematics Subject
Classification}. 49N10, 49N35  

\bs

\no{\bf Key Words:} Stochastic linear quadratic optimal control problem, stochastic evolution equation, random coefficients, backward stochastic Riccati equation, transposition solution

\section{Introduction}\label{s1}

In Control Theory, one of the fundamental issues
is to find optimal feedback controls, which are
particularly important in practical applications
since the main advantage of feedback controls is
to keep the corresponding control to be robust
with respect to (small) perturbation/
disturbance, which is usually unavoidable in
real world. Unfortunately, it is actually very
difficult to find optimal feedback controls  for
many optimal control problems. So far, the most
successful attempt in this respect is that for
linear quadratic control problems (LQ problems
for short) which are extensively studied in
Control Theory.  The study of LQ problems dates
back at least to
\cite{Bellman-Glicksberg-Gross-1958} in which
the system is governed by a linear ordinary
differential equation. Later, in the seminal
work of \cite{Kalman-1960}, the matrix-valued
Riccati equations were brought into LQ problems
to construct (linear) feedback controls. It is
well-known that, Kalman's theory for LQ problems
is among the three milestones in modern (finite
dimensional) optimal control theory. Due to the
elegant and fruitful mathematical structure,  LQ
problems were investigated extensively in the
literature for a variety of deterministic
control systems (e.g.,
\cite{Anderson-Moore-1989, Bensoussan-2007}).

LQ problems for controlled stochastic
differential equations (SDEs) was first studied
in \cite{Wonham-1968}. More precisely, given
proper $(s,\eta)$, the state equation is
described as
\bel{SDE-state-equation}
\left\{\ba{ll}
\ns\ds d X(t)=\big[A X(t)+B u(t)\big]dt+\big[C X(t)+D u(t)\big]dW(t),\ \ t \in (s,T],\\
\ns\ds X(s)=\eta,
\ea\right. \ee
and the cost functional is defined as
\bel{SDE-cost-functional}
 \ba{ll}
\ns\ds \cJ(s,\eta;u(\cd))=\frac 1 2 \dbE\Big[\int_s^T \big[\lan Q X(t),X(t)\ran
_{\dbR^n}+\lan R u(t),u(t)\ran_{\dbR^m}\big]dt+\lan GX(T),X(T)\ran_{\dbR^n}\Big].
\ea \ee
In the above, $W(\cd)$ is a standard Brownian
motion defined on complete filtered probability
space $(\Omega,\cF,\dbF,\dbP)$ with
$\dbF=\{\cF_t\}_{t\in[0,T]}$ generated by
$W(\cd)$ augmented by all $\dbP$-null sets. The
coefficients $A$, $B$, $C$, $D$, the control
variable $u(\cd)$, and the state variable
$X(\cd)$ are stochastic processes with suitable
measurability and integrability such that both
(\ref{SDE-state-equation}) and
(\ref{SDE-cost-functional}) are well-defined.
The optimal control problem is to find (if
possible) a $\bar u(\cd)$ to minimize the cost
functional.  There are a huge amount works
addressing the  LQ problems for controlled SDEs
(see \cite{Bismit-1978,
    Chen-Li-Zhou-1998-SICON,Lu-Wang-Zhang-PUQR,
    Sun-Li-Yong-2016-SICON,
    Sun-Yong-2020,Tang-2003-SICON,
    Wonham-1968} and the
rich references therein). To simplify the notations, here and henceforth the sample point $\omega(\in \Omega) $ and/or the
time variable $t(\in [0,T])$ in the coefficients are often suppressed, in the case that no
confusion would occur.

Similar as the deterministic setting, to find the optimal feedback
controls, people introduce the following Riccati equation:
\begin{equation}\label{Riccati-SDEs-1}
 \left\{
\begin{array}{ll}\ds
dP =-\big[ P A +
A^{\top} P + \L C + C^{\top} \L + C^{\top} PC  + Q - L^{\top} K^{\dagger} L
\big]dt+ \L dW(t) \q&\mbox{in }[0,T),\\
\ns\ds P(T)=G,
\end{array}
\right.
\end{equation}
where $K\equiv R+D^{\top}PD,$ $L=
B^{\top} P+D^{\top} (PC+\L),$ and
$K^{\dagger}$ is the Moore-Penrose
pseudo-inverse of $K$, is introduced as
fundamental tools to represent the
linear feedback controls.
We point out
that (\ref{Riccati-SDEs-1}) is actually
a backward stochastic differential
equation (BSDE), whose solvability is
highly challenging and nontrivial. Even
though the well-posedness is
guaranteed, as pointed in \cite[Remark
1.2]{Lu-Wang-Zhang-PUQR}, the solution is not fully qualified to
serve as the design of optimal feedback
controls. Naturally, one may ask a
question:

\it Is it possible to link
the existence of optimal feedback
control for LQ problems with the solvability of
(\ref{Riccati-SDEs-1})?

\rm Along this
line, an affirmative answer was given
in \cite{Lu-Wang-Zhang-PUQR}.

Now let us move to the linear quadratic  for
stochastic evolution equations (SEEs), which is
the main concern of the current paper. Let $H$
and $U$ be separable Hilbert spaces, and $A$ be
an unbounded linear operator (with domain $D(A)$
on $H$), which generates a $C_0$-semigroup
$\{e^{At}\}_{t\geq 0}$ on $H$. Denote by $A^*$
the adjoint operator of $A$, which generates the
adjoint $C_0$-semigroup $\{e^{A^*t}\}_{t\geq 0}$
on $H$. Given proper $(s,\eta)$, consider the
following controlled linear SEE:
\begin{equation}\label{SEE-state-equation}
\left\{\2n\begin{array}{ll}\ds
dx(t)=\big[(A+A_1) x  + B u \big]dt + \big(C
x +D
u \big)dW(t) \q \mbox{in }(s,T],\\
\ns\ds x(s)=\eta,
\end{array}
\right.
\end{equation}
with the quadratic cost functional
\begin{equation}\label{SEE-cost-functional}
\begin{array}{ll}\ds \cJ(s,\eta;u(\cd))
\3n&\ds=\frac{1}{2}\mE\Big[ \int_s^T
\big(\big\langle Q x(t),x(t)\big\rangle_H
+\big\langle R u(t),u(t)\big\rangle_U\big)dt +
\langle Gx(T),x(T)\rangle_H\Big].
\end{array}
\end{equation}
Here $ u(\cd) $  is the control variable,
$x(\cd)$ is the state variable, and the
coefficients (say, $A_1$, $B$, $C$ and $D$) are
stochastic processes. Under proper conditions
(See Section 2), (\ref{SEE-state-equation})
admits a unique solution and
(\ref{SEE-cost-functional}) becomes
well-defined. The optimal control problem is to
find (if possible) a control $\bar  u(\cd) $ to
minimize $\cJ$ in (\ref{SEE-cost-functional}).

SEEs are used to describe a lot of random
phenomena  appearing in physics, chemistry,
biology, and so on. In many situations SEEs are
more realistic mathematical models than the
deterministic ones (e.g.
\cite{Carmona-Rozovskii-1999,
Kotelenez-2008,Lu-Zhang-book1}). Consequently,
there are many works addressing the optimal
control problems for SEEs. In particular, we
refer the readers to \cite{Ahmed-1981-SICON,
Hafizoglu-et-al-2017-SICON,Lu-2019-JDE,Lu-Zhang-arxiv-2019,Lu-Zhang-book1}
and the rich references therein for LQ problems
for controlled SEEs.

To design optimal feedback controls for linear quadratic problem of
SEEs, similar to \eqref{Riccati-SDEs-1}, we introduce the following
operator-valued Riccati equation:
\begin{equation}\label{BSRE-SEEs}
 \left\{
\begin{array}{ll}\ds
dP =-\big[ P(A+A_1) +
(A+A_1)^* P + \L C + C^* \L\\
\ns\ds\qq\;\;\q + C^* PC  + Q - L^* K^{-1} L
\big]dt+ \L dW(t) \q&\mbox{in }[0,T),\\
\ns\ds P(T)=G,
\end{array}
\right.
\end{equation}
where $K\equiv R+D^*PD$, $L= B^* P+D^* (PC+\L).$
Besides the difficulties for the case
of finite dimensions, there exists
several new essentially ones in the
study of (\ref{BSRE-SEEs}) when
$\hbox{dim}H = \infty$:
\begin{itemize}
    \item There is no
    stochastic integration/evolution
    equation theory in general Banach
    spaces that can be employed to handle
    the (stochastic integral) term ``$\L
    dW(t)$'' in (\ref{BSRE-SEEs})
    effectively. Due to this difficulty,
    there exist only a quite limited number
    of works, see \cite{Ahmed-1981-SICON,
        Guatteri-Tessitor-2005-SICON,
        Guatteri-Tessitor-2014-SICON} and the
    papers therein. We refer to \cite{Lu-Zhang-arxiv-2019} for the detailed introduction.

    \item The appearance of $D$, and the resulting term $L^* K^{-1} L$ lead to another essential difficulty in the study of (\ref{BSRE-SEEs}). Notice that the second term is already quite subtle even in the finite dimensional case. Usually, people assume that $D=0$ to avoid that difficulty (e.g., \cite{        Guatteri-Tessitor-2005-SICON,
        Guatteri-Tessitor-2014-SICON}). We study the case that $D\neq 0$ according to the following reasons:

    On one hand, only when the controls/decisions could or would influence the scale of uncertainty do the stochastic problems differ from the deterministic
ones. On the other hand, once one put a control in the drift term, it will influence the diffusion term.
\end{itemize}

Generally speaking, to study the difficult
operator-valued stochastic differential
equations, people need to introduce suitable new
concept of solution. For example, in
\cite{Lu-Zhang-2018-MCRF, Lu-Zhang-book}, the
authors introduced the concept of transposition
solution  to operator-valued, backward
stochastic (linear) Lyapunov equations to study
the maximum principle for optimal control
problems of SEEs. In \cite{Lu-Zhang-arxiv-2019},
the authors defined a type of transposition
solution to treat (\ref{BSRE-SEEs}), and
obtained equivalence between the existence of
optimal feedback operator and the solvability of
BSRE (\ref{BSRE-SEEs}) in the sense of
transposition solution. We point out that such a
definition of transposition solution can bypass
the previous two difficulties. However, due to
the methodology limitation in obtaining BSREs by
the existence of optimal feedback operators, the
authors made a common yet technical assumption
of $A$, i.e., $A$ generates a $C_0$-group on
$H$. In fact, they followed the forward-backward
stochastic system ideas (see e.g.
\cite{Lu-Wang-Zhang-PUQR, Sun-Yong-2014-SICON}
for the finite dimensional case), and therefore
used the inverse of solution to certain forward
stochastic differential equation to construct
the solution of Riccati equation. In the finite
dimensional case, the inverse of the solution to
the forward stochastic system is easy to see.
However, in the infinite dimensions, to achieve
such a goal, it requires $-A^*$ generates a
$C_0$-semigroup as well, which means that $A^*$
generates a $C_0$-group.  As a result, their
equivalence ruled out some important stochastic
partial differential equations, say  stochastic
parabolic equation and stochastic stokes
equation for example.

In this paper, by supposing $A$ generating a contraction semigroup,
we establish the equivalence between optimal feedback controls and
backward stochastic Riccati equations in infinite dimensions.  As a
result, we not only have to deal with the two difficulties mentioned
in the last paragraph, but also need to overcome the methodology
limitation in \cite{Lu-2019-JDE, Lu-Zhang-arxiv-2019,
Lu-Zhang-book1}. Notice that all these obstacles can be well-treated
in the finite dimensional framework. To figure out these issues, let
us look at the following observation which holds in the finite
dimensional framework: the BSRE (\ref{Riccati-SDEs-1}) can be
rewritten as a linear BSDE with a parameter $\Th$:
\begin{equation}\label{BSRE-BSDE-1-transform}
\left\{\begin{array}{ll} \ds
dP=-\big[P(A+B\Th)+(A+B\Th)^{\top}P+(C+D\Th)^{\top}P(C+D\Th)+(C+D\Th)^{\top}\L\\
\ns\ds\qq\qq
+\L (C+D\Th) +Q + \Th^{\top}R\Th\big]dt  +\L dW(t) \qq\qq\qq \qq  \mbox{ in } [0,T],\\
\ns\ds P(T)=G,
\end{array}\right.
\end{equation}
where the parameter $\Th$ satisfies the following:
\bel{Constraint-BSRE-BSDE}\ba{ll}
\ns\ds D^{\top}PC+B^{\top}P+D^{\top}\L+(R+D^{\top}PD)\Th=0.
\ea\ee
This important fact was used in several existing papers under Markovian framework, see e.g. \cite[Theorem 3.3]{Li-Sun-Yong-2016-PUQR} and \cite[Theorem 2.3]{Wang-COCV-2020} for the mean-field stochastic LQ problems for SDEs,   and \cite[Lemma 2.8]{Lu-2020-COCV} for the mean-field stochastic LQ problem for SEEs.
Clearly, one advantage of (\ref{BSRE-BSDE-1-transform}) lies in its linearity in contrast with (\ref{Riccati-SDEs-1}). As a tradeoff, one has to explore the additional relation (\ref{Constraint-BSRE-BSDE}).
In this article, we successfully apply this idea
into our infinite dimensional framework with
random coefficients and derive the new
equivalence between optimal  feedback operators
and BSREs in the following three steps:
\begin{enumerate}
    \item  We establish the well-posedness of an infinite dimensional version of (\ref{BSRE-BSDE-1-transform}), i.e., operator-valued backward stochastic Lyapunov equation (BSLE),   in the sense of $H_\l$-transposition solution (inspired by \cite{Lu-Zhang-2018-MCRF}) (see Subsection 4.1).
    \item  We provide a relationship
    between the optimal feedback  operator and the
    operator-valued BSLE
    (see Subsection 4.2). More precisely,
     with
    subtle choice of $v$, small $\e>0$ such
    that $t+\e\leq T$, we introduce control
    variables $\Th \bar X$ and $\Th
    X^{v,\e}+vI_{[t,t+\e]}$ associated with
    optimal feedback operator $\Th$. At this moment, we point out one novel technical contribution. We refer to Remark \ref{Remark-novel-contribution} for more relevant details.
    \item Using the result established in step 2, we prove that the aforementioned $H_\l$-transposition solution of BSLE is just the transposition solution of BSREs.

\end{enumerate}

At this very moment, let us summarize our contributions of the current article:
\begin{itemize}
    \item  From the viewpoint of methodology, we develop some useful procedures to overcome the essential obstacles and limitation in the current relevant literature. Interestingly, some technologies happen to be new even in finite dimensional scenario.

    \item From the viewpoint of conclusions, we considerably relax the $C_0$-group assumption of operator $A$ and impose the $C_0$-semigroup assumption instead. This allows us to cover LQ problem for stochastic parabolic PDEs, which is one of the most important class of SPDEs.

\end{itemize}

It is worth noting that, besides the study of
optimal feedback controls for LQs,
operator-valued Riccati equations have other
applications (e.g.,\cite{Schuch}). We believe
that both the results and the methods in this
paper have other applications, such as the study
of open quantum systems with continuation
observation. Nevertheless, a detailed study for
it is beyond the scope of this paper and will be
done in our future works.

The rest of this paper is organized  as follows.
In Section 2, we present some necessary
notations, important remarks and a complete
formulation of our optimal control problems. In
Section 3, we state our main result, i.e., the
characterization of optimal feedback operators
by means of well-posedness of operator-valued
backward stochastic Riccati equation in the
suitable transposition sense. In Section 4, we
first establish the well-posedness of
operator-valued backward stochastic Lyapunov
equation  (BSLE for short) in Subsection 4.1,
then give a  characterization of optimal
feedback operator via the operator-valued BSLE
in Subsection 4.2 and then demonstrate the proof
of the main result. In Section 5, we present an
illustrative example for  LQ problem of
stochastic parabolic PDEs.

\section{Preliminary notations}
Let $T>0$ and $(\Om,\cF,\mathbf{F},\dbP)$ be a
complete filtered probability space (satisfying
the usual conditions), on which a
$1$-dimensional standard Brownian motion
$\{W(t)\}_{t\in[0,T]}$ is defined. Here
$\dbF\deq\{\cF_t\}_{t\in[0,T]}$ is the natural
filtration generated by $W(\cd)$. Denote by
$\dbF$ the progressive $\si$-field (in
$[0,T]\times\Omega$) with respect to
$\mathbf{F}$.

Let $\cX$ be a Banach space. For any $t\in[0,T]$
and $p\in [1,\infty)$, denote by
$L_{\cF_t}^p(\Om;\cX)$ the Banach space of all
$\cF_t$-measurable random variables $\xi:\Om\to
\cX$ such that $\mathbb{E}|\xi|_\cX^p < \infty$,
with the canonical  norm. Denote by
$L^{p}_{\dbF}(\Om;C([t,T];\cX))$ the Banach
space of all $\cX$-valued $\mathbf{F}$-adapted
continuous processes $\phi(\cdot)$, with the
norm\vspace{-3mm}
$$
|\phi(\cd)|_{L^{p}_{\dbF}(\Om;C([t,T];\cX))} \=
\[\mE\sup_{\tau\in
[t,T]}|\phi(\tau)|_\cX^p\]^{1/p}.  $$
Similarly, one can define $L^{p}_{\dbF}(\Om;C([\tau_1,\tau_2];\cX))$
for two stopping times $\tau_1$ and $\tau_2$ with $\tau_1\leq
\tau_2$, $\dbP$-a.s. Also, denote by
$C_{\dbF}([t,T];L^{p}(\Om;\cX))$ the Banach space of all
$\cX$-valued $\mathbf{F}$-adapted processes $\phi(\cdot)$ such that
$\phi(\cdot):[t,T] \to L^{p}_{\cF_T}(\Om;\cX)$ is continuous, with
the norm
$$
|\phi(\cd)|_{C_{\dbF}([t,T];L^{p}(\Om;\cX))} \=
\sup_{\tau\in
[t,T]}\left[\mE|\phi(\tau)|_\cX^p\right]^{1/p}.
$$
Fix any $p_1,p_2,p_3,p_4\in[1,\infty]$.
Put
$$
\begin{array}{ll}
\ds L^{p_1}_\dbF(\Om;L^{p_2}(t,T;\cX))
=\Big\{\f:(t,T)\times\Om\to
\cX\;\Big|\;\f(\cd)\hb{
is $\mathbf{F}$-adapted and }\dbE\(\int_t^T|\f(\tau)|_\cX^{p_2}d\tau\)^{\frac{p_1}{p_2}}<\infty\Big\},\\
\ns\ds
 L^{p_2}_\dbF(t,T;L^{p_1}(\Om;\cX)) =\Big\{\f:(t,T)\times\Om\to
\cX\;\Big|\;\f(\cd)\hb{ is $\mathbf{F}$-adapted
and
}\int_t^T\(\dbE|\f(\tau)|_X^{p_1}\)^{\frac{p_2}
{p_1}}d\tau<\infty\Big\}.
 \end{array}
 $$
(When any one of $p_j$ ($j=1,2,3,4$) is equal to
$\infty$, it is needed to make the usual
modifications in the above definitions of
$L^{p_1}_\dbF(\Om;L^{p_2}(t,T;\cX))$ and
$L^{p_2}_\dbF(0,T;L^{p_1}(\Om;\cX))$). If
$\cX=\dbR$, then we may omit the term $\cX$ in
the above spaces. Clearly, both
$L^{p_1}_\dbF(\Om;L^{p_2}(t,T;\cX))$ and
$L^{p_2}_\dbF(t,T;L^{p_1}(\Om;\cX))$ are Banach
spaces with the canonical norms. If $p_1=p_2$,
we simply write the above spaces as
$L^{p_1}_\dbF(t,T;\cX)$. For $r\in[0,T]$ and
$f\in L^1_{\cF_T}(\Om;\cX)$, denote by $\dbE_r
f$ the conditional expectation of $f$ with
respect to $\cF_r$ and by $\dbE f$ the
mathematical expectation of $f$.

\ms
Let $\cY$ be another Banach space. Denote by
$\cL(\cX; \cY)$ the Banach space of all bounded
linear operators from $\cX$ to $\cY$, with the
usual operator norm (When $\cY=\cX$, we simply
write $\cL(\cX)$ instead of $\cL(\cX; \cX)$).
Suppose $\cX_j$ and $\cY_j$ ($j=1,2$) are Banach
spaces satisfying $\cX_1\subset\cX\subset\cX_2$
and $\cY_1\subset\cY\subset\cY_2$. If
$M\in\cL(\cX; \cY)$ can be extended as an
operator $\wt M\in\cL(\cX_2; \cY_2)$, then, to
simplify the notations, (formally) we also write
$M\in\cL(\cX_2; \cY_2)$. Similarly, if
$M|_{\cX_1}\in\cL(\cX_1; \cY_1)$, then, we write
$M\in \cL(\cX_1; \cY_1)$.

\ms

Throughout this paper, for any operator-valued
process/random variable $M$, we write $M^*$ for
its pointwise dual. For example, if $M\in
L^{r_1}_\dbF(0,T; L^{r_2}(\Om; \cL(H)))$, then
$M^*\in L^{r_1}_\dbF(0,T; L^{r_2}(\Om;
\cL(H)))$, and $|M|_{L^{r_1}_\dbF(0,T;
L^{r_2}(\Om; \cL(H)))}=|M^*|_{L^{r_1}_\dbF(0,T;
L^{r_2}(\Om; \cL(H)))}$. When $\cX$ is a Hilbert
space, denote by  $\cS(\cX)$ the set of all
self-adjoint operators on $\cX$. Further, fix
any $p_3,p_4\in[1,\infty]$, we put
\begin{equation}\label{3.14-eq4}
\begin{array}{ll}\ds
\mathbf{L}\big(L^{p_2}_{\dbF}(\Om;L^{p_1}(t_1,t_2;\cX));\;L^{p_4}_{\dbF}(\Om;L^{p_3}(t_1,t_2;\cY))
\big)
\\[1mm]
\triangleq \big\{\mathrm{L}\in
\cL\big(L^{p_2}_{\dbF}(\Om;L^{p_1}(t_1,t_2;\cX));L^{p_4}_{\dbF}(\Om;L^{p_3}(t_1,t_2;\cY))\big)\,\big|
\mbox{ for a.e. } \\[1mm]
\ns\ds\q (t,\omega)\in
(t_1,t_2)\times\Omega, \mbox{
    there exists } L(t,\omega)\in\cL (\cX;\cY)
\mbox{ verifying } \\[1mm]
\ns\ds\q\mbox{ that }  \big(\mathrm{L}
f(\cd)\big)(t,\omega)=L
(t,\omega)f(t,\omega),\; \forall\;
f(\cd)\in
L^{p_2}_{\dbF}(\Om;L^{p_1}(t_1,t_2;\cX))\}
\end{array}
\end{equation}
and
\begin{equation}\label{3.14-eq5}
\begin{array}{ll}\ds
\mathbf{L}\big(\cX;\;L^{p_3}_{\dbF}(t_1,t_2;
L^{p_4}(\Om;\cY))\big)
\\[1mm]
\triangleq \big\{\mathrm{L}\in
\cL\big(\cX;L^{p_3}_{\dbF}(t_1,t_2;L^{p_4}(\Om;\cY))\big)\,\big|\,
\mbox{for a.e.  } (t,\omega)\in
(t_1,t_2)\times\Omega,\\[1mm]
\ns\ds\q \mbox{ there }\mbox{ exists }
L(t,\omega)\in\cL (\cX;\cY) \mbox{
verifying that } \big(\mathrm{L}
x\big)(t,\omega)=L (t,\omega)x, \;\;
\forall\; x \in \cX\big\}.
\end{array}
\end{equation}
To simplify the notations, in what
follows we shall identify the above
$\mathrm{L}$ with $L(\cd,\cd)$.
Similarly, one can define the spaces
$\mathbf{L}\big(L^{p_2}(\Om;\cX); L^{p_4}_{\dbF}(\Om;L^{p_3}(t_1,t_2;\cY))\big)$
and $\mathbf{L}\big(L^{p_2}(\Om;\cX); L^{p_4}(\Om;\cY)\big)$, etc. In the
sequel, we shall call an element in the
sets
$\mathbf{L}\big(L^{p_2}_{\dbF}(\Om;L^{p_1}(t_1,t_2;\cX)); L^{p_4}_{\dbF}(\Om;L^{p_3}(t_1,t_2;\cY))\big)$,
$\mathbf{L}\big(\cX; L^{p_3}_{\dbF}(t_1,t_2;
L^{p_4}(\Om;\cY))\big)$,
$\mathbf{L}\big(L^{p_2}(\Om;\cX); L^{p_4}_{\dbF}(\Om;L^{p_3}(t_1,t_2;\cY))\big),
\mathbf{L}\big(L^{p_2}(\Om;\cX);
L^{p_4}(\Om;\cY)\big)$ and so on a
pointwise defined operator.

\ss

Put
\begin{equation}\label{10.10-eq30}
\begin{array}{lll}\ds
\Upsilon_p(\cX;\cY) \ds \triangleq
\big\{L(\cd,\cd)\in
\mathbf{L}\big(L^2_{\dbF}(\Om;L^\infty(0,T;\cX));L_{\dbF}^2(\Om;L^p(0,T;\cY))\big) |\\[1mm]
\ns \ds \qq\qq\qq\qq\qq\qq
|L(\cd,\cd)|_{\cL(\cX;\cY)}\in
L^\infty_\dbF(\Om;L^p(0,T))\big\}.
\end{array}
\end{equation}
We shall simply denote
$\Upsilon_p(\cX;\cX)$ by
$\Upsilon_p(\cX)$.

\begin{lemma}\label{3.4-lm1}
Assume that $\cX$ and $\cY$ are separable
Hilbert spaces.  For any $L\in
\Upsilon_p(\cX;\cY)$, there exist a sequence
$\{L_n\}_{n=1}^\infty\subset
L_{\dbF}^\infty\big(\Om;L^p(0,T;\cL(\cX;\cY))\big)$
such that
\begin{equation}\label{3.4-lm1-eq1}
\lim_{n\to\infty} L_n v = Lv \q \mbox{ in
}L_{\dbF}^2(\Om;L^p(0,T;\cY)),\q \forall\, v\in
L^2_{\dbF}(\Om;L^\infty(0,T;\cX)).
\end{equation}
\end{lemma}

{\it Proof}. Since $\cX$ (\resp $\cY$) is a
separable Hilbert space, there is an orthomormal
basis $\{e_k\}_{k=1}^\infty$ (\resp $\{\tilde
e_k\}_{k=1}^\infty$)  of $\cX$
 (\resp $\cY$). Denote by $\cX_n$  (\resp $\cY_n$) the
subspace spanned by $\{e_k\}_{k=1}^n$(\resp
$\{\tilde e_k\}_{k=1}^n$) and by $\G_n$  (\resp
$\wt \G_n$) the orthogonal projection operator
from $\cX$  (\resp $\cY$) to $\cX_n$  (\resp
$\cY_n$).

Let $L_n = \wt \G_n L \G_n$. We first prove that
$\{L_n\}_{n=1}^\infty\subset
L_{\dbF}^\infty(\Om;L^p(0,T;\cL(\cX;\cY)))$.

Since $L_n(t,\om)\in \cL(\cX;\cY)$ for
a.e. $(t,\om)\in [0,T]\times\Om$, we
have  $L_n(t,\om)\in \cL(\cX;\cY)$ for
a.e. $(t,\om)\in [0,T]\times\Om$.
Moreover, from the definition of $L_n$,
we know that  for a.e. $(t,\om)\in
[0,T]\times\Om$, $L_n(t,\om)$ can be
regarded as an element in
$\cL(\cX_n;\cY_n)$.

Since $L(\cd,\cd)\!\in
\mathbf{L}(L^2_{\dbF}(\Om;L^\infty(0,T;\cX));L_{\dbF}^2(\Om;L^p(0,T;$
$\cY)))$, for any $x\in \cX$, we have that
$L_n(\cd,\cd)x\in
L_{\dbF}^2(\Om;L^p(0,T;\cY_n))$. Noting that
$\cL(\cX_n;\cY_n)$ is isomorphic to
$\dbR^{n\times n}$, we get that $L_n(\cd,\cd)
\in L_{\dbF}^2(\Om;L^p(0,$
$T;\cL(\cX_n;\cY_n)))$. Further, noting that
$$
|L_n(\cd,\cd)|_{\cL(\cX_n;\cY_n)}=|L_n(\cd,\cd)|_{\cL(\cX;\cY)}\leq
|L(\cd,\cd)|_{\cL(\cX;\cY)},
$$
it follows from
$|L(\cd,\cd)|_{\cL(\cX;\cY)}\in
L^\infty_\dbF(\Om;L^p(0,T))$ that
$|L_n(\cd,\cd)|_{\cL(\cX;\cY)}\in
L^\infty_\dbF(\Om;L^p(0,T))$. This,
together with $L_n(\cd,\cd) \in
L_{\dbF}^2(\Om;L^p(0,T;
\cL(\cX_n;\cY_n)))$, implies that
$L_n(\cd,\cd) \in
L_{\dbF}^\infty(\Om;L^p(0,T;
\cL(\cX_n;\cY_n)))$. Hence,
$L_n(\cd,\cd) \in
L_{\dbF}^\infty(\Om;L^p(0,T;
\cL(\cX;\cY)))$.

For any $v\in L^2_{\dbF}(\Om;L^\infty(0,T;\cX))$,
\begin{equation}\label{3.4-eq40}
\begin{array}{ll}\ds
|L_nv-Lv|_{L_{\dbF}^2(\Om;L^p(0,T;\cY))}\\
\ns\ds  = |\wt\G_n L \G_n v -
Lv|_{L_{\dbF}^2(\Om;L^p(0,T;\cY))}\\
\ns\ds \leq  |\wt\G_n L \G_n v -
\wt\G_n L
v|_{L_{\dbF}^2(\Om;L^p(0,T;\cY))} +
|\wt\G_n L v -
Lv|_{L_{\dbF}^2(\Om;L^p(0,T;\cY))}\\
\ns\ds \leq  | L \G_n v -  L
v|_{L_{\dbF}^2(\Om;L^p(0,T;\cY))} + |\wt\G_n L v
- Lv|_{L_{\dbF}^2(\Om;L^p(0,T;\cY))}.
\end{array}
\end{equation}
Since
$$
\lim_{n\to\infty}| L \G_n v -  L v|_{\cY}=0
$$
and
$$
|L \G_n v - L v|_{\cY} \leq
2|L|_{\cL(\cX;\cY)}|v|_{\cX},
$$
by Lebesgue's dominated convergence theorem, we get that
\begin{equation}\label{3.4-eq41}
\lim_{n\to\infty}| L \G_n v -  L
v|_{L_{\dbF}^2(\Om;L^p(0,T;\cY))}=0.
\end{equation}
Similarly, we can prove that
\begin{equation}\label{3.4-eq42}
\lim_{n\to\infty}|\wt\G_n L v -
Lv|_{L_{\dbF}^2(\Om;L^p(0,T;\cY))}=0.
\end{equation}
Combining \eqref{3.4-eq40}--\eqref{3.4-eq42}, we
obtain \eqref{3.4-lm1-eq1}.
\endpf

\begin{remark}\label{3.7-rmk1}
One can show that
$L^\infty_{\dbF}\big(\Om;L^p(0,T;\cL(\cX;\cY))\big)\!\subset\!
\Upsilon_p(\cX;\cY)$. For any $L(\cd,\cd)\!\in
\Upsilon_p(\cX;\cY)$, one does not need to have
$L(\cd,\cd)\in
L^\infty_{\dbF}\big(\Om;L^p(0,T;\cL(\cX;\cY))\big)$.
The reason for introducing such set is the lack
of measurability for operator-valued functions
and processes appeared in concrete problems. For
example, generally speaking, a $C_0$-semigroup
is not measurable with respect to the time
variable (e.g.,\cite[Subsection
1.1.c]{Hytonen-et-al-2016}). In such case, one
has to relax the condition on the measurability
of such functions and processes. Nevertheless,
as we shall see later, in some sense
$\Upsilon_p(\cX; \cY)$ is a nice ``replacement"
of the space
$L^\infty_{\dbF}\big(\Om;L^p(0,T;\cL(\cX;\cY))\big)$
in the study of stochastic LQ problems.
\end{remark}

\ms Now we  put the following assumptions on the
operator $A$ in the control system
\eqref{SEE-state-equation}.

\ss

({\bf AS0}) \it The $C_0$-semigroup $\{e^{At}\}_{t\geq 0}$ generated
by $A$ is contractive in the sense that $\big|S(t)\big|_{\cL(H)}\leq
e^{kt}$ for some $k\in\dbR$ and all $t\in[0,T]$.

\rm

\begin{remark}\label{3.7-rmk1-1}
    ({\bf AS0}) is not restrictive in the sense that, as far as we know, all the linear SPDEs satisfies it.
\end{remark}

({\bf AS1}) \it The eigenvectors
$\{e_j\}_{j=1}^\infty$ of $A$ such that
$|e_j|_H=1$ for all $j\in\dbN$ constitute an
orthonormal basis of $H$. \rm

\begin{remark}\label{3.7-rmk1-2}
    ({\bf AS1}) is not restrictive in the sense that, many important SPDEs, such as stochastic wave equations, stochastic heat equations and stochastic Schr\"odinger equations, involved on bounded domains fulfill it..
\end{remark}

The following assumption is put on the
coefficients of the control system
\eqref{SEE-state-equation} and the cost
functional \eqref{SEE-cost-functional}.

\vspace{0.1cm}

({\bf AS2})  {\it The coefficients
satisfy that
$$
\begin{cases}\ds
A_1(\cd)\in
L^1_\dbF(0,T;L^\infty(\Om;\cL(H))),\q
B(\cd)\in L^\infty_\dbF(\Om;L^2(0,T;
\cL(U;H))),\\
\ns\ds C(\cd)\in
L^2_\dbF(0,T;L^\infty(\Om;\cL(H))),\q D(\cd)\in L^\infty_\dbF(0,T;\cL(U;H)),\\
\ns\ds Q(\cd)\in L^\infty_\dbF(\Om;L^2(0,T; \cS(H))),\q R(\cd) \in L^\infty_\dbF(0,T;\cS(U)),\q G\in L^\infty_{\cF_T}(\Om;$ $\cS(H)),\\
\ns\ds G\geq 0, \q R>0,\q Q\geq0 \mbox{
for a.e. }(t,\om)\in [0,T]\times\Om.
\end{cases}
$$
}

Under ({\bf AS1}), the control system
\eqref{SEE-state-equation} admits a unique
solution $x(\cd)\in
C_{\dbF}([s,T];L^2(\Omega;H))$, and the cost
functional \eqref{SEE-cost-functional} becomes
well-defined. If in addition ({\bf AS0}) is
true, then $x(\cd)\in
L^2_{\dbF}(\Omega;C([s,T];H))$. Let us consider
the following optimal control problem:

\ms

\no\bf Problem (SLQ): \rm For each $(s,\eta)\in
[0,T)\times L^2_{\cF_s}(\Om;H)$, find (if
possible) a control $\bar  u(\cd)\in
L^2_\dbF(s,T;U)$ such that
\begin{equation}\label{5.2-eq3}
\cJ\big(s,\eta;\bar  u(\cd)\big)=\inf_{u(\cd)\in
L^2_\dbF(s,T;U)}\cJ\big(s,\eta;u(\cd)\big).
\end{equation}
If the above is possible, then {\bf Problem
(SLQ)} is called {\it solvable}. Any $\bar
u(\cdot)$ satisfying (\ref{5.2-eq3}) is called
an {\it optimal control}. If the $\bar u(\cd)$
which fulfills \eqref{5.2-eq3} is unique, then
{\bf Problem (SLQ)}  is called {\it uniquely
solvable}. The corresponding state $\bar
x(\cdot)$ is called an {\it optimal state}, and
$\big(\cl x(\cdot),\bar u(\cdot)\big)$ is called
an {\it optimal pair}. Under ({\bf AS2})), we
know that {\bf Problem (SLQ)}   admits at most
one optimal control (see \cite[Proposition
2.1]{Lu-Zhang-arxiv-2019}).

Set
$$
\ell^2_+\triangleq
\Big\{\{\lambda_j\}_{j=1}^\infty\in
\ell^2\;\Big|\;\lambda_j\in
(0,+\infty)\mbox{  for each } j\in\dbN
\Big\}.
$$
For any given $\l=\{\lambda_j\}_{j=1}^\infty\in \ell^2_+$, define a
norm $|\cd|_{H_\l'}$ on $H$ as follows (Recall that $|\cd|_{D(A)}$
stands for the graph norm of the operator $A$):
$$
|h|_{H_\l'} = \sqrt{\sum_{j=1}^\infty \lambda_j^2
|e_j|_{D(A)}^{-2}|h_j|^2},\qq \forall\; h=\sum_{j=1}^\infty h_j e_j
\in H.
$$
Denote by $H_\l'$ the completion of $H$
with respect to the norm
$|\cd|_{H_\l'}$. Then, $H_\l'$ is a
Hilbert space, $H\subset H_\l'$ and
$\{\lambda_j^{-1}|e_j|_{D(A)}e_j\}_{j=1}^\infty$
is an orthonormal basis of $H_\l'$.
Write $H_\l$ for the dual space of
$H_\l'$ with respect to the pivot space
$H\equiv H'$. Hence, $H_\l\subset
H\equiv H'\subset H_\l'$.

\vspace{0.1cm}

For any $\l=\{\lambda_j\}_{j=1}^\infty\in \ell^2_+$, from the
definition of $H_\l'$, we see that $\{e_j\}_{j=1}^\infty\subset
H_\l$ and the norm on $H_\l$ is given as follows:
$$
|\xi|_{H_\l} = \sqrt{\sum_{j=1}^\infty
    |\xi_j|^2|e_j|_{D(A)}^2\lambda_j^{-2}},\qq
\forall\; \xi \in H_\l,
$$
where $\xi_j=\langle \xi,e_j\rangle_H$ for $j\in\dbN$.
Furthermore,
$\{\lambda_j|e_j|_{D(A)}^{-1}e_j\}_{j=1}^\infty$
is an orthonormal basis of $H_\l$.

\begin{remark}
Note the set $H'_\lambda$ depends on the
parameter $\lambda$. Once a
$\l=\{\lambda_j\}_{j=1}^\infty\in \ell^2_+$   is
given, we get an $H_\l'$, namely,   $H_\l'$
depends on the choice of $\l$.  Generally
speaking, for two different $\l,\hat \l\in
\ell^2_+$, the corresponding spaces $H_\l'$ and
$H_{\hat\l}'$ are different.
\end{remark}

\vspace{0.1cm}

Denote by $\cL_2(H;H_\l')$ the set of all Hilbert-Schmidt operators
from $H$ to $H_\l'$. Namely,
$$
\cL_2(H;H_\l')\triangleq \Big\{F\in \cL(H;H_\l')\;\Big|\;
\sum_{j=1}^\infty |F e_j|_{H_\l'}^2<\infty \Big\}.
$$
It is well-known that $\cL_2(H;H_\l')$ is a Hilbert space itself
with the inner product
$$
\lan F_1, F_2 \ran_{\cL_2(H;H_\l')}\triangleq \sum_{j=1}^\infty \lan
F_1 e_j, F_2 e_j\ran_{H_\l'}.
$$
Denote by $\cL_2(H_\l;H)$ the set of all
Hilbert-Schmidt operators from $H_\l$ to $H$.
Similarly, $\cL_2(H_\l;H)$ is also a Hilbert
space itself. We refer to \cite{Schatten} for
more details on Hilbert-Schmidt operators.

We also need the following technical assumption:

\ms

({\bf AS3}) {\it There exists $\l=\{\lambda_j\}_{j=1}^\infty\in
\ell^2_+$ such that
\begin{equation}\label{3.8-eq1}
\left\{
\begin{array}{ll}\ds
A_1 \in L^1_\dbF(0,T;L^\infty(\Om;\cL(H_\l'))),
\q C \in
L^2_\dbF(0,T;L^\infty(\Om;\cL(H_\l'))),  \\
\ns\ds  Q\in
L^\infty_\dbF(\Om;L^2(0,T;\cL(H_\l'))),\q G\in
L^\infty_{\cF_T}(\Omega; \cL(H_\l')),\\
\ns\ds A_1 \in
L^1_\dbF(0,T;L^\infty(\Om;\cL(H_\l))), \q C \in
L^2_\dbF(0,T;L^\infty(\Om;\cL(H_\l))), \\
\ns\ds  Q\in
L^\infty_\dbF(\Om;L^2(0,T;\cL(H_\l))),\q G\in
L^\infty_{\cF_T}(\Omega; \cL(H_\l)).
\end{array}
\right.
\end{equation}
}

\begin{remark}
In Assumption {\bf (AS3)}, we assume
the existence of one $\l\in \ell^2_+$
so that \eqref{3.8-eq1} holds.
Nevertheless, there may exist another
$\hat \l \in \ell^2_+$ such that
\eqref{3.8-eq1} does not hold when $\l$
is replaced by $\hat \l$. Fortunately,
this will not influence the main
results of this paper. In other words,
we only need the existence of one
$\l\in \ell^2_+$ such that
\eqref{3.8-eq1} holds to prove our main
results.
\end{remark}
\begin{remark}
For an LQ problem of a given  SPDE,
Assumption {\bf (AS3)} means that some
coefficients in the equation and the
cost functional are smooth in some
sense (See Section \ref{sec-example} for example). This is not very restrictive
for many controlled SPDEs.
\end{remark}

Following
\cite{Lu-Zhang-arxiv-2019},
we  also introduce the notion of
optimal feedback operator.
\begin{definition}\label{5.7-def1}
An operator $\Th(\cd)\in \Upsilon_2(H;U)$ is called an {\it optimal
feedback  operator} for Problem (SLQ) if
\begin{equation}\label{5.7-eq2}
\begin{array}{ll}\ds \cJ(s,\eta;\Th(\cd)\bar x(\cd))\leq
\cJ(s,\eta;u(\cd)),\\[2mm]
\ns\ds\qq \forall\; (s,\eta)\in [0,T)\times
L^2_{\cF_s}(\Om;H),\;\; u(\cd)\in
L^2_\dbF(s,T;U),
\end{array}
\end{equation}
where $\bar x(\cd)=\bar x(\cd\,;s,\eta, \Th(\cd)\bar x(\cd))$ solves
the following equation:
\begin{equation}\label{5.2-eq1.1}
\left\{\begin{array}{ll}\ds d\bar
x(t)=\big[(A+A_1)\bar x(t) + B\Th \bar
x(t)\big]dt + \big(C \bar x(t)+D \Th \bar
x(t)\big)dW(t) &\mbox{ in
}(s,T],\\
\ns\ds \bar x(s)=\eta.
\end{array}
\right.
\end{equation}
\end{definition}

\medskip

In Definition \ref{5.7-def1}, $\Th(\cd)$ is
independent of $\eta\in H$. For a fixed $\eta\in
H$, the inequality \eqref{5.7-eq2} implies that
the control $\bar u(\cd)\equiv \Th(\cd)\bar
x(\cd)\in L^2_\dbF(0,T;U)$
is optimal for Problem (SLQ). Therefore, for
Problem (SLQ), the existence of an optimal
feedback  operator on $[0,T]$ implies the
existence of an optimal control of any
$\eta\in  H$.

\section{The statement of the main results }

To characterize the optimal feedback operator,
we need the following operator-valued, backward
stochastic Riccati equation:
\begin{equation}\label{5.5-eq6}
 \left\{
\begin{array}{ll}\ds
dP =-\big[ P(A+A_1) +
(A+A_1)^* P + \L C + C^* \L\\
\ns\ds\qq\;\;\q + C^* PC  + Q - L^* K^{-1} L
\big]dt+ \L dW(t) \q&\mbox{in }[0,T),\\
\ns\ds P(T)=G,
\end{array}
\right.
\end{equation}
where
\begin{equation}\label{9.7-eq10}
K\equiv R+D^*PD, \qq L= B^* P+D^* (PC+\L).
\end{equation}

Before presenting the main result, we need to recall the definition
of transposition solution $(P,\L)$ of \eqref{5.5-eq6}. To this end,
let us consider the following two (forward) SEEs:
\begin{equation}\label{op-fsystem1}
\left\{
\begin{array}{ll}
\ds dx_1 = \big[(A+A_1) x_1 + u_1\big]d\tau + \big(C x_1 + v_1\big)dW(\tau) &\mbox{ in } (t,T],\\
\ns\ds x_1(t)=\xi_1
\end{array}
\right.
\end{equation}
and
\begin{equation}\label{op-fsystem2}
\left\{
\begin{array}{ll}
\ds dx_2 =\big[ (A+A_1) x_2  + u_2\big]d\tau +\big( C x_2 + v_2 \big)dW(\tau) &\mbox{ in } (t,T],\\
\ns\ds x_2(t)=\xi_2.
\end{array}
\right.
\end{equation}
Here $t\in [0,T)$, $\xi_1,\xi_2$ are suitable
random variables and $u_1,u_2,v_1,v_2$ are
suitable stochastic processes.

 Put
$$
\begin{array}{ll} \ds C_{\dbF,w}\big([0,T];
L^{\infty}(\Om;\cL(H))\big)\\
\ns \ds\=\!\Big\{P\!:\![0,T]\!\times\!\Om\to\!
\cL(H) \,\Big|\, |P|_{\cL(H)}\!\in\!
L^\infty_\dbF(0,T),\,\chi_{[t,T]}P(\cd)\zeta\!\in\!
C_\dbF([t,T];L^2(\Om;H)), \
\forall\,\zeta\!\in\! L^2_{\cF_t}(\Om;H)\Big\},
\end{array}
$$
$$
\begin{array}{ll} \ds C_{\dbF,w}\big([0,T];
L^{\infty}(\Om;\cS(H))\big)\\
\ns \ds\=\Big\{P\in C_{\dbF,w}\big([0,T]; L^{\infty}(\Om;\cL(H))\big)
\;\Big|\;  P(t,\om)\in \cS(H) \mbox{ for a.e. }(t,\om)\in
[0,T]\times\Om\Big\},
\end{array}
$$
and
$$
\begin{array}{ll}\ds
L^2_{\dbF,D,w}(0,T;\cS(H))
\=\Big\{\L:[0,T]\times\Om\to
\dbS_2(H;H_\l')\;\Big|\; |\L|_{\cL_2(H;H_\l')}\in L^2_\dbF(0,T),\\
\ns\ds\hspace{3.8cm}
|D^*\L|_{\cL(H;U)}\in
L^\infty_\dbF(\Om;L^2(0,T)) \Big\},
\end{array}
$$
where
\bel{self-2}
\begin{array}{ll}\ds
\dbS_2(H;H_\l') \deq\big\{F\!\in\!\cL_2(H;H_\l')\;\big|\; F|_{H_\l},
 \mbox{ the restriction of }F
\mbox{ on }H_\l,\mbox{ is  a Hilbert-Schmidt  } \\
\ns\ds \hspace{4.85cm}\mbox{  operator from }H_\l \mbox{ to }H,
\mbox{ and }(F|_{H_\l})^*=F \big\},
\end{array}
\ee
is a Hilbert space with the inner product inherited from
$\cL_2(H;H_\l')$.

These spaces are Banach spaces  with
the norms
$$
|P|_{C_{\dbF,w}\big([0,T];
L^{\infty}(\Om;\cS(H))\big)}\deq
\big||P|_{\cL(H)}\big|_{L^\infty_\dbF(0,T)},\q
\forall P\in C_{\dbF,w}\big([0,T];
L^{\infty}(\Om;\cS(H))\big)
$$
and
$$
|\L|_{L^2_{\dbF,D,w}(0,T;\cS(H))}\deq
\big||\L|_{\cL_2(H;H_\l')}\big|_{L^2_\dbF(0,T)}
+ \big||D^*\L|_{\cL(H;U)}\big|_{
L^\infty_\dbF(\Om;L^2(0,T))},\q \forall
\L\in L^2_{\dbF,D,w}(0,T;\cS(H)),
$$
respectively.

\begin{definition}\label{4.8-def2}
We call $\big(P(\cd),\L(\cd)\big)\in
C_{\dbF,w}\big([0,T];
L^{\infty}(\Om;\cS(H))\big) \times
L^2_{\dbF,D,w}(0,T;\cS(H))$ a
transposition solution to
\eqref{5.5-eq6} if the following
 conditions hold:

\ms

{\rm 1)}  $K(t,\om)\big(\equiv R(t,\om) +
D(t,\om)^*P(t,\om)D(t,\om)\big)> 0$ and its left
inverse $K(t,\om)^{-1}$ is a densely defined
closed operator for a.e. $(t,\om)\in
[0,T]\times\Om$;

\ms

{\rm 2)}  For any $t\in [0,T]$,
$\xi_1,\xi_2\in L^4_{\cF_t}(\Om;H)$,
$u_1(\cd), u_2(\cd)\in
L^4_\dbF(\Om;L^2(t,T;H))$ and
$v_1(\cd), v_2(\cd)$ $\in
L^4_\dbF(\Om;L^2(t,T;H_\l))$, with
$L\=D^* PC+B^*P+D^*\L,$ it holds that
\begin{eqnarray}\label{6.18-eq1}
&& \3n\3n\3n \mE\langle
Gx_{1}(T),x_{2}(T)\rangle_{H} \!+\mE
\!\int_t^T\!\! \big\langle Q(\tau) x_{1}(\tau),
x_{2}(\tau) \big\rangle_{H}d\tau \!- \mE
\!\int_t^T\!\! \big\langle K(\tau)^{-1} L(\tau)
x_{1}(\tau), L(\tau)x_{2}(\tau)
\big\rangle_{H}d\tau\nonumber
\\
&&\3n\3n\3n = \mE\big\langle P(t)
\xi_{1},\xi_{2} \big\rangle_{H} + \mE \int_t^T
\big\langle P(\tau)u_{1}(\tau),
x_{2}(\tau)\big\rangle_{H}d\tau + \mE \int_t^T
\big\langle P(\tau)x_{1}(\tau),
u_{2}(\tau)\big\rangle_{H}d\tau \\
&& \3n\3n\3n\q  + \mE \int_t^T\big\langle
P(\tau)C(\tau)x_{1}(\tau),
v_{2}(\tau)\big\rangle_{H}d\tau + \mE
\int_t^T \big\langle  P(\tau)v_{1}(\tau), C(\tau)x_{2}(\tau)+v_{2}(\tau)\big\rangle_{H}d\tau\nonumber\\
&& \3n\3n\3n\q + \mE \int_t^T
\big\langle v_{1}(\tau),
\L(\tau)x_2(\tau)\big\rangle_{H_\l,
H_\l'}d\tau+ \mE \int_t^T \big\langle
\L(\tau)x_1(\tau), v_{2}(\tau)
\big\rangle_{H_\l',H_\l}d\tau,\nonumber
\end{eqnarray}
where $x_1(\cd)$ and $x_2(\cd)$ solve
\eqref{op-fsystem1} and \eqref{op-fsystem2},
respectively.
\end{definition}

Before stating the main result, which reveals
the relationship between the existence of
optimal feedback operator for Problem (SLQ) and
the well-posedness of \eqref{5.5-eq6} in the
sense of transposition solution, we also need
the following technical condition:

\vspace{0.2cm}

({\bf AS4}) {\it There is a dense
    subspace $\wt
U$ of $U$ such that  $R\in
L^\infty_\dbF(0,T;\cL(\wt U))$ and $B,D\in
L^\infty_\dbF(0,T;$ $\cL(\wt U;H_\l))$, where
$H_\l$ is given in Assumption ({\bf AS3}).}

\begin{remark}
Similar to Assumption {\bf (AS3)},  for an LQ problem of a given  SPDE,
Assumption {\bf (AS4)} means that some
coefficients in the equation and the
cost functional are smooth in some
sense (See Section \ref{sec-example} for example). This is not very restrictive
for many controlled SPDEs.
\end{remark}

\vspace{0.15cm}
%

Put
$$
\begin{array}{ll}\ds
\Upsilon_2(H;U)\cap \Upsilon_2(H_\l;\wt
U) \deq \{L\in  \Upsilon_2(H;U) |
\mbox{For any }f\in
L^2_{\dbF}(\Om;L^\infty(0,T;H_\l)),\;
Lf\in L_{\dbF}^2(0,T;\wt U)\big),\\ \ns
\ds \qq\qq\qq\qq\qq\qq\qq\qq\qq
|L(\cd,\cd)|_{\cL(H_\l;\wt U)}\in
L^\infty_\dbF(\Om;L^p(0,T)) \}.
\end{array}
$$

\begin{theorem}\label{5.7-th1}
Let  {\bf (AS0)}-- {\bf (AS4)} hold. Then, Problem (SLQ) admits a
unique optimal feedback operator $\Th(\cd)\in \Upsilon_2(H;U)\cap
\Upsilon_2(H_\l;\wt U)$ if and only if the Riccati equation
\eqref{5.5-eq6} admits a unique transposition solution
$\big(P(\cd),\L(\cd)\big)\in C_{\dbF,w}\big([0,T];
L^{\infty}(\Om;\cS(H))\big) \times L^2_{\dbF,D,w}(0,T; \cS(H)) $
such that
\begin{equation}\label{5.7-eq5}
\begin{array}{ll}\ds
K(\cd)^{-1}\big[B(\cd)^* P(\cd)
+D(\cd)^* P(\cd)C(\cd) +
D(\cd)^*\L(\cd)\big]\in
\Upsilon_2(H;U)\cap \Upsilon_2(H_\l;\wt
U).
\end{array}
\end{equation}
In this case, the optimal feedback operator $\Th(\cd)$ is given by
\begin{equation}\label{5.10}
\begin{array}{ll}
\ns\ds\Th(\cd)=-K(\cd)^{-1}[B(\cd)^* P(\cd)
+D(\cd)^* P(\cd)C(\cd) + D(\cd)^*\L(\cd)].
\end{array}
\end{equation}
Furthermore,
\begin{equation}\label{Value}
\inf_{u\in
L^2_\dbF(s,T;U)}\cJ(s,\eta;u)=\frac{1}{2}\,\dbE\langle
P(s)\eta,\eta\rangle_H.
\end{equation}
\end{theorem}

\br
If $H$, $U$ are finite dimensional, then {\bf
(AS1)}, {\bf (AS3)} and {\bf (AS4)} hold true.
The transposition solution becomes equivalent
with the classical adapted solution of
\eqref{5.5-eq6}, and our conclusion is reduced
to \cite[Theorem 2.1]{Lu-Wang-Zhang-PUQR}. \er
\begin{remark}
It is more natural to require the optimal
feedback operator $\Th(\cd)\in \Upsilon_2(H;U)$
other than $\Th(\cd)\in \Upsilon_2(H;U)\cap
\Upsilon_2(H_\l;\wt U)$.  Nevertheless, we do
not know how to prove Theorem \ref{5.7-th1} with
such condition now.
\end{remark}
In the above, we require $ X\in L^2_{\dbF}(\Omega;C([0,T]; H))$ and
$\Th \in\Upsilon_2(H;U)\cap \Upsilon_2(H_\l;\wt U)$ to guarantee
that $\Th X\in L^2_{\dbF}(0,T;U)$. By relaxing the condition on $X$
and strengthening that on $\Th$ as $X\in
C_{\dbF}([0,T];L^2(\Omega;H))$ and $\Th \in \wt\Upsilon_2(H;U)\cap
\wt\Upsilon_2(H_\l;\wt U)$, where
$$ \wt\Upsilon_2(H;U)\deq \big\{L \in
\mathbf{L}(L^\infty_{\dbF}(0,T;L^2(\Om;H));L_{\dbF}^2(0,T;U)) |  |L
|_{\cL(H;U)}\in L^2_\dbF(0,T;L^\infty(\Om))\big\}
$$
and
$$
\begin{array}{ll}\ds
\wt\Upsilon_2(H;U)\cap
\wt\Upsilon_2(H_\l;\wt
U) \deq \{L\in  \wt\Upsilon_2(H;U) |
\mbox{For any }f\in
L^\infty_{\dbF}(0,T;L^2(\Om;H_\l)),\;
Lf\in L_{\dbF}^2(0,T;\wt U)\big),\\ \ns
\ds \qq\qq\qq\qq\qq\qq\qq\qq\qq
|L(\cd,\cd)|_{\cL(H_\l;\wt U)}\in
L^2_\dbF(0,T;L^\infty(\Om)) \}.
\end{array}
$$
one can also obtain  $\Th X\in
L^2_{\dbF}(0,T;U)$. In this case, we can relax
({\bf AS0}) as $A$ generates a $C_0$-semigroup
on $H$. We state the analogue result of Theorem
\ref{5.7-th1} and omit the proof,  which is
similar to the one for Theorem \ref{5.7-th1}.

\begin{theorem}\label{5.7-th1-2}
Let  {\bf (AS1)}-- {\bf (AS4)} hold and $A$
generates a $C_0$-semigroup on $H$. Then,
Problem (SLQ) admits a unique optimal feedback
operator $\Th(\cd)\in \wt\Upsilon_2(H,U)\cap
\wt\Upsilon_2(H_\l;\wt U)$ if and only if the
Riccati equation \eqref{5.5-eq6} has a unique
transposition solution
$\big(P(\cd),\L(\cd)\big)\in C_{\dbF,w}([0,T];
L^{\infty}(\Om;\cS(H))) \times
L^2_{\dbF,D,w}(0,T; \cS(H))$ such that
$$\ba{ll}
\ns\ds K(\cd)^{-1}\big[B(\cd)^* P(\cd)
+D(\cd)^* P(\cd)C(\cd) +
D(\cd)^*\L(\cd)\big]\in
\wt\Upsilon_2(H;U)\cap
\wt\Upsilon_2(H_\l;\wt
U).
\ea
$$
In this case, the optimal feedback operator $\Th(\cd)$ is given by
$$\ba{ll}
\ns\ds \Th(\cd)=-K(\cd)^{-1}[B(\cd)^* P(\cd)
+D(\cd)^* P(\cd)C(\cd) + D(\cd)^*\L(\cd)].
\ea
$$
Furthermore,\vspace{-4mm}
$$\ba{ll}
\ns\ds
\inf_{u\in
L^2_\dbF(s,T;U)}\cJ(s,\eta;u)=\frac{1}{2}\,\dbE\langle
P(s)\eta,\eta\rangle_H.
\ea
$$
\end{theorem}
\begin{remark}
At a first glance, it seems that Theorem \ref{5.7-th1-2} is more general than Theorem \ref{5.7-th1}. However, we believe that Theorem \ref{5.7-th1} is more interesting due to the following reasons: first, as we have explained in Remark \ref{3.7-rmk1-1}, {\bf (AS0)} is always fulfilled by controlled SPDEs; second, the space $\Upsilon_2(H;U)\cap \Upsilon_2(H_\l;\wt
U)$ is larger than the space $\wt\Upsilon_2(H;U)\cap
\wt\Upsilon_2(H_\l;\wt
U)$.
\end{remark}

\section{Proof of the main result}

This section is devoted to the proof of Theorem \ref{5.7-th1}.
Inspired by \cite[Theorem 3.3]{Li-Sun-Yong-2016-PUQR},  we transform
the nonlinear BSRE \eqref{5.5-eq6} into a proper linear BSEE with
the solution satisfying an appropriate equality condition. To this
end, we need two preliminary results, i.e., Lemma \ref{3.14-lm1} and
Lemma \ref{th1} which are obtained in the following two subsections
respectively.

\subsection{A new transposition solution for backward stochastic Lyapunov equations in infinite dimensions}

Consider the following BSEE:
\begin{equation}\label{Relaxed-t-solution-Riccati}
\left\{\begin{array}{ll} \ds
dP=-\big[P(A+A_\Th)+(A+A_\Th)^{*}P+C_\Th^{*}PC_\Th+C_\Th^{*}\L
+\L C_\Th +Q + \Th^{*}R\Th\big]ds\\
\ns\ds\qq\qq +\L dW(s) \q\mbox{ in }  [0,T],\\
\ns\ds P(T)=G,
\end{array}\right.
\end{equation}
where $ A_\Th\=A_1+B\Th$ and $C_\Th\=C+D\Th $ for some
$\Th\in\Upsilon_2(H;U)\cap \Upsilon_2(H_\l;\wt U)$ Since $A_\Th$ and
$C_\Th$  may not belong to $L^\infty_\dbF(0,T;\cL(H))$, we cannot
employ the existing result to obtain the well-posedness of
\eqref{Relaxed-t-solution-Riccati} (see \cite{Lu-Zhang-2018-MCRF}
for example).

Inspired by \cite{Lu-Zhang-2018-MCRF}, we first define the
$H_\l$-transposition solution of (\ref{Relaxed-t-solution-Riccati}).
To this end, we consider the following two (forward) SEEs:
\begin{equation}\label{op-fsystem3}
\left\{
\begin{array}{ll}
\ds d\tilde x_1 = \big[(A+A_\Th) \tilde x_1 + \tilde u_1\big]dr + \big(C_\Th \tilde x_1 + \tilde v_1\big)dW(r) &\mbox{ in } (t,T],\\
\ns\ds \tilde x_1(t)=\tilde \xi_1
\end{array}
\right.
\end{equation}
and
\begin{equation}\label{op-fsystem4}
\left\{
\begin{array}{ll}
\ds d\tilde x_2 =\big[ (A+A_\Th) \tilde x_2  + \tilde u_2\big]dr +\big( C_\Th\tilde x_2 + \tilde v_2 \big)dW(r) &\mbox{ in } (t,T],\\
\ns\ds \tilde x_2(t)=\tilde \xi_2,
\end{array}
\right.
\end{equation}
where  $t\in [0,T)$.

Let $p>1$. If $\tilde \xi_1,\tilde \xi_2\in L^p_{\cF_t}(\Om;H_\l)$
and $\tilde u_1(\cd), \tilde u_2(\cd),\tilde v_1(\cd), \tilde
v_2(\cd)\in L^p_\dbF(\Om;L^2(t,T;H_\l))$, by \cite[Theorem
3.20]{Lu-Zhang-book1}, the equation \eqref{op-fsystem3} (\resp
\eqref{op-fsystem4}) admits a unique solution $\tilde x_{1}=\tilde
x_{1}(\cd;t,\tilde \xi_1,\tilde u_1,\tilde v_1)$ (\resp $\tilde
x_{2}=\tilde x_{2}(\cd;t,\tilde \xi_2,\tilde u_2,\tilde v_2)$) in
$L^p_\dbF(\Om;C([t,T];H_\l))\subset L^p_\dbF(\Om;C([t,T];H))$.
Further, for $j=1,2$, we have\vspace{-3mm}
\begin{equation}\label{3.14-eq32}
|\tilde x_{j}|_{L^p_\dbF(\Om;C([t,T];H_\l))}\leq \cC \big(|\tilde
\xi_j|_{L^p_{\cF_t}(\Om;H_\l)}+|\tilde
u_j|_{L^p_\dbF(\Om;L^2(t,T;H_\l))}+|\tilde
v_j|_{L^p_\dbF(\Om;L^2(t,T;H_\l))}\big).
\end{equation}
Here and in what follows, we denote by $\cC$ a
generic positive constant, which may be
different from one place to another.

If $\tilde \xi_1,\tilde \xi_2\in L^p_{\cF_t}(\Om;H)$, $\tilde
u_1(\cd), \tilde u_2(\cd)\in L^p_\dbF(\Om;L^2(t,T;H))$ and $\tilde
v_1(\cd), \tilde v_2(\cd)$ $\in L^p_\dbF(\Om;L^2(t,T;$ $H))$, by
\cite[Theorem 3.20]{Lu-Zhang-book1} again, the equation
\eqref{op-fsystem3} (\resp \eqref{op-fsystem4}) admits a unique
solution $\tilde x_{1}$ (\resp $\tilde x_{2}$) in
$L^p_\dbF(\Om;C([t,T];H))$ and for $j=1,2$, we have\vspace{-1mm}
\begin{equation}\label{3.14-eq33}
|\tilde x_{j}|_{L^p_\dbF(\Om;C([t,T];H))}\leq \cC \big(|\tilde
\xi_j|_{L^p_{\cF_t}(\Om;H)}+|\tilde
u_j|_{L^p_\dbF(\Om;L^2(t,T;H))}+|\tilde
v_j|_{L^p_\dbF(\Om;L^2(t,T;H))}\big).
\end{equation}

\begin{definition}
We call $(P(\cd),\L(\cd))\in
C_{\dbF,w}\big([0,T];L^\infty(\Omega;\cL(H))\big)\times
L^2_{\dbF}\big(0,T;L^2(\Omega;\dbS_2(H;H_\l'))\big)$
an $H_\l$-transposition solution to the equation
\eqref{Relaxed-t-solution-Riccati} if for any $t\in[0,T]$,   $\tilde
\xi_1,\tilde \xi_2\in L^4_{\cF_t}(\Om;H)$, $\tilde u_1(\cd),$
$\tilde u_2(\cd)\in L^4_\dbF(\Om;L^2(t,T;H))$ and $\tilde v_1(\cd),
\tilde v_2(\cd)$ $\in L^4_\dbF(\Om;L^2(t,T;H_\l))$, it holds that
\begin{equation}\label{3.4-eq31}
\begin{array}{ll}
\ds \dbE\lan G\tilde x_1(T),\tilde x_2(T)\ran_{H}+ \dbE\int_t^T \lan [Q(r)+\Th(r)^*R(r)\Th(r)]\tilde x_1(r),\tilde x_2(r)\ran_H dr\\
\ns\ds =\dbE\lan P(t)\tilde \xi_1,\tilde \xi_2\ran_{H}+\dbE\int_t^T \lan P(r)\tilde u_1(r),\tilde x_2(r)\ran_{H} dr+\dbE\int_t^T \lan P(r)\tilde x_1(r),\tilde u_2(r)\ran_H dr\\
\ns\ds \q +\dbE\int_t^T \lan P(r)C_\Th(r)\tilde x_1(r),\tilde
v_2(r)\ran _{H} dr+\dbE\int_t^T
\lan P(r)\tilde v_1(r), C_\Th(r)\tilde x_2(r)+\tilde v_2(r)\ran_H dr\\
\ns\ds\q +\dbE\int_t^T \lan \tilde v_1(r),\L^*(r) \tilde
x_2(r)\ran_{H_\l,H_\l'} dr +\dbE\int_t^T \lan \tilde
v_2(r),\L(r)\tilde x_1(r)\ran_{H_\l,H_\l'} dr.
\end{array}
\end{equation}
Here $\tilde x_1(\cd)$ and $\tilde x_2(\cd)$
solve \eqref{op-fsystem3} and
\eqref{op-fsystem4}, respectively.
\end{definition}
\begin{lemma}\label{3.14-lm1}
The equation \eqref{Relaxed-t-solution-Riccati} admits a unique
$H_\l$-transposition solution $(P(\cd),\L(\cd))$. Moreover,
$$
\begin{array}{ll}\ds
|(P(\cd),\L(\cd))|_{C_{\dbF,w}([0,T];L^\infty(\Omega; \cS(H) ))\times
L^2_{\dbF}(0,T;L^2(\Omega; \dbS_2(H;H_\l') ))}\\
\ns\ds \leq \cC\big(\big||Q +
\Th^{*}R\Th|_{\cL(H)}\big|_{L^\infty_\dbF(\Omega;L^1(0,T))} +
|G|_{L^\infty_{\cF_T}(\Om;\cL(H))}\big).
\end{array}
$$
\end{lemma}

Lemma \ref{3.14-lm1} is a slight modification of
\cite[Theorem 1.2]{Lu-Zhang-2018-MCRF}. Here we
relax the assumption on the measurability of the
nonhomogeneous terms and coefficients, due to
the appearance of the feedback operator $\Th$ in
the nonhomogeneous terms and coefficients of the
equation \eqref{Relaxed-t-solution-Riccati}. The
proof of Lemma \ref{3.14-lm1} is similar as the
one for  \cite[Theorem 1.2]{Lu-Zhang-2018-MCRF}.
Hence, we provide details for the different part
and give a sketch of the similar part.

We first present the following preliminary
results.
\begin{lemma}\label{lm14}\cite[Lemma 3.5]{Lu-Zhang-arxiv-2019}
If  $\{S(t)\}_{t\geq 0}$  is  a
contraction semigroup on $H$,
then it is also a  c
 semigroup on $H_\l$, and it can
be uniquely extended to be  a
contraction semigroup on $H_\l'$.
\end{lemma}
\begin{lemma}\label{lemma5}
For each $t\in[0,T]$, if $u_2=v_2=0$ in the equation
\eqref{op-fsystem4}, then there exists an operator $\Phi(\cd,t)\in
\cL\big(L^{4}_{\cF_t}(\Om;H),$ $ L^{4}_\dbF(\Om;C([t,T];H))\big)$
such that the solution to \eqref{op-fsystem4} can be represented as
$\tilde x_2(\cd) = \Phi(\cd,t)\tilde \xi_2$. Further, for any $t\in
[0,T)$, $\xi\in L^{4}_{\cF_t}(\Om;H)$ and $\e>0$, there is a $\d\in
(0,T-t)$ such that for any $s\in [t, t+\d]$, it holds that
\begin{equation}\label{10.31eq1}
|\Phi(\cd,t)\xi-\Phi(\cd,s)\xi|_{L^{4}_\dbF(\Om;C([s,T];H))} <\e.
\end{equation}
\end{lemma}
The proof of Lemma \ref{lemma5} is very similar to that of
\cite[Lemma 2.6]{Lu-Zhang-book}. Hence we omit it here.

\begin{lemma}\label{lm11-1}
The  set
$$
\begin{array}{ll}\ds
\big\{\tilde x_2(\cd)\;\big|\; \tilde x_2(\cd)\mbox{ solves }
\eqref{op-fsystem4} \mbox{ with }t=0,\;\tilde \xi_2=0,\; \tilde
v_2=0 \mbox{ and }\tilde u_2\in L^4_{\dbF}(\Om;L^2(0,T;H)) \big\}
\end{array}
$$
is dense in $L^2_{\dbF}(0,T;H)$.
\end{lemma}
The proof of Lemma \ref{lm11-1} is almost the same as the one for
\cite[Lemma 3.9]{Lu-Zhang-arxiv-2019}. We omit it.

\begin{lemma}\label{lemma2.1}\cite[Lemma 2.5]{Lu-Zhang-book}
Assume that $p\in(1,\infty]$, $q=\left\{\ba{ll}\frac{p}{p-1}&\hb{if
}\ p\in (1,\infty),\\[2mm] 1&\hb{if }\ p=\infty,\ea\right.$ $f_1\in
L^p_{\dbF}(0,T;$ $L^2(\Om;H))$ and $f_2\in
L^q_{\dbF}(0,T;L^2(\Om;H))$. Then there exists a monotonic sequence
$\{h_n\}_{n=1}^\infty$ of positive numbers such that
$\ds\lim_{n\to\infty}h_n=0$, and
\begin{equation}\label{2.21}
\lim_{n\to\infty}\frac{1}{h_n}\int_t^{t+h_n}\mathbb{E} \langle
f_1(t),f_2(\tau)\rangle_{H} d\tau=\mathbb{E} \langle
f_1(t),f_2(t)\rangle_{H},\qq\ae\, t\in
 [0,T].
 \end{equation}
\end{lemma}
\begin{lemma}\label{3.14-lm2}
Let $L\in \cL(H)$. Then $L\in \cL_2(H;H_\l')$.
\end{lemma}
{\it Proof}. For $j\in\dbN$, let $\ell_{jk}=\lan
L e_j, e_k\ran_H$ for $k\in\dbN$. Then $L e_j =
\ds\sum_{k=1}^\infty \ell_{jk}e_j$
and\vspace{-4mm}
$$
\sum_{k=1}^\infty|\ell_{jk}|^2 = |L e_j|_H^2
\leq |L|_{\cL(H)}^2.
$$
Noting that $\ell_{jk}=\lan L e_j,
e_k\ran_H=\lan  e_j, L^*e_k\ran_H$, we have
$L^*e_k = \ds\sum_{k=1}^\infty \ell_{jk}e_j$
and\vspace{-4mm}
\begin{equation}\label{3.14-eq15}
\sum_{j=1}^\infty|\ell_{jk}|^2 = |L^*e_k|_H^2
\leq |L^*|_{\cL(H)}^2=|L|_{\cL(H)}^2.
\end{equation}
By \eqref{3.14-eq15}, we have that
$$
|L|_{\cL_2(H;H_\l')}^2= \sum_{j=1}^\infty|L
e_j|_{H_\l'}^2 = \sum_{j=1}^\infty
\sum_{k=1}^\infty\l_k^2
|e_k|_{D(A)}^{-2}|\ell_{jk}|^2 \leq
\sum_{k=1}^\infty\l_k^2
\sum_{j=1}^\infty|\ell_{jk}|^2\leq
\(\sum_{k=1}^\infty\l_k^2\)|L|_{\cL(H)}.
$$
\endpf

Now we are in a position to prove Lemma \ref{3.14-lm1}.

\it Proof of Lemma \ref{3.14-lm1}. \rm
We divide the proof into six steps.

{\bf Step 1}. By Lemma \ref{3.4-lm1},
there is $\{\Th_n\}_{n=1}^\infty\subset
L^\infty_\dbF(\Om;L^2(0,T;\cL(H;U)))$
such that
\begin{equation}\label{3.4-eq43}
\lim_{n\to\infty} \Th_n v = \Th v \q
\mbox{ in }L_{\dbF}^2(0,T;U),\q
\forall\, v\in
L^2_{\dbF}(\Om;L^\infty(0,T;H)).
\end{equation}

Define a family of operators $\{\cT
(t)\}_{t\geq 0}$ on $\cL_2(H; H_\l')$
as follows:
$$
\cT (t)O = e^{At}Oe^{A^*t}, \q\forall\; O\in
\cL_2(H;H_\l').
$$

Similar to Step 1 in the proof of
\cite[Theorem 1.2]{Lu-Zhang-2018-MCRF},
we can show that $\{\cT (t)\}_{t\geq
0}$ is a  contraction semigroup
on $\cL_2(H;H_\l')$, and its
restriction on $\cL_2(H_\l;H)$ is a
contraction semigroup on
$\cL_2(H_\l;H)$.

Denote by $\cA$ the infinitesimal
generater of $\{\cT (t)\}_{t\geq 0}$.
Consider the following $\cL_2(H;H_\l')$-valued
BSEE:
\begin{equation}\label{5.19-eq1}
\left\{
\begin{array}{ll}\ds
dP_n = -\cA^* P_ndt + f_n(t,P_n,\L_n)dt
+ \L_ndW(t) &\mbox{ in
}[0,T),\\
\ns\ds P(T)=G,
\end{array}
\right.
\end{equation}
where
\begin{equation}\label{s5f}
f_n(t,P_n,\L_n) = -(P_n A_{\Th_n} + A_{\Th_n}^{*}P_n+C_{\Th_n}^{*}
P_nC_{\Th_n} +C_{\Th_n}^{*}\L_n +\L_n C_{\Th_n}  +Q +
\Th_n^{*}R\Th_n)
\end{equation}
with $ A_{\Th_n}\=A_1+B\Th_n$ and $C_{\Th_n}\=C+D\Th_n$.

Since $G\in L^2_{\cF_T}(\Om;\cL(H))$, by Lemma \ref{3.14-lm2}, we
get that
\begin{equation}\label{3.14-eq19}
|G|_{L^2_{\cF_T}(\Om;\cL_2(H;H_\l'))} =
\big||G|_{\cL_2(H;H_\l')}\big|_{L^2_{\cF_T}(\Om)} \leq
\big|\cC|G|_{\cL(H)}\big|_{L^2_{\cF_T}(\Om)} =
\cC\big||G|_{\cL(H)}\big|_{L^2_{\cF_T}(\Om)}.
\end{equation}
Similarly, we have that\vspace{-4mm}
\begin{equation}\label{3.14-eq18}
|Q|_{L^\infty_\dbF(\Om;L^2(0,T; \cL_2(H;H_\l')))} \leq \cC
|Q|_{L^\infty_\dbF(\Om;L^2(0,T; \cL(H)))}
\end{equation}
and
\begin{equation}\label{3.14-eq16}
|\Th_n^{*}R\Th_n|_{L^\infty_\dbF(\Om;L^1(0,T;
\cL_2(H;H_\l')))}\leq \cC
\big||\Th_n^{*}R\Th_n|_{\cL(H)}\big|_{L^\infty_\dbF(\Om;L^2(0,T
))}.
\end{equation}
Since $|\Th_n|_{\cL(H;U)}\leq |\Th|_{\cL(H;U)}$, we get from
\eqref{3.14-eq16} that
\begin{equation}\label{3.14-eq17}
\begin{array}{ll}\ds
|\Th_n^{*}R\Th_n|_{L^\infty_\dbF(\Om;L^1(0,T;
\cL_2(H;H_\l')))}\\
\ns\ds \leq \cC  \big||\Th_n^{*}|_{\cL( U,H)}|R|_{\cL(
U)}|\Th_n|_{\cL(H, U)}\big|_{L^\infty_\dbF(\Om;L^1(0,T
))} \\
\ns\ds \leq \cC \big||\Th^{*}|_{\cL(
U,H)}\big|_{L^\infty_\dbF(\Om;L^2(0,T ))}|R|_{L^\infty_\dbF(\cL(
U))}\big||\Th|_{\cL(H, U)}\big|_{L^\infty_\dbF(\Om;L^2(0,T ))}.
\end{array}
\end{equation}

For any $\eta_j\in \cL_2(H;H_\l')$
($j=1,2,3,4$), we have
$$
\begin{array}{ll}\ds
\big|f(t,\eta_1,\eta_2)-
f(t,\eta_3,\eta_4)\big|_{\cL_2(H;H_\l')}\\
\ns\ds \leq  \big(|A_1|_{\cL(H)} +|A_1^*|_{\cL(H_\l')}
+|C|_{\cL(H)}|C^*|_{\cL(H_\l')} + |B|_{\cL(U;H)} |\Th_n|_{\cL(H;U)}+
|B^*|_{\cL(H_\l';\wt U')} |\Th_n^*|_{\cL(\wt U';H_\l')}
\\
\ns\ds\qq +
|C^*|_{\cL(H_\l')}|D|_{\cL(U;H)}
|\Th_n|_{\cL(H;U)} +
|C|_{\cL(H)}|D^*|_{\cL(H_\l';\wt U')}
|\Th_n^*|_{\cL(\wt U';H_\l')}\\
\ns\ds\qq +|D|_{\cL(U;H)}
|\Th_n|_{\cL(H;U)} |D^*|_{\cL(H_\l';\wt
U')} |\Th_n^*|_{\cL(\wt U';H_\l')}
\big)\big|\eta_1-\eta_2\big|_{\cL_2(H;H_\l')} \\
\ns\ds\q + \big(|C|_{\cL(H)} +|C^*|_{\cL(H_\l')}+|D|_{\cL(U;H)}
|\Th_n|_{\cL(H;U)}+|D^*|_{\cL(H_\l';\wt U')} |\Th_n^*|_{\cL(\wt
U';H_\l')} \big)\big|\eta_3-\eta_4\big|_{\cL_2(H;H_\l')}.
\end{array}
$$
From ({\bf AS2}), ({\bf AS3}), \eqref{10.10-eq30} and noting that
\begin{equation}\label{3.14-eq22}
|\Th_n|_{\cL(H;U)}\leq
|\Th|_{\cL(H;U)},\q |\Th_n^*|_{\cL(\wt
U';H_\l')}\leq |\Th^*|_{\cL(\wt
U';H_\l')},
\end{equation}
we see that there exist $L_1\in L^1(0,T)$ (independent of $n$) and
$L_2\in L^2(0,T)$ such that
\begin{equation}\label{3.14-eq12}
\begin{array}{ll}\ds
\big|f(t,\eta_1,\eta_2)-
f(t,\eta_3,\eta_4)\big|_{\cL_2(H;H_\l')}
\leq
L_1(t)\big|\eta_1-\eta_2\big|_{\cL_2(H;H_\l')}+L_2(t)
\big|\eta_3-\eta_4\big|_{\cL_2(H;H_\l')},\q
\dbP\mbox{-a.s.}
\end{array}
\end{equation}

Noting \eqref{3.14-eq19}, \eqref{3.14-eq18}, \eqref{3.14-eq17} and
\eqref{3.14-eq12}, by \cite[Theorem 4.10]{Lu-Zhang-book1}, there
exists a unique mild solution $(P_n,\L_n)\in
L^2_{\dbF}(\Omega;C([0,T];\cL_2(H;H_\l')))\times
L^2_{\dbF}(0,T;\cL_2(H;H_\l'))$ of \eqref{5.19-eq1} such that
$$
\begin{array}{ll}\ds
|(P_n,Q_n)|_{ L^2_{\dbF}(\Omega;C([0,T];\cL_2(H;H_\l')))\times
L^2_{\dbF}(0,T;\cL_2(H;H_\l'))}\\
\ns\ds \leq \cC(L_1,L_2)
\big(|G|_{L^2_{\cF_T}(\Om;\cL_2(H;H_\l'))}
+ |Q +
\Th_n^{*}R\Th_n|_{L^2_{\dbF}(\Om;L^1(0,T;\cL_2(H;H_\l')))}
\big).
\end{array}
$$
Hence,
\begin{equation}\label{3.14-eq26}
\begin{array}{ll}\ds
|(P_n,Q_n)|_{
L^2_{\dbF}(\Omega;C([0,T];\cL_2(H;H_\l')))\times
L^2_{\dbF}(0,T;\cL_2(H;H_\l'))} \\
\ns\ds \leq\cC
\big(|G|_{L^2_{\cF_T}(\Om;\cL(H))} +
|Q|_{L^2_{\dbF}(\Om;L^1(0,T;\cL(H)))} +
\big||\Th_n^{*}R\Th_n|_{\cL(H)}\big|_{L^2_{\dbF}(\Om;L^1(0,T))}
\big)\\
\ns\ds \leq\cC
\big(|G|_{L^2_{\cF_T}(\Om;\cL(H))} +
|Q|_{L^2_{\dbF}(0,T;\cL(H))} +
\big||\Th^{*}_n|_{\cL(U;H)}|R|_{\cL(U)}|\Th_n|_{\cL(H;U)}\big|_{L^2_{\dbF}(\Om;L^1(0,T))}
\big)\\
\ns\ds \leq\cC
\big(|G|_{L^2_{\cF_T}(\Om;\cL(H))} +
|Q|_{L^2_{\dbF}(0,T;\cL(H))} +
\big||\Th^{*}|_{\cL(U;H)}|R|_{\cL(U)}|\Th|_{\cL(H;U)}\big|_{L^2_{\dbF}(\Om;L^1(0,T))}
\big),
\end{array}
\end{equation}
where the constant $\cC$ is independent
of $n\in\dbN$.

\ms

{\bf Step 2}. For $t\in [0,T)$, consider the following two (forward)
SEEs:
\begin{equation}\label{op-fsystem3-n}
\left\{
\begin{array}{ll}
\ds d\tilde x_{1,n} = \big[(A+A_{\Th_n}) \tilde x_{1,n} + \tilde u_1\big]d\tau + \big(C_{\Th_n} \tilde x_{1,n} + \tilde v_1\big)dW(\tau) &\mbox{ in } (t,T],\\
\ns\ds \tilde x_{1,n}(t)=\tilde \xi_1
\end{array}
\right.
\end{equation}
and
\begin{equation}\label{op-fsystem4-n}
\left\{
\begin{array}{ll}
\ds d\tilde x_{2,n} =\big[ (A+A_{\Th_n}) \tilde x_{2,n}  + \tilde u_2\big]d\tau +\big( C_{\Th_n}\tilde x_{2,n} + \tilde v_2 \big)dW(\tau) &\mbox{ in } (t,T],\\
\ns\ds \tilde x_{2,n}(t)=\tilde \xi_2.
\end{array}
\right.
\end{equation}
From \eqref{3.14-eq22}, we know there exist  $L_3\in L^1(0,T)$ and
$L_4\in L^2(0,T)$ such that for all $n\in\dbN$,
$$
|A_{\Th_n}|_{\cL(H)} + |A_{\Th_n}|_{\cL(H_\l)}\leq L_3(t), \q
\dbP\mbox{-a.s.,}
$$
and
$$
|C_{\Th_n}|_{\cL(H)} + |C_{\Th_n}|_{\cL(H_\l)}\leq L_4(t), \q
\dbP\mbox{-a.s.}
$$
Let $p>1$. If $\tilde \xi_1,\tilde \xi_2\in L^p_{\cF_t}(\Om;H_\l)$
and $\tilde u_1(\cd), \tilde u_2(\cd),\tilde v_1(\cd), \tilde
v_2(\cd)\in L^p_\dbF(\Om;L^2(t,T;H_\l))$, by \cite[Theorem
3.20]{Lu-Zhang-book1}, the equation \eqref{op-fsystem3-n} (\resp
\eqref{op-fsystem4-n}) admits a unique solution $\tilde
x_{1,n}=\tilde x_{1,n}(\cd;t,\tilde \xi_1,\tilde u_1,\tilde v_1)$
(\resp $\tilde x_{2,n}=\tilde x_{2,n}(\cd;t,\tilde \xi_2,\tilde
u_2,\tilde v_2)$) in $L^p_\dbF(\Om;C([t,T];H_\l))$ and for $j=1,2$,
\begin{equation}\label{3.14-eq20}
|\tilde x_{j,n}|_{L^p_\dbF(\Om;C([t,T];H_\l))}\leq \cC(L_3,
L_4)\big(|\tilde \xi_j|_{L^p_{\cF_t}(\Om;H_\l)}+|\tilde
u_j|_{L^p_\dbF(\Om;L^2(t,T;H_\l))}+|\tilde
v_j|_{L^p_\dbF(\Om;L^2(t,T;H_\l))}\big),
\end{equation}
where the constant $\cC(L_3,L_4)$ is
independent of $n\in\dbN$.

If $\tilde \xi_1,\tilde \xi_2\in L^p_{\cF_t}(\Om;H)$, $\tilde
u_1(\cd), \tilde u_2(\cd)\in L^p_\dbF(\Om;L^2(t,T;H))$ and $\tilde
v_1(\cd), \tilde v_2(\cd)$ $\in L^p_\dbF(\Om;L^2(t,T; H))$, by
\cite[Theorem 3.20]{Lu-Zhang-book1} again, the equation
\eqref{op-fsystem3-n} (\resp \eqref{op-fsystem4-n}) admits a unique
solution $\tilde x_{1,n}$ (\resp $\tilde x_{2,n}$) in
$L^p_\dbF(\Om;C([t,T];H))$ and for $j=1,2$,
\begin{equation}\label{3.14-eq21}
|\tilde
x_{j,n}|_{L^p_\dbF(\Om;C([t,T];H))}\leq
\cC(L_3, L_4)\big(|\tilde
\xi_j|_{L^p_{\cF_t}(\Om;H)}+|\tilde
u_j|_{L^p_\dbF(\Om;L^2(t,T;H))}+|\tilde
v_j|_{L^p_\dbF(\Om;L^2(t,T;H))}\big),
\end{equation}
where the constant $\cC(L_3,L_4)$ is
independent of $n\in\dbN$.

\ss

We claim that   for $j=1,2$ and $p=4$,
\begin{equation}\label{3.14-eq9}
\lim_{n\to\infty}|\tilde x_{j,n}-\tilde
x_j|_{L^4_\dbF(\Om;C([t,T];H))}=0.
\end{equation}

Denote by $\lceil a \rceil$ the integer
part of a number $a\in\dbR$. Write
$$
\begin{array}{ll}\ds
N=\Big\lceil\frac{1}{\e}\big(\big||A_1|_{\cL(H)}\big|_{L^\infty_\dbF(\Om;L^1(t,T))}
+\big||B|_{\cL(U;H)}|\Th|_{\cL(H;U)}\big|_{L^\infty_\dbF(\Om;L^1(t,T))}\big)^4\\
\ns\ds \qq +\big(
\big||C|_{\cL(H)}\big|_{L^\infty_\dbF(\Om;L^2(s,T))}
+\big||D|_{\cL(U;H)}|\Th
|_{\cL(H;U)}\big|_{L^\infty_\dbF(\Om;L^2(t,T))} \big)^4\Big\rceil
+1,
\end{array}
$$
where $\e>0$ is a constant to be
determined later. Define a sequence of
stopping times
$\{\tau_{j,\e}\}_{j=1}^N$ as follows:
\begin{equation*}\label{4.6-eq10}
\!\!\left\{\!\!\ba{ll}
\ds\tau_{1,\e}(\om)\!=\!\inf\left\{\tau\in
[t,T]\,\left|\,\[\int_t^\tau\big(|A_1(r,\om)|_{\cL(H)}+ |B(r,\om)|_{\cL(U;H)}|\Th(r,\om)|_{\cL(H;U)}\big)dr\]^4\right.\right.\\
\ns\ds\hspace{4cm}\left.+\[\int_t^r\big(|C(r,\om)|^2_{\cL(H)} + |D(r,\om)|_{\cL(U;H)}^2|\Th(r,\om)|_{\cL(H;U)}^2\big)dr\]^2=\e\right\},\\[3mm]
\ns\ds \tau_{k,\e}(\om)\!=\! \inf\left\{\tau\in
[\tau_{k-1,\e}(\om),T]\,\left|\,\[\int_{\tau_{k-1,\e}}^\tau\!\!\!\big(|A_1(r,\om)|_{\cL(H)}+ |B(r,\om)|_{\cL(U;H)}|\Th(r,\om)|_{\cL(H;U)}\big)dr\]^4\right.\right.\\
\ns\ds\hspace{4.5cm}\left.+\!\int_{\tau_{k-1,\e}}^\tau\!\!\big(|C(r,\om)|^2_{\cL(H)}
+
|D(r,\om)|_{\cL(U;H)}^2|\Th(r,\om)|_{\cL(H;U)}^2\big)dr\]^2\!=\e\right\},\\
\ns\ds \hspace{8cm} k= 2,\cdots,N.
\ea\right.
\end{equation*}
Here, we agree that $\inf\emptyset=T$.

For any $s\in [t,\tau_{1,\e}]$, by
Burkholder-Davis-Gundy's inequality, we
get that
\begin{eqnarray}\label{3.14-eq10}
&& \mE\sup_{\tau\in [t,\tau_{1,\e}]}\big|\tilde
x_{j,n}(\tau)-\tilde
x_j(\tau)\big|_H^4 \nonumber\\
&& = \mE\sup_{\tau\in
[t,\tau_{1,\e}]}\Big|\int_t^\tau
e^{A(\tau-r)}\big( A_ \Th \tilde x_j -
A_{\Th_n}\tilde x_{j,n}\big)dr  + \int_t^\tau
e^{A(\tau-r)}\big( C_\Th x_j
- C_{\Th_n} \tilde x_{j,n}\big)dW(r)\Big|_H^4\nonumber\\
&& \leq \cC\mE\sup_{\tau\in
[t,\tau_{1,\e}]}\[\Big|\int_t^\tau e^{A(\tau-r)}
A_{\Th_n}\big(\tilde x_j-\tilde
x_{j,n}\big)dr\Big|_H^4 + \Big|\int_t^\tau
e^{A(\tau-r)}
B(\Th-\Th_n) \tilde x_j dr\Big|_H^4 \nonumber\\
&&\qq  + \Big|\int_t^\tau\! e^{A(\tau-r)}
C_{\Th_n} \big(\tilde x_j\!-\!\tilde
x_{j,n}\big)dW(r)\Big|_H^4\! +
\Big|\int_t^\tau\! e^{A(\tau-r)}
D(\Th\!-\!\Th_n) \tilde x_j dW(r)\Big|_H^4  \] \nonumber\\
&& \leq \cC\mE\sup_{\tau\in
[t,\tau_{1,\e}]}\Big\{\[\int_t^\tau
\big(|A_1|_{\cL(H)}+|B|_{\cL(U;H)}|\Th_n|_{\cL(H;U)}\big)\big|\tilde
x_j-\tilde x_{j,n}\big|_H
dr\]^4 \nonumber\\
&&\qq  +
\[\int_t^\tau  \big(|C|_{\cL(H)}^2+|D|_{\cL(U;H)}^2|\Th_n|_{\cL(H;U)}\big)^2 \big|\tilde x_j-\tilde x_{j,n}\big|_H^2dr\]^2 \\
&&\qq +  \[ \int_t^\tau\! |B|_{\cL(U;H)}|\Th -
\Th_n|_{\cL(H;U)} |\tilde x_j|_H dr\]^4\! +\!
\[\int_t^\tau\! |D|_{\cL(U;H)}^2|\Th -
\Th_n|_{\cL(H;U)}^2
|\tilde x_j|_H^2 dr\]^2\Big\}\nonumber\\
&& \leq \cC\Big|\[\int_t^{\tau_{1,\e}}
\big(|A_1|_{\cL(H)}+|B|_{\cL(U;H)}|\Th|_{\cL(H;U)}\big)dr
\]^4 \nonumber\\
&&\qq   +
\[\int_t^{\tau_{1,\e}}  \big(|C|_{\cL(H)}^2+ |D|_{\cL(U;H)}^2|\Th|_{\cL(H;U)}^2\big)dr\]^2
\Big|_{L^\infty(\Om)}
\mE\sup_{\tau\in [t,\tau_{1,\e}]}\big|\tilde x_j(\tau)-\tilde x_{j,n}(\tau)\big|_H^4 \nonumber\\
&&\q + \cC\mE\[\( \int_t^{\tau_{1,\e}}\!
|B|_{\cL(U;H)}|(\Th\!-\!\Th_n)\tilde x_j|_{U}
dr\)^4 \!+\! \( \int_t^{\tau_{1,\e}}\!
|D|_{\cL(U;H)}^2|(\Th-\Th_n)\tilde x_j |_{U}^2
dr\)^2\], \nonumber
\end{eqnarray}
where the constant $\cC$ is independent
of $\tau_{1,\e}$.

Let us choose $\e=\frac{1}{\cC}$. We
find that
\begin{equation}\label{3.14-eq11}
\begin{array}{ll}\ds
\mE\sup_{\tau\in [t,\tau_{1,\e}]}\big|\tilde
x_{j,n}(\tau)-\tilde x_j(\tau)\big|_H^4 \\
\ns\ds \leq \cC\mE\[\( \int_t^{\tau_{1,\e}}
|B|_{\cL(U;H)}|(\Th\!-\!\Th_n)\tilde x_j |_{U}
dr\)^4 \!+\! \( \int_t^{\tau_{1,\e}}
|D|_{\cL(U;H)}^2|(\Th-\Th_n)\tilde x_j |_{U}^2
dr\)^2\].
\end{array}
\end{equation}
By \eqref{3.4-eq43} and
\eqref{3.14-eq11}, we see that
\begin{equation}\label{3.14-eq11-1}
\lim_{n\to\infty}\mE\sup_{\tau\in
[t,\tau_{1,\e}]}\big|\tilde x_{j,n}(\tau)-\tilde
x_j(\tau)\big|_H^4=0.
\end{equation}
If $\tau_{1,\e}=T$, $\dbP$-a.s., then
we get \eqref{3.14-eq9}. Otherwise, by
repeating the above argument $N$ times, we can
prove \eqref{3.14-eq9}.

\ms

\ms

{\bf Step 3}. By \eqref{3.14-eq26}, we know that
$\{(P_n,\L_n)\}_{n=1}^\infty$ is a bounded sequence in
$L^2_{\dbF}(0,T;\cL_2(H;H_\l'))\times
L^2_{\dbF}(0,T;\cL_2(H;H_\l'))$. Hence, there is a subsequence of
$\{(P_n,\L_n)\}_{n=1}^\infty$, denoted by
$\{(P_{n_k},\L_{n_k})\}_{k=1}^\infty$, such that
\begin{equation}\label{3.14-eq28-1}
\begin{array}{ll}\ds
\{(P_{n_k},\L_{n_k})\}_{n=1}^\infty \mbox{ converges weakly to some
} (P,\L)\\
\ns\ds\qq \mbox{ in }L^2_{\dbF}(0,T;\cL_2(H;H_\l'))\times
L^2_{\dbF}(0,T;\cL_2(H;H_\l'))\mbox{  as }k\to\infty.
\end{array}
\end{equation}
By \eqref{3.14-eq26} again, we know
that for any $t\in [0,T]$,
$\{P_{n_{k}}(t)\}_{n=1}^\infty$  is a
bounded sequence in
$L^2_{\cF_t}(\Om;\cL_2(H;$ $H_\l'))$.
Hence, for every $t\in [0,T)$, there is
a subsequence of
$\{P_{n_{k}}(t)\}_{n=1}^\infty$,
denoted by
$\{P_{n_{k,t}}(t)\}_{n=1}^\infty$, such
that
\begin{equation}\label{3.14-eq29}
\{P_{n_{k,t}}\}_{n=1}^\infty \mbox{ converges weakly to some }  P^t
\mbox{ in }L^2_{\cF_t}(\Om;\cL_2(H;H_\l'))\mbox{  as }k\to\infty.
\end{equation}

Similar to Steps 2-3 of the proof of
\cite[Theorem 1.2]{Lu-Zhang-2018-MCRF},
we have that
\begin{equation}\label{3.14-eq27}
\begin{array}{ll} \ds \dbE\lan
G\tilde x_{1,n}(T),\tilde
x_{2,n}(T)\ran_{H}
+ \dbE\int_t^T \lan [Q(r)+\Th_{n}(r)^*R(r)\Th_{n}(r)]\tilde x_{1,n}(r),\tilde x_{2,n}(r)\ran_H dr\\
\ns\ds =\dbE\lan P_{n}(t) \tilde \xi_1,\tilde
\xi_2\ran_{H_\l',H_\l}\!+\!\dbE\!\int_t^T\!\! \lan P_{n}(r)\tilde
u_1(r),\tilde x_{2,n}(r)\ran_{H_\l',H_\l} dr
\\
\ns\ds \q +\dbE \int_t^T \lan P_{n}(r)\tilde x_{1,n}(r),\tilde
u_2(r)\ran_{H_\l',H_\l} dr +\dbE \int_t^T \lan
P_{n}(r)C_{\Th_{n}}(r)\tilde x_{1,n}(r),\tilde v_2(r)\ran
_{H_\l',H_\l} dr\\
\ns\ds \q +\dbE \int_t^T  \lan P_{n}(r)\tilde v_1(r),
C_{\Th_{n}}(r)\tilde x_{2,n}(r) + \tilde v_2(r) \ran_{H_\l',H_\l} dr
\\
\ns\ds\q  +\dbE\int_t^T \lan \tilde v_1(r),\L_{n}^*(r) \tilde
x_{2,n}(r)\ran_{H_\l,H_\l'} dr +\dbE\int_t^T \lan \tilde
v_2(r),\L_{n}(r)\tilde x_{1,n}(r)\ran_{H_\l,H_\l'} dr.
\end{array}
\end{equation}
Hence, for any $t\in [0,T)$,
\begin{eqnarray}\label{3.14-eq27-1}
&&\3n\3n\3n\3n \dbE\lan
G\tilde x_{1,n_{k,t}}(T),\tilde
x_{2,n_{k,t}}(T)\ran_{H}
+ \dbE\int_t^T \lan [Q(r)+\Th_{n_{k,t}}(r)^*R(r)\Th_{n_{k,t}}(r)]\tilde x_{1,n_k}(r),\tilde x_{2,n_k}(r)\ran_H dr\nonumber\\
&&\3n\3n\3n\3n =\dbE\lan P_{n_{k,t}}(t) \tilde \xi_1,\tilde
\xi_2\ran_{H_\l',H_\l}\!+\!\dbE\!\int_t^T\!\! \lan
P_{n_{k,t}}(r)\tilde u_1(r),\tilde x_{2,n_{k,t}}(r)\ran_{H_\l',H_\l}
dr\nonumber
\\
&&\3n\3n  +\dbE \int_t^T \lan P_{n_k}(r)\tilde x_{1,n_k}(r),\tilde
u_2(r)\ran_{H_\l',H_\l} dr +\dbE \int_t^T \lan
P_{n_{k,t}}(r)C_{\Th_{n_{k,t}}}(r)\tilde x_{1,n_{k,t}}(r),\tilde
v_2(r)\ran
_{H_\l',H_\l} dr\nonumber\\
&&\3n\3n  +\dbE \int_t^T  \lan P_{n_{k,t}}(r)\tilde v_1(r),
C_{\Th_{n_{k,t}}}(r)\tilde x_{2,n_{k,t}}(r) + \tilde v_2(r)
\ran_{H_\l',H_\l} dr
\\
&& \3n\3n +\dbE\int_t^T \lan \tilde v_1(r),\L_{n_{k,t}}^*(r) \tilde
x_{2,n_{k,t}}(r)\ran_{H_\l,H_\l'} dr +\dbE\int_t^T \lan \tilde
v_2(r),\L_{n_{k,t}}(r)\tilde x_{1,n_{k,t}}(r)\ran_{H_\l,H_\l'} dr.\nonumber
\end{eqnarray}
Letting $k\to\infty$ in \eqref{3.14-eq27-1}, by \eqref{3.4-eq43},
\eqref{3.14-eq11-1} and \eqref{3.14-eq28-1}, we obtain that
\begin{equation}\label{3.14-eq28}
\begin{array}{ll} \ds \dbE\lan
G\tilde x_{1}(T),\tilde
x_{2}(T)\ran_{H}
+ \dbE\int_t^T \lan [Q(r)+\Th(r)^*R(r)\Th(r)]\tilde x_{1}(r),\tilde x_{2}(r)\ran_H dr\\
\ns\ds =\dbE\lan P^t \tilde \xi_1,\tilde
\xi_2\ran_{H_\l',H_\l}\!+\!\dbE\!\int_t^T\!\! \lan P(r)\tilde
u_1(r),\tilde x_{2}(r)\ran_{H_\l',H_\l} dr
\\
\ns\ds \q +\dbE \int_t^T \lan P(r)\tilde x_{1}(r),\tilde
u_2(r)\ran_{H_\l',H_\l} dr +\dbE \int_t^T \lan P(r)C_{\Th}(r)\tilde
x_{1}(r),\tilde v_2(r)\ran
_{H_\l',H_\l} dr\\
\ns\ds \q +\dbE \int_t^T  \lan P(r)\tilde v_1(r), C_{\Th}(r)\tilde
x_{2}(r) + \tilde v_2(r) \ran_{H_\l',H_\l} dr
\\
\ns\ds\q  +\dbE\int_t^T \lan \tilde v_1(r),\L^*(r) \tilde
x_{2}(r)\ran_{H_\l,H_\l'} dr +\dbE\int_t^T \lan \tilde
v_2(r),\L(r)\tilde x_{1}(r)\ran_{H_\l,H_\l'} dr.
\end{array}
\end{equation}

\ms

{\bf Step 4}. In this step, we prove
that  $P^{\cd}\in
C_{\dbF,w}([0,T];L^{\infty}(\Om;\cL(H)))$.

We first show that for any $t\in
[0,T]$, $|P^t|_{\cL(H)}\in
L^{\infty}_{\cF_t}(\Omega)$ by a
contradiction argument. Assume that this was
untrue for some $t\in [0,T]$. Noting that $H$ is separable, by \cite[Corollary 2.3]{Blasco-2010},  we
can find two sequences
$\{\eta_n\}_{n=1}^\infty,\{\hat\eta_n\}_{n=1}^\infty\subset
L^{4}_{\cF_t}(\Omega;H_\l)$ with
$|\eta_n|_{L^{2}_{\cF_t}(\Omega;H)}=|\hat\eta_n|_{L^{2}_{\cF_t}(\Omega;H)}=1$
for all $n\in\dbN$ such that
\begin{equation}\label{5.24-eq3}
\mE\langle P^t\eta_n, \hat
\eta_n\rangle_{H_\l',H_\l} \geq n\;
\mbox{ for all }n\in\dbN.
\end{equation}
Thanks to \eqref{3.14-eq28},  we get
that for any $n\in\dbN$,
\begin{equation}\label{5.25-eq3}
\!\!\begin{array}{ll}\ds \q\dbE
\big\langle
P^t\eta_n,\hat\eta_{n}\big\rangle_{H_\l',H_\l}
\\
\ns \ds = \dbE \big\langle G \tilde
x_{1}(T;t,\eta_{n},0,0),\tilde
x_{2}(T;t,\hat\eta_n,0,0)\big\rangle_{H}\\
\ns\ds\q
- \dbE \int_t^T \big\langle
\big(Q+\Th^*R\Th\big)\tilde
x_{1}(s;t,\eta_{n},0,0),\tilde
x_{2}(s;t,\hat\eta_n,0,0)\big\rangle_{H}
ds.
\end{array}
\end{equation}
By \eqref{5.25-eq3} and
\eqref{3.14-eq21}, for any $n\in\dbN$,
it holds that
\begin{eqnarray}\label{5.25-eq5}
&&\3n
 \big|\mE\langle P^t\eta_n, \hat
\eta_n\rangle_{H}\big|\nonumber\\
&&\3n\3n\3n\leq
|G|_{L^{\infty}_{\cF_T}(\Om;\cL(H))}|x_1(T;t,\eta_n,0,0)|_{L^{2}_{\cF_T}(\Omega;H)}
|x_2(T;t,\hat\eta_n,0,0)|_{L^{2}_{\cF_T}(\Omega;H)}\nonumber
\\
&& \3n +
\big||(Q\!+\!\Th^*R\Th)|_{\cL(H)}\big|_{L^{\infty}_\dbF(\Om;L^1(0,T))}
|x_1(\cd;t,\eta_n,0,0)|_{L^{2}_\dbF(\Omega;C([t,T];H))}
|x_2(\cd;t,\hat\eta_n,0,0)|_{L^{2}_\dbF(\Omega;C([t,T];H))}\nonumber
\\
&&\3n\3n\3n \leq
\cC\big(|G|_{L^{\infty}_{\cF_T}(\Om;\cL(H))}
+
\big||(Q+\Th^*R\Th)|_{\cL(H)}\big|_{L^{\infty}_\dbF(\Om;L^1(0,T))}\big)|\eta_n|_{L^{2}_{\cF_t}(\Omega;H)}
|\hat\eta_n|_{L^{2}_{\cF_t}(\Omega;H)}\\
&&\3n\3n\3n \leq
\cC\big(|G|_{L^{\infty}_{\cF_T}(\Om;\cL(H))}
+
\big||(Q+\Th^*R\Th)|_{\cL(H)}\big|_{L^{\infty}_\dbF(\Om;L^1(0,T))}\big).\nonumber
\end{eqnarray}
This contradicts \eqref{5.24-eq3}.
Hence, we get that $|P^t|_{\cL(H)}\in
L^{\infty}_{\cF_t}(\Omega)$ for all
$t\in [0,T]$. By \cite[Corollary 2.3]{Blasco-2010}, for any
$t\in [0,T]$, we can find
$\eta,\hat\eta\in L^2_{\cF_t}(\Om;H)$
such that
$|\eta|_{L^2_{\cF_t}(\Om;H)}=|\hat\eta|_{L^2_{\cF_t}(\Om;H)}=1$
and that
\begin{equation}\label{3.14-eq23}
\frac{1}{2}\big||P^t|_{\cL(H)}\big|_{L^{\infty}_{\cF_t}(\Omega)}\leq
\dbE \big\langle
P(t)\eta,\hat\eta\big\rangle_{H}.
\end{equation}
Similar to the proof of
\eqref{5.25-eq5}, we can obtain that
\begin{equation}\label{3.14-eq24}
\big|\mE\langle P^t\eta, \hat
\eta\rangle_{H}\big|\leq
\cC\big(|G|_{L^{\infty}_{\cF_T}(\Om;\cL(H))}
+
\big||(Q+\Th^*R\Th)|_{\cL(H)}\big|_{L^{\infty}_\dbF(\Om;L^1(0,T))}\big),\q
\forall t\in [0,T].
\end{equation}
From \eqref{3.14-eq23} and
\eqref{3.14-eq24}, we find that
\begin{equation}\label{3.14-eq25}
\big||P^t|_{\cL(H)}\big|_{L^{\infty}_{\cF_t}(\Omega)}
 \leq
\cC\big(|G|_{L^{\infty}_{\cF_T}(\Om;\cL(H))}
+
\big||(Q+\Th^*R\Th)|_{\cL(H)}\big|_{L^{\infty}_\dbF(\Om;L^1(0,T))}\big),\q
\forall t\in [0,T].
\end{equation}

\ms

By \eqref{3.14-eq28} again, we see that
for any $\tilde\xi_1, \tilde\xi_2\in
L^2_{\cF_t}(\Om;H)$,
\begin{equation}\label{3.14-eq30}
\begin{array}{ll}\ds
\q\dbE \big\langle
P^t\tilde\xi_1,\tilde\xi_2\big\rangle_{H}
\3n&\ds =  \dbE \big\langle G \tilde
x_{1}(T;t,\tilde\xi_1,0,0),\tilde
x_{2}(T;t,\tilde\xi_2,0,0)\big\rangle_{H}
\\
\ns& \ds \q- \dbE \int_t^T \big\langle
\big(Q+\Th^*R\Th\big)\tilde
x_{1}(s;t,\tilde\xi_1,0,0),\tilde
x_{2}(s;t,\tilde\xi_2,0,0)\big\rangle_{H}
ds.
\end{array}
\end{equation}
From \eqref{3.14-eq30} and Lemma \ref{lemma5}, we get that
$$
\mE\big\langle
P^t\tilde\xi_1,\tilde\xi_2
\big\rangle_{H} =\mE\Big\langle \Phi(T,t)^*G
\Phi(T,t)\tilde\xi_1 - \int_{t}^T
\Phi(r,t)^*\big(Q+\Th^*R\Th\big)
\Phi(r,t)\tilde\xi_1dr, \tilde\xi_2
\Big\rangle_{H}.
$$
This, together with the arbitrariness of $\tilde\xi_2\in
L^2_{\cF_t}(\Om;H)$, implies that
\begin{equation}\label{5.25-eq8}
P^t\xi_1 =
\mE_t\[\Phi(T,t)^*G \Phi(T,t)\tilde\xi_1 -
\int_{t}^T
\Phi(r,t)^*\big(Q+\Th^*R\Th\big)
\Phi(r,t)\tilde\xi_1dr\].
\end{equation}
Similarly, for any $\tau\in [t,T]$, we can
obtain that
\begin{equation}\label{5.25-eq9}
P^\tau\tilde\xi_1 =
\mE_\tau \[\Phi(T,\tau)^*G
\Phi(T,\tau)\big(Q+\Th^*R\Th\big)\tilde\xi_1 -
\int_{\tau}^T
\Phi(r,\tau)^*\big(Q+\Th^*R\Th\big)
\Phi(r,\tau)\tilde\xi_1dr\].
\end{equation}
It follows from \eqref{5.25-eq8} and
\eqref{5.25-eq9} that for any $\xi\in
L^2_{\cF_t}(\Om;H)$,
\begin{equation}\label{5.25-eq10}
\begin{array}{ll}\ds
\q\mE\big|P^\tau\tilde\xi - P^t\tilde\xi\big|_{H}^{2} \\
\ns\ds \le \cC\Big\{\mE\Big| \mE_\tau
\[\Phi(T,\tau)^*G \Phi(T,\tau) \xi
- \int_{\tau}^T  \Phi(r,\tau)^*\big(Q+\Th^*R\Th\big) \Phi(r,\tau) \xi dr\]  \\
\ns\ds \qq\q - \mE_\tau \[\Phi(T,t)^*G
\Phi(T,t)\xi -
\int_{t}^T  \Phi(r,t)^*\big(Q+\Th^*R\Th\big) \Phi(r,t) \xi dr\] \Big|_{H}^{2}\\
\ns\ds \qq + \mE\Big| \mE_\tau \[\Phi(T,t)^*G
\Phi(T,t)\xi - \int_{t}^T
\Phi(r,t)^*\big(Q+\Th^*R\Th\big)
\Phi(r,t)\xi dr\]
\\
\ns\ds \qq\qq - \mE_t \[\Phi(T,t)^*G
\Phi(T,t)\xi - \int_{t}^T
\Phi(r,t)^*\big(Q+\Th^*R\Th\big)
\Phi(r,t)\xi dr\]\Big|_{H}^{2}\Big\}.
\end{array}
\end{equation}
Since $\dbF$ is the natural filtration
of $W(\cd)$, by the Martingale
representation theorem, we know that
any martingale on $(\Om,\cF,\dbF,\dbP)$
is continuous. Thus,
\begin{equation}\label{5.25-eq11}
\begin{array}{ll}\ds
\lim_{\tau\to t+0}\mE\Big|\mE_\tau
\(\Phi(T,t)^* G \Phi(T,t)\xi -
\int_{t}^T
\Phi(r,t)^*\big(Q+\Th^*R\Th\big)
\Phi(r,t)\xi dr \)  \\
\ns\ds \qq\q\, - \mE_t \(\Phi(T,t)^*G
\Phi(T,t)\xi - \int_{t}^T
\Phi(r,t)^*\big(Q+\Th^*R\Th\big)
\Phi(r,t)\xi dr \)\Big|_{H}^{2} =0.
\end{array}
\end{equation}
On the other hand,
\begin{equation}\label{5.25-eq12}
\begin{array}{ll}\ds
\q \mE\Big|\mE_\tau \(\Phi(T,\tau)^*G
\Phi(T,\tau)\xi - \int_{\tau}^T
\Phi(r,\tau)^*\big(Q+\Th^*R\Th\big)
\Phi(r,\tau)\xi dr\)
\\
\ns\ds \qq - \mE_\tau \(\Phi(T,t)^*G
\Phi(T,t)\xi
- \int_{t}^T  \Phi(r,t)^*\big(Q+\Th^*R\Th\big) \Phi(r,t)\xi dr\)\Big|_{H}^{2}\\
\leq \cC\Big\{\mE\Big|\Phi(T,\tau)^*G \Phi(T,\tau)\xi - \Phi(T,t)^*G\Phi(T,t)\xi \Big|_{H}^{2} \\
\ns\ds \qq + \mE\Big|\int_{\tau}^T
\[\Phi(r,\tau)^*\big(Q+\Th^*R\Th\big)\Phi(r,\tau)\xi -
\Phi(r,t)^*\big(Q+\Th^*R\Th\big)
\Phi(r,t)\xi\]dr \Big|_{H}^{2}
\\
\ns\ds \qq  +  \mE\Big|\int_{t}^{\tau}
\Phi(r,\tau)^*\big(Q+\Th^*R\Th\big)
\Phi(r,\tau)\xi dr\Big|_{H}^{2}\Big\}.
\end{array}
\end{equation}
By Lemma \ref{lemma5}, we know that the three
terms in the right hand side of
\eqref{5.25-eq12} tend to zero as $\tau\to t+0$.
Consequently,   we obtain that
\begin{equation}\label{5.25-eq13}
\begin{array}{ll}\ds
\lim_{\tau\to t+0}\mE\Big|\mE_\tau
\(\Phi(T,\tau)^*G \Phi(T,\tau)\xi -
\int_{\tau}^T
\Phi(r,\tau)^*\big(Q+\Th^*R\Th\big)
\Phi(r,\tau)\xi dr\)
\\
\ns\ds \qq\qq - \mE_\tau \(\Phi(T,t)^*G
\Phi(T,t)\xi - \int_{t}^T
\Phi(r,t)^*\big(Q+\Th^*R\Th\big) \Phi(r,t)\xi
dr\)\Big|_{H}^{2} =0.
\end{array}
\end{equation}
From \eqref{5.25-eq10}, \eqref{5.25-eq11} and
\eqref{5.25-eq13}, we see that for any
$\xi\in L^{2}_{\cF_t}(\Om;H)$,
$$
P^{\cd}\xi\in C_{\dbF}([0,T];L^{2}(\Om;H))
$$
and
$$
|P^{\cd}\xi|_{C_{\dbF}([0,T];L^{2}(\Om;H))}
\leq
\cC\big(\big||Q+\Th^*R\Th|_{\cL(H)}\big|_{L^{\infty}_\dbF(\Om;L^1(0,T))}
+
|G|_{L^{\infty}_{\cF_T}(\Om;\cL(H))}\big)|\xi|_{L^{2}_{\cF_t}(\Om;H)}.
$$
Hence, we get that $P^{\cd}\in
C_{\dbF,w}([0,T];L^{\infty}(\Om;\cL(H)))$.

\ms

{\bf Step 5}. In this step,  we prove
that
\begin{equation}\label{P1}
P(t,\om)=P^{t}(\om), \ \ \ae (t,\om)\in
[0,T]\times\Om.
\end{equation}
To show this, for any $0\leq t_1<t_2<T$, we choose $\tilde
x_1(t_1)=\tilde \xi_1\in L^4_{\cF_{t_1}}(\Om;H)$ and $\tilde
u_1=\tilde v_1=0$ in the equation \eqref{op-fsystem3}, and $\tilde
x_2(t_1)=0$, $\tilde u_2(\cd)=\frac{\chi_{[t_1,t_2]}}{t_2-t_1}\tilde
\xi_2$ with $\tilde \xi_2\in L^4_{\cF_{t_1}}(\Om;H)$ and $\tilde
v_2=0$ in the equation \eqref{op-fsystem4}. By \eqref{3.14-eq28} and
recalling the definition of the evolution operator $\Phi(\cd,\cd)$
(see Lemma \ref{lemma5}), we see that
\begin{equation}\label{s5eq11}
\begin{array}{ll}\ds
 \frac{1}{t_2-t_1}\mE\int_{t_1}^{t_2}\Big\langle
P(r)
\Phi(r,t_1)\tilde \xi_1, \tilde \xi_2 \Big\rangle_H dr\\
\ns\ds =\mE\Big\langle G \Phi(T,t_1)\tilde \xi_1, \tilde
x_{2}(T;t_2,0,\tilde u_2,0) \Big\rangle_H - \mE\int_{t_1}^T
\Big\langle \big(Q+\Th^*R\Th\big)\Phi(r,t_1)\tilde \xi_1, \tilde
x_{2}(r;t_1,0,\tilde u_2,0) \Big\rangle_H dr.
\end{array}
\end{equation}
It is clear that
\begin{equation}\label{vz1}
\tilde x_{2}(r;t_1,0,\tilde
u_2,0)=\left\{
\ba{ll}\ds\int_{t_1}^{r}e^{A(r-\tau)}A_\Th(\tau)\tilde
x_{2}(\tau)d\tau
+\int_{t_1}^{r}e^{A(r-\tau)}C_\Th(\tau)\tilde x_{2}(\tau)dW(\tau)\\
\ns\ds\q +\frac{1}{t_2-t_1}\int_{t_1}^{s}e^{A(r-\tau)}\tilde \xi_2
d\tau,\qq\; r\in [t_1,t_2],\\ \ns\ds \Phi(r,t_2)\tilde x_{2}(t_2),
\qq\qq\qq\qq\; r\in [t_2,T]. \ea\right.
\end{equation}
Hence, we have
$$
\begin{array}{ll}\ds
\q\mE \big|\tilde x_{2}(t_2)-\tilde
\xi_2\big|_H^4\\
\ns\ds \leq
\cC\Big\{\[\(\int_{t_1}^{t_2}\big||A_\Th(\tau)|_{\cL(H)}\big|_{
L^\infty(\Om)}d\tau\)^4
+\(\int_{t_1}^{t_2}\big||C_\Th(\tau)|_{\cL(H)}\big|_{
L^\infty(\Om)}^2d\tau\)^2\]
\mE\sup_{s\in[t_1,t_2]}\big|\tilde
x_{2}(s)\big|_H^4
\\\ns\ds\qq  +
\mE\Big|\frac{1}{t_2-t_1}\int_{t_1}^{t_2}e^{A(t_2-\tau)} \tilde\xi_2
d\tau -\tilde\xi_2\Big|_H^4\Big\}.
\end{array}
$$
This, together with \eqref{3.14-eq33},
implies that
$$
\begin{array}{ll}\ds
\lim_{t_2\to t_1+0}\mE\big|\tilde
x_{2}(t_2)-\tilde\xi_2\big|_H^4\3n&\ds\leq \cC\lim_{t_2\to
t_1+0}\mE\Big|\frac{1}{t_2-t_1}\int_{t_1}^{t_2}e^{A(t_2-\tau)}\tilde\xi_2
d\tau -\tilde\xi_2\Big|_H^4=0.
\end{array}
$$
Therefore, for any $r\in [t_2,T]$, by Lemma \ref{lemma5}, we get that
$$
\begin{array}{ll}\ds
\q\lim_{t_2\to
t_1+0}\mE\big|\Phi(r,t_2)\tilde x_{2}(t_2)-\Phi(r,t_1)\tilde \xi_2\big|_H^4\\
\ns\ds \leq \cC\lim_{t_2\to t_1+0}\(\mE\big|\Phi(r,t_2)\tilde
x_{2}(t_2)-\Phi(r,t_2)\tilde \xi_2\big|_H^4
+\mE\big|\Phi(r,t_2)\tilde \xi_2-\Phi(r,t_1)\tilde \xi_2\big|_H^4\)\\
\ns\ds \leq \cC\lim_{t_2\to t_1+0}\(\mE\big|\tilde x_{2}(t_2)-
\tilde \xi_2\big|_H^4 +\mE\big|\Phi(r,t_2)\tilde
\xi_2-\Phi(r,t_1)\tilde \xi_2\big|_H^4\)=0.
\end{array}
$$
Hence, we obtain that
\begin{equation}\label{vz1q} \lim_{t_2\to t_1+0}
\tilde x_{2}(s) = \Phi(s,t_1)\tilde \xi_2 \ \mbox{ in }
L^4_{\cF_s}(\Om;H),\qq\forall\;s\in [t_2,T].
\end{equation}
By \eqref{3.14-eq33} and \eqref{vz1q},
we conclude that
\begin{equation}\label{10.20eq2}
\begin{array}{ll}\ds
\lim_{t_2\to t_1+0}\[\mE\big\langle G \Phi(T,t_1)\tilde\xi_1,
\tilde x_{2}(T) \big\rangle_H   - \mE\int_{t_1}^T
\big\langle\big(Q+\Th^*R\Th\big)\Phi(s,t_1)\tilde\xi_1, \tilde
x_{2}(s)
\big\rangle_H ds\] \\
\ns\ds = \mE\big\langle G \Phi(T,t_1)\tilde\xi_1,
\Phi(T,t_1)\tilde\xi_2 \big\rangle_H - \mE\int_{t_1}^T \big\langle
\big(Q+\Th^*R\Th\big)\Phi(s,t_1)\tilde\xi_1, \Phi(s,t_1)\tilde\xi_2
\big\rangle_H ds.
\end{array}
\end{equation}

By choosing $\tilde x_1(t_1)=\tilde \xi_1$ and $\tilde u_1=\tilde
v_1=0$ in \eqref{op-fsystem3}, and $\tilde x_2(t_1)=\tilde \xi_2$
and $\tilde u_2=\tilde v_2=0$ in \eqref{op-fsystem4}, from
\eqref{3.14-eq28}, we find that
\begin{equation}\label{s5eq12}
\begin{array}{ll}\ds
\q\mE\big\langle P^{t_1} \tilde \xi_1, \tilde \xi_2
\big\rangle_H\\
\ns\ds =\mE\big\langle G \Phi(T,t_1)\tilde \xi_1, \Phi(T,t_1) \tilde
\xi_2 \big\rangle_H  -  \mE \int_{t_1}^T \big\langle
\big(Q+\Th^*R\Th\big)\Phi(s,t_1)\tilde \xi_1, \Phi(s,t_1) \tilde
\xi_2 \big\rangle_H ds.
\end{array}
\end{equation}

Combining \eqref{s5eq11},
\eqref{10.20eq2} and \eqref{s5eq12}, we
obtain that
\begin{equation}\label{s5eq15}
\lim_{t_2\to t_1+0} \frac{1}{t_2 - t_1}\mE \int_{t_1}^{t_2}
\big\langle P(s) \Phi(s,t_1)\tilde\xi_1, \tilde\xi_2 \big\rangle_H
ds = \mE\big\langle P^{t_1} \tilde\xi_1, \tilde\xi_2 \big\rangle_H.
\end{equation}
By Lemma \ref{lemma2.1}, we see that
there is a monotonically decreasing
sequence $\{t_2^{(k)}\}_{k=1}^\infty$
with $t_2^{(k)}>t_1$ for every
$k\in\dbN$, such that
$$
\lim_{t_2^{(k)}\to
t_1+0}\frac{1}{t_2^{(k)}-t_1}\mE\int_{t_1}^{t_2^{(k)}}\big\langle
P(s) \Phi(s,t_1)\tilde\xi_1, \tilde\xi_2 \big\rangle_H ds =
\mE\big\langle P(t_1)\tilde\xi_1, \tilde\xi_2 \big\rangle_H,\q\ae
t_1\in [0,T).
$$
This, together with \eqref{s5eq15},
implies that
$$
\mE\big\langle P(t_1) \tilde\xi_1, \tilde\xi_2 \big\rangle_H =
\mE\big\langle P^{t_1}\tilde\xi_1, \tilde\xi_2 \big\rangle_H,
\;\mbox{ for }\ae t_1\in [0,T).
$$
Since $\tilde\xi_1$ and $\tilde\xi_2$ are arbitrary elements in
$L^4_{\cF_{t_1}}(\Om;H)$, we conclude \eqref{P1}.  By modifying
$P(\cd,\cd)$ on a null measure subset of $[0,T]\times\Om$,  we get
that $P(\cd,\cd)=P^{t}(\cd)$ for all $t \in [0,T]$. Hence,  $P \in
C_{\dbF,w}\big([0,T];L^{\infty}(\Om;\cL(H))\big)$. Further, we get
from \eqref{3.14-eq28} that
\begin{equation}\label{3.14-eq31}
\begin{array}{ll} \ds \dbE\lan
G\tilde x_{1}(T),\tilde
x_{2}(T)\ran_{H}
+ \dbE\int_t^T \lan [Q(r)+\Th(r)^*R(r)\Th(r)]\tilde x_{1}(r),\tilde x_{2}(r)\ran_H dr\\
\ns\ds =\dbE\lan P(t) \tilde \xi_1,\tilde
\xi_2\ran_{H_\l',H_\l}\!+\!\dbE\!\int_t^T\!\!
\lan P(r)\tilde u_1(r),\tilde
x_{2}(r)\ran_{H_\l',H_\l} dr
\\
\ns\ds \q +\dbE \int_t^T \lan P(r)\tilde
x_{1}(r),\tilde u_2(r)\ran_{H_\l',H_\l} dr +\dbE
\int_t^T \lan P(r)C_{\Th}(r)\tilde
x_{1}(r),\tilde v_2(r)\ran
_{H_\l',H_\l} dr\\
\ns\ds \q +\dbE \int_t^T  \lan P(r)\tilde
v_1(r), C_{\Th}(r)\tilde x_{2}(r) + \tilde
v_2(r) \ran_{H_\l',H_\l} dr
\\
\ns\ds\q  +\dbE\int_t^T \lan \tilde
v_1(r),\L^*(r) \tilde x_{2}(r)\ran_{H_\l,H_\l'}
dr +\dbE\int_t^T \lan \tilde v_2(r),\L(r)\tilde
x_{1}(r)\ran_{H_\l,H_\l'} dr.
\end{array}
\end{equation}

\ss

{\bf Step 6}. In this step, we complete the proof. For any given
$\tilde\xi_{1},\tilde\xi_{2}\in L^4_{\cF_t}(\Om;H)$, let
$\{\tilde\xi_{1,k}\}_{k=1}^\infty,$
$\{\tilde\xi_{2,k}\}_{k=1}^\infty\subset L^4_{\cF_t}(\Om;H_\l)$ such
that for $j=1,2$,
\begin{equation}\label{3.14-eq34}
\lim_{k\to\infty}\tilde\xi_{j,k} = \tilde\xi_{j}\q\mbox{ in
}L^4_{\cF_t}(\Om;H).
\end{equation}
For any given $\tilde u_{1},\tilde
u_{2}\in L^4_{\cF_t}(\Om;L^2(0,T;H))$,
let $\{\tilde u_{1,k}\}_{k=1}^\infty,$
$\{\tilde u_{2,k}\}_{k=1}^\infty\subset
L^4_{\cF_t}(\Om;L^2(0,T;H_\l))$ such
that for $j=1,2$,
\begin{equation}\label{3.14-eq34-1}
\lim_{k\to\infty}\tilde u_{j,k} =
\tilde u_{j}\q\mbox{ in
}L^4_{\cF_t}(\Om;L^2(0,T;H)).
\end{equation}
From \eqref{3.14-eq31}, we have that
\begin{equation}\label{3.14-eq35}
\begin{array}{ll} \ds \dbE\lan
G\tilde x_{1,k}(T),\tilde
x_{2,k}(T)\ran_{H}
+ \dbE\int_t^T \lan [Q(r)+\Th(r)^*R(r)\Th(r)]\tilde x_{1,k}(r),\tilde x_{2,k}(r)\ran_H dr\\
\ns\ds =\dbE\lan P(t) \tilde \xi_{1,k},\tilde
\xi_{2,k}\ran_{H_\l',H_\l}\!+\!\dbE\!\int_t^T\!\!
\lan P(r)\tilde u_{1,k}(r),\tilde
x_{2,k}(r)\ran_{H_\l',H_\l} dr
\\
\ns\ds \q +\dbE \int_t^T \lan P(r)\tilde
x_{1,k}(r),\tilde u_{2,k}(r)\ran_{H_\l',H_\l} dr
+\dbE \int_t^T \lan P(r)C_{\Th}(r)\tilde
x_{1,k}(r),\tilde v_2(r)\ran
_{H_\l',H_\l} dr\\
\ns\ds \q +\dbE \int_t^T  \lan P(r)\tilde
v_1(r), C_{\Th}(r)\tilde x_{2,k}(r) + \tilde
v_2(r) \ran_{H_\l',H_\l} dr
\\
\ns\ds\q  +\dbE\int_t^T \lan \tilde
v_1(r),\L^*(r) \tilde
x_{2,k}(r)\ran_{H_\l,H_\l'} dr +\dbE\int_t^T
\lan \tilde v_2(r),\L(r)\tilde
x_{1,k}(r)\ran_{H_\l,H_\l'} dr,
\end{array}
\end{equation}
where $\tilde x_{1,k}$ (\resp $\tilde x_{2,k}$)
solves \eqref{op-fsystem3} (\resp
\eqref{op-fsystem4}) with $\tilde\eta_1$ (\resp
$\tilde\eta_2$) and $\tilde u_1$ (\resp $\tilde
u_2$) replaced by $\tilde\eta_{1,k}$ (\resp
$\tilde\eta_{2,k}$) and $\tilde u_{1,k}$ (\resp
$\tilde u_{2,k}$).

Similar but simpler to the proof of
\eqref{3.14-eq9}, for $j=1,2$, we can
get that
\begin{equation}\label{3.14-eq36}
\lim_{k\to\infty} \tilde x_{j,k}=\tilde
x_{j} \q\mbox{ in
}L^4_\dbF(\Om;C([t,T];H)).
\end{equation}
Noting that $P\in
C_{\dbF,w}\big([0,T];L^{2}(\Om;\cL(H))\big)$,
from \eqref{3.14-eq34},
\eqref{3.14-eq34-1} and
\eqref{3.14-eq36}, by taking
$k\to\infty$ in \eqref{3.14-eq35}, we
obtain that
\begin{equation}\label{3.14-eq31-1}
\begin{array}{ll} \ds \dbE\lan
G\tilde x_{1}(T),\tilde
x_{2}(T)\ran_{H}
+ \dbE\int_t^T \lan [Q(r)+\Th(r)^*R(r)\Th(r)]\tilde x_{1}(r),\tilde x_{2}(r)\ran_H dr\\
\ns\ds =\dbE\lan P(t) \tilde \xi_1,\tilde
\xi_2\ran_{H}\!+\!\dbE\!\int_t^T\!\! \lan
P(r)\tilde u_1(r),\tilde x_{2}(r)\ran_{H} dr
\\
\ns\ds \q +\dbE \int_t^T \lan P(r)\tilde
x_{1}(r),\tilde u_2(r)\ran_{H} dr +\dbE \int_t^T
\lan P(r)C_{\Th}(r)\tilde x_{1}(r),\tilde
v_2(r)\ran
_{H} dr\\
\ns\ds \q +\dbE \int_t^T  \lan P(r)\tilde
v_1(r), C_{\Th}(r)\tilde x_{2}(r) + \tilde
v_2(r) \ran_{H} dr
\\
\ns\ds\q  +\dbE\int_t^T \lan \tilde
v_1(r),\L^*(r) \tilde x_{2}(r)\ran_{H_\l,H_\l'}
dr +\dbE\int_t^T \lan \tilde v_2(r),\L(r)\tilde
x_{1}(r)\ran_{H_\l,H_\l'} dr.
\end{array}
\end{equation}

  Assume that
$$
(P_1(\cd),\Lambda_1(\cd)), (P_2(\cd),\Lambda_{2}(\cd))\in
C_{\dbF,w}([0,T]; L^{\infty}(\Om;\cL(H))) \times
L^2_{\dbF}(0,T;\cL_2(H;H_\l))
$$
satisfy \eqref{3.14-eq31-1}.

\ms

By choosing $\tilde u_1=\tilde u_2=0$
and $\tilde v_1=\tilde v_2=0$, from
\eqref{3.14-eq31-1}, we have that for
any $t\in[0,T)$ and
$\tilde\xi_1,\tilde\xi_2 \in
L^4_{\cF_4}(\Om;H)$,
\begin{equation}\label{8.22-eq2}
 \mE\langle P_1(s)\tilde\xi_1,\tilde\xi_2 \rangle_H
= \mE\langle P_2(s)\tilde\xi_1,\tilde\xi_2\rangle_H.
\end{equation}
This implies
\begin{equation}\label{3.14-eq38}
P_1(\cd)=P_2(\cd).
\end{equation}
\vspace{0.1cm}

Let $\tilde v_2=0$ in \eqref{op-fsystem4}. By \eqref{3.14-eq38},  we
see that for any $\tilde\xi_1\in L^4_{\cF_0}(\Om;H)$, $\tilde u_1\in
L^4_\dbF(\Om;L^2(0,T;H))$, $\tilde v_1 \in
L^4_\dbF(\Om;L^2(0,T;H_\l))$,
\begin{equation}\label{8.22-eq9}
0=\mE \int_0^T \big\langle\tilde v_1(r),
\big(\Lambda_1(r) - \Lambda_2(r)\big)\tilde
x_2(r)\big\rangle_{H_\l,H_\l'}dr.
\end{equation}
Consequently,
\begin{equation}\label{9.7-eq13}
\big(\Lambda_1 - \Lambda_2 \big)\tilde x_2 = 0 \;\mbox{ in }\;
L^{\frac{4}{3}}_\dbF(\Om;L^2(0,T;H_\l')).
\end{equation}
This, together with Lemma \ref{lm11-1}, implies that\vspace{1mm}
$$\Lambda_1=\Lambda_2 \mbox{ in
}L^2_\dbF(0,T;\cL_2(H;H_\l')).$$ Hence, the desired uniqueness
follows, which also implies the uniqueness of the
$H_\l$-transposition solution to \eqref{Relaxed-t-solution-Riccati}.

From \eqref{3.14-eq31-1}, we see that if $(P,\L)$ satisfies
\eqref{3.14-eq31-1}, then $(P^*,\L^*)$ also satisfies
\eqref{3.14-eq31-1}. Hence, $(P,\L)=(P^*,\L^*)$. This concludes that $(P(\cd),\L(\cd))\in
C_{\dbF,w}\big([0,T];L^\infty(\Omega;\cL(H))\big)\times
L^2_{\dbF}\big(0,T;L^2(\Omega;$ $\dbS_2(H;H_\l'))\big)$.
\endpf

\subsection{A novel characterization of optimal feedback operator via $H_\l$-transposition solution}

The following result, which plays indispensable role in proving Theorem \ref{5.7-th1}, reveals another characterization of the optimal feedback operator $\Th$.
\begin{lemma}\label{th1}
Let ({\bf AS0})--({\bf AS4}) hold and $\Th\in \Upsilon_2(H;U)\cap
\Upsilon_2(H_\l;\wt U)$, and $(\overline P,\overline\L)$ be the
$H_\l$-transposition solution of (\ref{Relaxed-t-solution-Riccati})
corresponding to $\Th$. Let
\begin{equation}\label{th1-eq1}
\overline K=R+D^{*}\overline PD,\ \ \overline L=D^{*} \overline P C
+B^*\overline P +D^*\overline\L.
\end{equation}
Then $\Th$ is an optimal feedback
operator for Problem (SLQ) if and only if
\begin{equation}\label{Nece-1}
\overline K(t,\om)\geq0, \mbox{ for a.e. }(t,\om)\in (0,T)\times\Om,
\end{equation}
\begin{equation}\label{Nece-3}
\overline L \in  \Upsilon_2(H;U)\cap \Upsilon_2(H_\l;\wt U),
\end{equation}
and
\begin{equation}\label{Nece-2}
\overline K(t,\om)\Th(t,\om)+\overline L(t,\om)=0 \q\mbox{ in
}\;\cL(H;U) \mbox{ for a.e. }(t,\om)\in (0,T)\times\Om.
\end{equation}
\end{lemma}
%
%


Before proving Lemma \ref{th1}, we need another preliminary result.
The proof is almost
the same as the one for \cite[Lemma
3.4]{Hu-et-al-2017}, where a very similar result
is proved when $H=\dbR^l$. Here we provide it
for the convenience of readers.
\begin{lemma}\label{Lebesgue-differential-theorem}
Suppose $Y(\cd)\in L^2_{\dbF}(0,T; H)$ is a
given process satisfying that for a.e. $t\in
[0,T]$, there is a sequence of decreasing
positive numbers $\{\e_n\}_{n=1}^\infty$ such
that
\begin{equation}\label{3.4-eq16}
\begin{array}{ll}
\ns\ds \lim_{\e_n\to 0}
\int_t^{t+\e_n}\frac{\dbE_t Y(r)}{\e_n}dr=0
\mbox{ in }H,\ \ \q\dbP\mbox{-a.s.}
\end{array}
\end{equation}
Then we have $Y=0$, for a.e.
$(t,\om)\in[0,T)\times\Om$.
\end{lemma}

\ms

\it Proof. \rm Since  $L^2_{\cF_T}(\Omega;H)$ is
separable, it follows from the (deterministic)
Lebesgue differentiation theorem that there is a
countable dense subset $\cD\subset
L^{\infty}_{\cF_T}(\Omega;H)$ of
$L^2_{\cF_T}(\Omega;H)$ and a full measure set
$\G_1\subset [0,T)$ (depending on $Y(\cd)$) such
that  for  all $t\in \G_1$, we have
$$\begin{array}{ll}
\ns\ds \lim_{\e\to+0}\int_t^{t+\e}\frac{\dbE\lan
Y(r),\eta\ran_H dr}{\e}=\dbE\lan
Y(t),\eta\ran_H,\ \ \forall \eta\in \cD,
\end{array}
$$
and
\begin{equation}\label{3.4-eq14}
\begin{array}{ll}
\ds \lim_{\e\to
0}\int_t^{t+\e}\frac{\dbE|Y(r)|^2}{\e}dr
=\dbE|Y(t)|^2.
\end{array}
\end{equation}

For any $\eta\in\cD$, put
\begin{equation}\label{3.4-eq14-1}
\tilde\eta(r)\=\dbE_r\eta \mbox{ for }r\in[0,T].
\end{equation}
Since $\eta$ is essentially bounded, so is
$\tilde\eta(\cd)$. Further, since
$\tilde\eta(\cd)$ is a martingale on the natural
filtration of the Brownian motion $W(\cd)$, it
is continuous.

From \eqref{3.4-eq14-1}, we have that
$$\dbE\lan
Y(r),\eta\ran_H=\dbE\lan
Y(r),\tilde\eta(r)\ran_H  \mbox{ for }r\in[0,T],
$$
and  for
all $t\in \G_1$,
$$\begin{array}{ll}
\ds \Big|\lim_{\e\to 0}\int_t^{t+\e}
\frac{\dbE\lan Y(r),\tilde\eta(r)-\tilde\eta(t)\ran_H}{\e}dr\Big|\\
\ns\ds \leq \lim_{\e\to
0}\Big[\int_t^{t+\e}\frac{\dbE|Y(r)|_H^2}{\e}dr
\Big]^{\frac 1
2}\Big[\int_t^{t+\e}\frac{\dbE|\tilde\eta(r)-\tilde\eta(t)|_H^2}{\e_n}dr
\Big]^{\frac 1 2}\\
\ns\ds \leq2\lim_{\e\to
0}\Big[\int_t^{t+\e}\frac{\dbE|Y(r)|_H^2}{\e}dr
\Big]^{\frac 1 2}
\sqrt{\sup_{s\in[t,t+\e]}\dbE|\tilde\eta
(r)-\tilde\eta(t)|_H^2} \\
\ns\ds \leq2\big(\dbE|Y(t)|_H^2
\big)^{\frac 1 2} \lim_{\e\to
0}\sqrt{\sup_{s\in[t,t+\e]}\dbE|\tilde\eta
(r)-\tilde\eta(t)|_H^2}=0,
\end{array}
$$
where the last inequality is due to the
continuity of $\tilde\eta(\cd)$. Hence, for every
$\eta\in\cD$, we have
$$
\begin{array}{ll}
\ds \dbE\big[\lan
Y(t),\tilde\eta(t)\ran_H\big]\3n&\ds
=\lim_{\e\to
0}\int_t^{t+\e}\frac{\dbE\lan
Y(r),\tilde\eta(r)\ran_H}{\e}dr=\lim_{\e\to
0}\int_t^{t+\e}\frac{\dbE\lan
Y(r),\tilde\eta(t)\ran_H}{\e}dr\\
\ns&\ds =\lim_{\e\to 0}\dbE\Big\langle
\int_t^{t+\e}\frac{ \dbE_tY(r)}{\e}dr
,\tilde\eta(t)\Big\rangle_H.
\end{array}
$$
For any $t\in[0,T)$ such that $t+\e\leq T$,
\begin{equation}\label{3.4-eq15}
\begin{array}{ll}
\ds \dbE \Big|\int_t^{t+\e}\frac{\dbE_t
Y(r)}{\e}dr\Big|_H^2 \leq
\(\int_t^{t+\e}\frac{1}{\e}dr\)
\mE\int_t^{t+\e}\frac{|\dbE_t
Y(r)|_H^2}{\e}dr \leq\frac
{1}{\e}\dbE\int_t^{t+\e}|Y(r)|_H^2dr.
\end{array}
\end{equation}
From \eqref{3.4-eq14} and \eqref{3.4-eq15}, we
know that for any $t\in\G_1$, there exists a
constant $\d(t)>0$ such that
$$
\begin{array}{ll}
\ds \dbE\Big|\int_t^{t+\e}\frac{\dbE_t
Y(r)}{\e}dr\Big|_H^2 < \dbE|Y(t)|_H^2+1,\
\ \forall \e\in(0,\d(t)).
\end{array}
$$
This implies that
$$
\cY(\e;t)\=\int_t^{t+\e}\frac{\dbE_tY(r)}{\e}dr
$$
is uniformly integrable in $\e_n\in(0,\d(t))$.
Hence, by \eqref{3.4-eq16}, for a.e. $t\in\G_1$,
there is a sequence $\{\e_n\}_{n=1}^\infty$ of
decreasing positive numbers such that
$$
\begin{array}{ll}
\ds \lim_{\e_n\to
0}\dbE\Big|\int_t^{t+\e_n}\frac{\dbE_t
Y(r)}{\e_n}dr\Big|_H=\dbE\lim_{\e_n\to
0} \Big|\int_t^{t+\e_n}\frac{\dbE_t
Y(r)}{\e_n}dr\Big|_H=\dbE
\Big|\lim_{\e_n\to
0}\int_t^{t+\e_n}\frac{\dbE_t
Y(r)}{\e_n}dr\Big|_H=0.
\end{array}
$$
It follows from the essential  boundedness of
$\tilde\eta(\cd)$ that there exists a constant
$\cC>0$ such that
$$\begin{array}{ll}
\ds \lim_{\e_n\to
0}\Big|\dbE\Big\langle\int_t^{t+\e_n} \frac{
\dbE_tY(r)}{\e_n}dr
,\tilde\eta(t)\Big\rangle_H\Big|\leq
\cC\lim_{\e_n\to 0}\dbE
\Big|\int_t^{t+\e_n} \frac{
\dbE_tY(r)}{\e_n}dr \Big|_H= 0.
\end{array}
$$
Consequently, for any $\eta\in\cD$,
$$
\dbE \lan Y(t),\eta\ran_H =0 \; \text{
for a.e. }t\in\G_1.
$$
This concludes that
$$Y(t)=0 \; \text{ for a.e.
} (t,\om)\in[0,T]\times\Om.
$$
\endpf

{\it Proof of Lemma \ref{th1}}:\, The ``if" part.   Let us divide the proof
into several steps.

\ms

{\bf Step 1}. For any $\tilde v\in \wt U$ and given $t\in[0,T)$, let
\begin{equation}\label{3.4-eq19}
\xi\=\dbE_ t\int_{ t}^{ t+\e}\big[\overline K(r)\Th(r)+\overline
L(r)\big]^* \tilde v dr\in
 L^2_{\cF_t}(\Omega;H).
\end{equation}
Given $\f \in L^2_{\dbF}(0,T;\wt
U)$, consider the following equation:
\begin{equation}\label{3.4-eq26}
\left\{
\begin{array}{ll}
\ds d X =\big[\big(A+ A_{\Th} \big)
X +B \f \big]dr +\big(C_{\Th} X +D \f \big)dW(r) & \mbox{ in }(t,T],\\
\ns\ds X(t)=\xi.
\end{array}
\right.
\end{equation}
From ({\bf AS4}), we have that $D\f\in
L^2_{\dbF}(0,T;H_\l)$. Thanks to ({\bf AS0}),
the solution $X\in L^2_{\dbF}(\Omega;C([t,T];$
$H))$. By the definition of the
$H_\l$-transposition solution to
\eqref{Relaxed-t-solution-Riccati} and Lemma
\ref{3.14-lm1}, we have
\begin{eqnarray}\label{Applied-to-SLQ}
&&\dbE\lan  G X(T),X(T)\ran_{H}+ \dbE\int_t^T \lan [Q(r)+\Th(r)^*R(r)\Th(r)]X(r),X(r)\ran_H  dr\nonumber\\
&& =\dbE\lan \overline P(t) \xi, \xi\ran_{H}+\dbE\int_t^T \lan \overline P(r)B(r)\f(r),X(r)\ran_{H}  dr+\dbE\int_t^T \lan \overline P(r)X(r),B(r)\f(r)\ran_H  dr\nonumber\\
&& \q +\dbE\!\int_t^T\!\! \lan \overline
P(r)C_{\Th}(r)X(r),D(r)\f(r)\ran _{H} dr\!+\!\dbE\!\int_t^T\!\!
\lan \overline P(r)D(r)\f(r), C_{\Th}(r)X(r)\!+\!D(r)\f(r)\ran_H dr\nonumber\\
&&\q +\dbE\int_t^T \lan
D(r)\f(r),\overline\L^*(r)X(r)\ran_{H_\l,H_\l'} dr +\dbE\int_t^T
\lan
D(r)\f(r),\overline\L(r)X(r)\ran_{H_\l,H_\l'}  dr\\
&& =\dbE\lan \overline P(t) \xi, \xi\ran_{H}+\dbE\int_t^T \lan \overline P(r)B(r)\f(r),X(r)\ran_{H}  dr+\dbE\int_t^T \lan \overline P(r)X(r),B(r)\f(r)\ran_H  dr\nonumber\\
&& \q +2\dbE\!\int_t^T\!\! \lan \overline
P(r)C_{\Th}(r)X(r),D(r)\f(r)\ran _{H} dr\!+\!\dbE\!\int_t^T\!\!
\lan \overline P(r)D(r)\f(r), D(r)\f(r)\ran_H dr\nonumber\\
&&\q +2\dbE\int_t^T \lan
D(r)\f(r),\overline\L^*(r)X(r)\ran_{H_\l,H_\l'} dr .\nonumber
\end{eqnarray}

Denote by $\overline x(\cd)\in L^2_{\dbF}(\Omega;C([t,T];H))$ the
solution to \eqref{3.4-eq26} with $\xi$ given by \eqref{3.4-eq19}
and $\f=0$. Let $\e>0$ be such that $t+\e\leq T$. Write
$x^{v,\e}(\cd)$ for the solution to \eqref{3.4-eq26} with  $\xi$ given
by \eqref{3.4-eq19} and $\f(\cd)= \chi_{[t,t+\e]}(\cd)v$ for some
$v\in\wt U$.
It follows from \eqref{3.4-eq26} and the choice of $\f$ that
$x^{v,\e}(t)=\overline x(t)$.

Let $t_\e\in [t,t+\e]$ be a stopping
time to be fixed later. For
$r\in[t,t_\e]$,
\begin{equation}\label{3.14-eq37}
\begin{array}{ll}\ds
\mE\sup_{r\in [t,t_\e]}|x^{v,\e}(r)-\overline x(r)|_H^2\\
\ns\ds = \mE\sup_{r\in [t,t_\e]}\Big|\int_t^{r}e^{A(\tau-t)}A_\Th(\tau)
\big(x^{v,\e}(\tau)-\overline x(\tau)\big)d\tau +
\int_t^{r}e^{A(\tau-t)}C_\Th(\tau)
\big(x^{v,\e}(\tau)-\overline x(\tau)\big)dW(\tau)\\
\ns\ds\qq\qq\; + \int_t^{r}e^{A(\tau-t)}B(\tau)\f(\tau)d\tau+
\int_t^{r}e^{A(s-t)}D(\tau)\f(\tau)dW(\tau) \Big|_H^2\\
\ns\ds \leq \cC\mE\sup_{r\in [t,t_\e]}\[
\Big|\int_t^{r}\!e^{A(\tau-t)}A_\Th(\tau)
\big(x^{v,\e}(\tau)\!-\!\overline x(\tau)\big)d\tau \Big|_H^2 \!+
 \Big|\int_t^{r}\!e^{A(\tau-t)}C_\Th(\tau)
\big(x^{v,\e}(\tau)\!-\!\overline x(\tau)\big)dW(\tau)
\Big|_H^2\\
\ns\ds \qq\qq\q  +
\Big|\int_t^{r}e^{A(\tau-t)}B(\tau)\f(\tau)d\tau\Big|_H^2+
\Big|\int_t^{r}e^{A(\tau-t)}D(\tau)\f(\tau)dW(\tau)
\Big|_H^2\]\\
\ns\ds \leq \cC \Big|\(\int_t^{t_\e} |A_\Th(\tau)|_{\cL(H)}d\tau\)^2
+
\int_t^{t_\e}|C_\Th(\tau)|_{\cL(H)}^2d\tau\Big|_{L^\infty_{\cF_T}(\Om)}
\mE\sup_{r\in [t,t_\e]}\big|x^{ v,\e}(r)-\overline x(r)\big|_H^2
d\tau + \cC(t_\e-t),
\end{array}
\end{equation}
where the constant $\cC$ is independent
of the stopping time $t_\e$. Let us
choose $t_\e$ as follows:
$$
t_\e=\begin{cases}\ds\inf\Big\{r\in [t,t+\e]\Big| \(\int_t^{r}
|A_\Th(\tau)|_{\cL(H)}d\tau\)^2 +
 \int_t^{t_\e}|C_\Th(\tau)|_{\cL(H)}^2d\tau=\frac{1}{2\cC}\Big\},\\
\ns\ds\qq\q \mbox{ if } \Big\{r\in [t,t+\e]\Big| \(\int_t^{r}
|A_\Th(\tau)|_{\cL(H)}d\tau\)^2 +
 \int_t^{t_\e}|C_\Th(\tau)|_{\cL(H)}^2d\tau=\frac{1}{2\cC}\Big\}\neq\emptyset,
\\
\ns\ds t+\e, \q \mbox{otherwise}.
\end{cases}
$$
From \eqref{3.14-eq37}, we have that
$$
\mE\sup_{r\in [t,t_\e]}|x^{ v,\e}(r)-\overline x(r)|_H^2\leq
\cC(t_\e-t),\qq \forall\, r\in [t, t_\e].
$$
If $t_\e=\e$, then we get
\begin{equation}\label{3.4-eq17}
\mE\sup_{r\in [t,\e]}|x^{ v,\e}(r)-\overline x(r)|_H^2\leq \cC\e,\qq
\forall\, r\in [t, t+\e].
\end{equation}
Otherwise, noting that
$\big||A_\Th|_{\cL(H)}\big|\in
L^\infty_\dbF(\Om;L^1(0,T))$ and
$\big||C_\Th|_{\cL(H)}\big|\in
L^\infty_\dbF(\Om;L^2(0,T))$, we can repeat the above
argument  in
finite steps to get
\eqref{3.4-eq17}.

\ms

{\bf Step 2}.
From \eqref{Applied-to-SLQ}, we
have that
\begin{equation}\label{3.4-eq29}
\begin{array}{ll}\ds
2\cJ(t,\xi;\Th \overline x) \\
\ns\ds =\dbE\lan G\overline x(T),\overline x(T)\ran_{H}+
\dbE\int_t^T \big[\lan
Q(r)\overline x(r),\overline x(r)\ran_H +\langle R(r)\Th(r)\overline x(r),\Th(r)\overline x(r)\ran_U \big]dr\\
\ns \ds=\dbE\lan G\overline x(T),\overline x(T)\ran_{H}+
\dbE\int_t^T \lan
[Q(r)+\Th(r)^*R(r)\Th(r)]\overline x(r),\overline x(r)\ran_H dr\\
\ns \ds = \dbE\lan \overline P(t)\xi,\xi\ran_H
\end{array}
\end{equation}
and that
\begin{equation}\label{3.4-eq27}
\begin{array}{ll}\ds
2\cJ(t,\xi;\Th x^{ v,\e}+\f)\\
\ns\ds  =\dbE\lan
Gx^{ v,\e}(T),x^{v,\e}(T)\ran_{H} + \dbE\int_t^T \Big[\lan
Q(r)x^{ v,\e}(r),x^{ v,\e}(r)\ran_H \\
\ns\ds \q  +\lan R(r)\big(\Th(r)x^{ v,\e}(r)+\f(r)\big),\Th(r)x^{ v,\e}(r)+\f(r)\ran_H\Big] dr\\
\ns\ds  =\dbE\lan
Gx^{ v,\e}(T),x^{ v,\e}(T)\ran_{H}+
\dbE\int_t^T \lan
\big(Q(r)+ \Th(r)^*\overline P(r)\Th(r) \big)x^{ v,\e}(r),x^{ v,\e}(r)\ran_H dr\\
\ns\ds \q +  \dbE\int_t^T \Big[ 2\lan
R(r) \Th(r)x^{ v,\e}(r), \f(r)\ran_U +
\lan R(r)\f(r), \f(r)\ran_U\Big] dr.
\end{array}
\end{equation}
By \eqref{Applied-to-SLQ} again, we get that
\begin{eqnarray}\label{3.4-eq28}
&&\dbE\lan
Gx^{ v,\e}(T),x^{ v,\e}(T)\ran_{H}+
\dbE\int_t^T \lan
\big(Q(r)+ \Th(r)^*\overline P(r)\Th(r) \big)x^{ v,\e}(r),x^{ v,\e}(r)\ran_H dr \nonumber\\
&&\3n \q +  \dbE\int_t^T \Big[ 2\langle R(r) \Th(r)x^{ v,\e}(r), \f(r)\ran_U + \langle R(r)\f(r), \f(r)\ran_U\Big] dr\nonumber\\
&&\3n  =\dbE\lan \overline P(t) \xi, \xi\ran_{H}+\dbE\int_t^T \lan \overline P(r)B(r)\f(r),x^{ v,\e}(r)\ran_{H} dr \\
&& \3n \q +\dbE\!\int_t^T\!\! \lan \overline P(r)x^{
v,\e}(r),B(r)\f(r)\ran_H dr \!+\!2\dbE\!\int_t^T\!\! \lan \overline
P(r)\big(C(r)\!+\!D(r)\Th(r)\big)x^{v,\e}(r),D(r)\f(r)\ran
_{H} dr\nonumber\\
&&\3n \q + \dbE \int_t^T \lan \overline P(r)D(r)\f(r),
D(r)\f(r)\ran_H dr +2\dbE\int_t^T \lan D(r)\f(r),\overline\L^*(r)x^{
v,\e}(r)\ran_{H_\l,H_\l'} dr\nonumber
\\
&&\3n \q + \dbE\int_t^T \Big[ 2\langle
R(r) \Th(r)x^{ v,\e}(r), \f(r)\ran_U +
\lan R(r)\f(r), \f(r)\ran_U\Big]
dr\nonumber
\\
&& \3n = \dbE\big[\lan \overline P(t)\xi,\xi\ran_H\big]+\dbE\int_
t^{ t+\e}\Big[ \lan \f(r), \overline K(r)  \f(r)\ran_U
+2\lan\big[\overline K(r)\Th(r) +\overline L(r)\big] x^{ v,\e}(r),
\f(r)\ran_U \Big]dr.\nonumber
\end{eqnarray}
Combining
\eqref{3.4-eq29}--\eqref{3.4-eq28}, we
obtain that
\begin{eqnarray*}
&& 2\cJ(t,\xi;\Th x^{ v,\e}+\f)
\\
&&   =2\cJ(t,\xi;\Th \overline x) +\dbE\int_t^{t+\e}\Big[ \lan
\f(r), \overline K(r)
 \f(r)\ran_U +2\lan\big[\overline K(r)\Th(r) +\overline
L(r)\big] x^{v,\e}(r), \f(r)\ran_U \Big]dr.
\end{eqnarray*}
Since $\Th$ is an optimal feedback operator for
Problem (SLQ), for any $ t\in[0,T)$, we have
\begin{equation}\label{3.4-eq3}
\begin{array}{ll}\ds
2\cJ(t,\xi;\Th x^{v,\e}+\f)-2\cJ(t,\xi;\Th \overline x)\\
\ns\ds  =\dbE\int_t^{t+\e}\Big[ \lan  \f(r),\overline K(r)
\f(r)\ran_U+2 \lan\big[\overline K(r)\Th(r)+\overline
L(r)\big]x^{v,\e}(r), \f(r)\ran_{\wt U', \wt U}\Big]dr\geq 0.
\end{array}
\end{equation}
Since $\f(\cd)$ is a bounded process, we have that
\begin{equation}\label{3.4-eq4}
\begin{array}{ll}
\ds \Big|\dbE\int_t^{t+\e} \frac{\lan\big[\overline
K(r)\Th(r)+\overline L(r)\big]\big[x^{v,\e}(r)-\overline x(r)\big],
\f(r)\ran_{\wt U', \wt
U}}{\e} dr\Big|\\
\ns\ds \leq \cC\dbE\int_t^{t+\e} \frac{\big(|\overline
K(r)\Th(r)|_{\cL(H;\wt U')}+|\overline L(r)|_{\cL(H;\wt U')}\big)
|x^{v,\e}(r)-\overline x(r)|_H}{\e}  dr
\\
\ns\ds \leq \frac{\cC}{\e}\Big[\int_t^{t+\e} \dbE\big(|\overline
K(r)\Th(r)|_{\cL(H;\wt U')}+|\overline L(r)|_{\cL(H;\wt U')}\big)^2
dr\Big]^{\frac 1 2}\Big[\int_t^{t+\e}\dbE|x^{v,\e}(r)-\overline
x(r)|_H^2dr\Big]^{\frac 1 2} \\
\ns\ds \leq  \cC \Big[\int_t^{t+\e} \dbE\big(|\overline
K(r)\Th(r)|_{\cL(H;\wt U')}+|\overline L(r)|_{\cL(H;\wt U')}\big)^2
dr\Big]^{\frac 1
2}\Big[\frac{1}{\e^2}\int_t^{t+\e}\dbE|x^{v,\e}(r)-\overline
x(r)|_H^2dr\Big]^{\frac 1 2}
\\
\ns\ds \leq  \cC \Big[\int_t^{t+\e} \dbE\big(|\overline
K(r)\Th(r)|_{\cL(H;\wt U')}+|\overline L(r)|_{\cL(H;\wt U')}\big)^2
dr\Big]^{\frac 1
2}\Big[\frac{1}{\e}\sup_{r\in[t,t+\e]}\dbE|x^{v,\e}(r)-\overline
x(r)|_H^2\Big]^{\frac 1 2}.
\end{array}
\end{equation}
This, together with \eqref{3.4-eq17}, implies
that
\begin{equation}\label{3.4-eq18}
\lim_{\e\to 0}\Big|\dbE\int_t^{t+\e} \frac{\lan\big[\overline
K(r)\Th(r)+\overline L(r)\big]\big[x^{v,\e}(r)-\overline x(r)\big],
\f(r)\ran_{\wt U', \wt U}}{\e} dr\Big|=0.
\end{equation}
Therefore, we
get that for any $v\in\wt U$,
\begin{equation}\label{3.4-eq21}
f(t,v)\=\limsup_{\e\rightarrow0}\frac{1}{\e}\dbE\int_
t^{t+\e}\big\langle v, \overline K(r)v\big\rangle_Udr
+\limsup_{\e\rightarrow0}\frac{2}{\e}
\dbE\int_t^{t+\e}\Big\langle\big[\overline K(r)\Th(r)+\overline
L(r)\big]\xi,v\Big\rangle_{\wt U', \wt U } dr\geq 0.
\end{equation}
By \eqref{3.4-eq21}, we see that for all
$n\in\dbN$ and $v\in \wt U$,
\begin{equation}\label{3.4-eq11}
nf\(t,\frac{v}{n}\)\geq 0, \q \mbox{ for a.e.
}t\in[0,T].
\end{equation}
Letting $n\to\infty$ in \eqref{3.4-eq11}, we have
that
\begin{equation}\label{3.4-eq12}
\limsup_{\e\to0}\dbE \int_t^{t+\e}\Big\langle \frac{\big[\overline
K(r)\Th(r)+\overline L(r)\big]\xi}{\e}, v\Big\rangle_{\wt U', \wt U
} dr\geq 0,\q \forall\,v\in \wt U.
\end{equation}
By the arbitrariness of $v\in \wt U$, one has
\bel{A-key-contribution}\begin{array}{ll}
\ns\ds \limsup_{\e\rightarrow0}\dbE \int_t^{t+\e}\frac{\big[
\overline K(r)\Th(r)+\overline L(r)\big]\xi}{\e} dr =0\; \mbox{ in
}\;\wt U'.
\end{array}
\ee
By the definition of $\xi$ (see
\eqref{3.4-eq19}), we know that for any $\tilde
v\in \wt U$,
$$
\limsup_{\e\to0}\dbE\Big[\dbE_t\int_t^{t+\e}\frac{\big[ \overline
K(r)\Th(r) +\overline L(r)\big]}{\e}
dr\(\dbE_t\int_t^{t+\e}\frac{\big[\overline K(r)\Th(r)+\overline
L(r)\big]^*\tilde v}{\e} dr\Big)\] =0.
$$
Consequently,
$$
\begin{array}{ll}\ds
\Big\langle\limsup_{\e\to0}\dbE\Big[\dbE_t\int_t^{t+\e}\frac{\big[
\overline K(r)\Th(r) +\overline L(r)\big]}{\e}
dr\(\dbE_t\int_t^{t+\e}\frac{\big[\overline K(r)\Th(r)+\overline
L(r)\big]^*\tilde v}{\e}
dr\Big)\], \tilde v\Big\rangle_{\wt U',\wt U} \\
\ns\ds =
\limsup_{\e\to0}\dbE\Big[\dbE_t\Big\langle\int_t^{t+\e}\frac{\big[
\overline K(r)\Th(r) +\overline L(r)\big]}{\e}
dr\(\dbE_t\int_t^{t+\e}\frac{\big[\overline K(r)\Th(r)+\overline
L(r)\big]^*\tilde v}{\e}
dr\Big), \tilde v\Big\rangle_{\wt U',\wt U}\]\\
\ns\ds =
\limsup_{\e\to0}\dbE\Big|\dbE_t\int_t^{t+\e}\frac{\big[\overline
K(r)\Th(r)+\overline L(r)\big]^*\tilde v}{\e} dr\Big|^2_H =0 .
\end{array}
$$
Hence there exists a  sequence
$\{\e_n\}_{n=1}^\infty$ of positive numbers
(depending on $\wt v$ and $t$) such that
$$\begin{array}{ll}
\ds \lim_{\e_n\rightarrow0^+}
\dbE_t\int_t^{t+\e_n}\frac{\big[\overline K(r)\Th(r)+\overline
L(r)\big]^*\tilde v}{\e_n} dr=0.
\end{array}
$$
This, together with Lemma
\ref{Lebesgue-differential-theorem}, implies
that
\begin{equation*}
\big[\overline K(t)\Th(t)+\overline L(t)\big]^*\tilde v=0\; \mbox{
for a.e. } t \in [0,T],\q\dbP\mbox{-a.s.}
\end{equation*}
By separability of $\wt U$, we can find a
density countable subset $\wt U_1$ of $\wt U$
and a full measurable subset $\G_2$ of $[0,T]$
such that such that for all $\tilde v\in \wt
U_1$,
\begin{equation*}
\big[\overline K(t)\Th(t)+\overline L(t)\big]^*\tilde v=0\; \mbox{
for all } t \in \G_2,\q\dbP\mbox{-a.s.}
\end{equation*}
This implies that for all $\tilde v\in \wt U$,
\begin{equation*}
\big[\overline K(t)\Th(t)+\overline L(t)\big]^*\tilde v=0\; \mbox{
for all } t \in \G_2,\q\dbP\mbox{-a.s.}
\end{equation*}
Consequently,
\begin{equation*}
\big[\overline K(t)\Th(t)+\overline L(t)\big]^* =0\; \mbox{ for a.e.
} t \in [0,T],\q\dbP\mbox{-a.s.,}
\end{equation*}
Thus,
\begin{equation*}
\overline K(t)\Th(t)+\overline L(t) =0 \mbox{ in }\cL(H;\wt U')
\;\mbox{ for a.e. } t \in [0,T],\q\dbP\mbox{-a.s.}
\end{equation*}
This, together with \eqref{th1-eq1}, implies
that
\begin{equation}\label{3.4-eq34}
D(t)^*\overline\L(t)=-  D(t)^{*} \overline P(t) C(t)
-B(t)^*\overline P(t) -\overline K(t)\Th(t).
\end{equation}
Noting that the right hand side of
\eqref{3.4-eq34} belongs to
$\Upsilon_2(H;U)\cap
\Upsilon_2(H_\l;\wt U)$, we have
\begin{equation}\label{3.4-eq35}
D^*\overline\L\in \Upsilon_2(H;U)\cap \Upsilon_2(H_\l;\wt U).
\end{equation}
Consequently,  \eqref{Nece-3} and \eqref{Nece-2} hold.

\ms

{\bf Step 3}. In this step, we   prove \eqref{Nece-1}.

 For any $t\in[0,T)$, $\xi\in
L^2_{\cF_t}(\Omega;H)$ and $u\in
L^2_{\dbF}(t,T;U)$,  consider the
following equation:
\bel{Used-in-later-sufficiency}\left\{\ba{ll}
\ds  d y =\big[\big(A+A_1 \big)y +B u \big]dr+\big(C y +D u\big) dW(r) &\mbox{ in }(t,T],\\
\ns\ds y(t)=\xi,
\ea\right.
\ee
which admits a unique solution
$y(\cd)\in
L^2_{\dbF}(\Omega;C([t,T];H))$ under
({\bf AS0}) and ({\bf AS2}). Let $\f\=u-\Th y$. We
can rewrite
\eqref{Used-in-later-sufficiency} as
\bel{Second-cX-equation}\left\{\ba{ll}
\ds  d y =\big[\big(A+A_\Th \big)y +B \f \big]dr
+\big[C_\Th y +D \f \big]dW(r) &\mbox{ in }(t,T], \\
\ns\ds y(t)=\xi.
\ea\right.
\ee
It follows from \eqref{SEE-cost-functional} that
\begin{equation}\label{3.4-eq24}
\begin{array}{ll}\ds
2\cJ(t,\xi;u(\cd))\\
\ns\ds = \dbE\lan G y(T),y(T)\ran_{H}+
\dbE\int_t^T \Big[\lan
Q(r)y(r),y(r)\ran_H
+\langle R(r)u(r),u(r)\ran_{U} \Big]dr\\
\ns\ds = \dbE\lan G y(T),y(T)\ran_{H}\\
\ns\ds \q+ \dbE\int_t^T \Big[\lan
Q(r)y(r),y(r)\ran_H +\langle
R(r)\big(\f(r)+\Th(r)
y(r)\big),\big(\f(r)+\Th(r)
y(r)\big)\ran_{U} \Big]dr\\
\ns\ds = \dbE\lan Gy(T),y(T)\ran_{H}+
\dbE\int_t^T \lan
[Q(r)+\Th(r)^*R(r)\Th(r)]y(r),y(r)\ran_H dr\\
\ns\ds\q + \dbE\int_t^T \big[ \langle
R(r) \f(r),\f(r)\ran_{ U } + 2\langle
R(r)\Th(r) y(r),\f(r)\ran_{ U }\big]
dr.
\end{array}
\end{equation}
Since $\f\in L^2_{\dbF}(0,T;U)$, and
$\wt U$ is a dense  subspace of $U$,
there exists a sequence of
$\{\f_n\}_{n=1}^\infty\subset
L^2_{\dbF}(0,T;\wt U)$ such that
\begin{equation}\label{10.25-eq1}
\lim_{n\to\infty}\f_n=\f \q\mbox{ in }
L^2_{\dbF}(0,T;U).
\end{equation}
Let $y_n$ be the solution of
(\ref{Second-cX-equation})
corresponding with $\f_n$. Similar to
the proof of  (\ref{3.14-eq9}), one can
show that
$$
\lim_{n\to\infty}y_n= y \q\mbox{ in }
L^2_{\dbF}(\Omega;C([0,T];H)).
$$
Hence,
\begin{equation}\label{3.4-eq24-1}
\begin{array}{ll}\ds
2\cJ(t,\xi;u(\cd))\\
\ns\ds =
\lim_{n\rightarrow\infty}\Big\{\dbE\lan
Gy_n(T),y_n(T)\ran_{H}+ \dbE\int_t^T
\lan
[Q(r)+\Th(r)^*R(r)\Th(r)]y_n(r),y_n(r)\ran_H dr\\
\ns\ds\qq\qq + \dbE\int_t^T \big[
\langle R(r) \f_n(r),\f_n(r)\ran_{U} +
2\langle R(r)\Th(r)
y_n(r),\f_n(r)\ran_{U}\big] dr\Big\}.
\end{array}
\end{equation}
By \eqref{Applied-to-SLQ} and \eqref{3.4-eq34},
we have that
\begin{eqnarray}\label{3.4-eq23}
&&\dbE\lan G y_n(T),y_n(T)\ran_{H}+ \dbE\int_t^T \lan [Q(r)+\Th(r)^*R(r)\Th(r)]y_n(r),y_n(r)\ran_H dr\nonumber\\
&&   =\dbE\lan \overline P(t)  \xi, \xi \ran_{H}+\dbE\int_t^T \lan \overline P(r)B(r)\f_n(r),y_n(r)\ran_{H} dr+\dbE\int_t^T \lan \overline P(r)y_n(r),B(r)\f_n(r)\ran_H dr\nonumber\\
&& \q +2\dbE\!\int_t^T\!\! \lan \overline
P(r)\big(C(r)+D(r)\Th(r)\big)y_n(r),D(r)\f_n(r)\ran
_{H} dr\\
&&\q  + \dbE\!\int_t^T\!\! \lan \overline P(r)D(r)\f_n(r),
D(r)\f_n(r)\ran_H dr+2\dbE\int_t^T \lan
D(r)\f_n(r),\overline \L(r)y_n(r)\ran_{H_\l,H_\l'} dr\nonumber\\
&& = \dbE\lan  \overline P(t) \xi, \xi\ran_{H}+\!\dbE\!\int_t^T\!\!
\lan \overline P(r)D(r)\f_n(r), D(r)\f_n(r)\ran_H dr\! -\!2\dbE\!
\int_t^T\!\! \lan R(r)\Th(r) y_n(r), \f_n(r) \ran_{H} dr.\nonumber
\end{eqnarray}
Combining \eqref{3.4-eq24-1} and
\eqref{3.4-eq23}, and recalling
\eqref{10.25-eq1}, we obtain that
\begin{eqnarray}\label{3.4-eq25}
&&
2\cJ(t,\xi;u(\cd))\nonumber\\
&=&\3n\lim_{n\rightarrow\infty}\!\Big( \dbE\lan \overline  P(t) \xi,
\xi \ran_{H}\! +\!\dbE\!\int_t^T\!\! \lan \overline P(r)D(r)\f_n(r),
D(r)\f_n(r)\ran_H dr \!+\!\dbE\!\int_t^T\!\! \langle R(r)
\f_n(r),\f_n(r)\ran_U   dr\Big) \nonumber\\
& =&\3n\dbE\lan \overline  P(t) \xi, \xi \ran_{H}
+\dbE\int_t^T \lan \overline K(r)\f(r),\f(r)\ran dr\\
& = &\3n\dbE\lan \overline P(t) \xi, \xi\ran_{H}
+\dbE\int_t^T\big\langle \overline K(u(r)-\Th (r) y(r)),u(r)-\Th(r)
y(r)\big\rangle_Udr.\nonumber
\end{eqnarray}
Noting that $\Th$ is an optimal feedback operator,  we get that
\begin{equation}\label{8.22-eq1}
\begin{array}{ll}\ds
\frac{1}{2}\mE\langle \overline P(t)\xi,\xi
\rangle_H=\cJ(t,\xi;\Th(\cd)\overline x(\cd))\leq
\cJ(t,\xi;u(\cd)),\qq\forall\; u(\cd)\in L^2_\dbF(t,T;U).
\end{array}
\end{equation}
Combining \eqref{3.4-eq25} and \eqref{8.22-eq1}, we get \eqref{Nece-1}.

\ms

The ``only if" part. Suppose $\Th\in
\Upsilon_2(H;U)$ is an operator
satisfying (\ref{Nece-1}) and
(\ref{Nece-2}). Following the exact
same procedures as the above Step 3, we
can show that $\Th$ is optimal. \endpf

\begin{remark}\label{Remark-novel-contribution}
We point out one essential difference from  the
Markovian study in   \cite[Theorem
3.3]{Wang-COCV-2020}, \cite[Lemma
2.8]{Lu-2020-COCV} due to the appearance of
random coefficients. To this end, let us take a
closer look at (\ref{A-key-contribution}). In
the Markovian setting, the expectation $\dbE$
can be omitted. Then by the separable property
of Hilbert space $H$, $U$ or the Euclidean
space, and the classical Lebesgue's differential
theorem, one can obtain the pointwise limit of
(\ref{A-key-contribution}). This method does not
work in our setting. To overcome this obstacle,
we introduce a proper form of $\xi$ in
\eqref{3.4-eq19}, and suitable Lebesgue's
differential theorem (Lemma
\ref{Lebesgue-differential-theorem}). To the
best of our knowledge, these two technologies
are even new in the counterpart finite
dimensional situation.
\end{remark}


\subsection{Proof of Theorem \ref{5.7-th1}}

Based on the preliminary results in the previous subsections, we are ready to prove the main conclusion of the current paper, i.e. Theorem \ref{5.7-th1}. To begin with, let us recall the following result.
\begin{lemma}\cite[Lemma 3.9]{Lu-Zhang-arxiv-2019}\label{lm11}
The set
$$
\begin{array}{ll}\ds
\big\{x_2(\cd)\;\big|\; x_2(\cd)\mbox{ solves }
\eqref{op-fsystem2}
\mbox{ with }t=0,\; \xi_2=0,\;  v_2=0\mbox{ and } u_2\in
L^4_{\dbF}(\Om;L^2(0,T;H)) \big\}
\end{array}
$$
is dense in $L^2_{\dbF}(0,T;H)$.
\end{lemma}

{\it Proof of Theorem
    \ref{5.7-th1}}\,: We first prove the ``only if" part. The basic idea is
to prove that the $H_\l$-transposition solution $(\overline
P,\overline\L)$ obtained in Theorem \ref{th1} is also the unique
transposition solution in terms of Definition \ref{4.8-def2}
satisfying (\ref{5.7-eq5}). We divide the proof of the ``only if"
part into two steps.

\ss

{\bf Step 1}. In this step, we prove
the existence of a transposition
solution.  From ({\bf AS2}) and
\eqref{Nece-1}, we know that
\begin{equation}\label{3.4-eq38}
\overline K>0 \mbox{ for a.e. } (t,\omega)\in [0,T]\times\Omega.
\end{equation}

By the same argument in {\bf Step 7} of \cite[Theorem
2.2]{Lu-Zhang-arxiv-2019},  we can show that the domain of
$\overline K(t,\om)^{-1}$ is dense in $U$ for a.e. $(t,\om)\in
[0,T]\times\Om$ and $\overline K(t,\om)^{-1}$ is a closed operator
for a.e. $(t,\om)\in (0,T)\times\Om$.

Let $\xi_1\in L^4_{\cF_t}(\Omega;H_\l)$
(\resp $\xi_2\in
L^4_{\cF_t}(\Omega;H_\l)$) and $u_1,
v_1\in
L^4_{\dbF}(\Omega;L^2(t,T;H_\l))$
(\resp $u_2, v_2\in
L^4_{\dbF}(\Omega;L^2(t,T;H_\l))$) in
\eqref{op-fsystem1} (\resp
\eqref{op-fsystem2}).  Let $\tilde
\xi_1=\xi_1$ (\resp $\tilde \xi_2 =
\xi_2$), $\tilde u_1 = u_1 - B\Th\tilde
x_1$ (\resp $\tilde u_2 = u_2 -
B\Th\tilde x_2$) and $\tilde v_1 = v_1
- D\Th\tilde x_1$ (\resp $\tilde v_2 =
v_2 - D\Th\tilde x_2$)  in
\eqref{op-fsystem3} (\resp
\eqref{op-fsystem4}). Then $\tilde x_1
= x_1$ and $\tilde x_2 = x_2$, where
$x_1$ and $x_2$ are solutions to
\eqref{op-fsystem1} and
\eqref{op-fsystem2}, respectively. From
\eqref{3.4-eq31}, we have that
\begin{eqnarray}\label{3.4-eq32}
&& \3n\3n\dbE\lan Gx_1(T),
x_2(T)\ran_{H}+ \dbE\int_t^T \lan
[Q(r)+\Th(r)^*R(r)\Th(r)] x_1(r),
x_2(r)\ran_H dr
\nonumber \\
&&\3n\3n =\dbE\lan \overline P(t) \xi_1, \xi_2\ran_{H}
+\dbE\int_t^T \lan \overline P(r)\big(u_1(r) - B(r)\Th(r) x_1(r)\big), x_2(r)\ran_{H} dr \nonumber \\
&&\3n\3n\q +\dbE\int_t^T \lan \overline P(r) x_1(r),\big(u_2(r) - B(r)\Th(r) x_2(r)\big)\ran_H dr\nonumber \\
&&\3n\3n \q +\dbE\int_t^T \lan \overline
P(r)\big(C(r)+D(r)\Th(r)\big) x_1(r),v_2(r) -
D(r)\Th(r) x_2(r) \ran _{H} dr\nonumber\\
&&\3n\3n \q  +\dbE\!\int_t^T\! \lan \overline P(r)\big(v_1(r) \!-\!
D(r)\Th(r) x_1(r)\big), \big(C(r)\!+\!D(r)\Th(r)\big)
x_2(r)\!+\!v_2(r) \!-\!
D(r)\Th(r) x_2(r)\ran_H dr\nonumber \\
&&\3n\3n\q +\dbE\int_t^T \lan v_1(r) -
D(r)\Th(r)
x_1(r),\overline \L^*(r) x_2(r)\ran_{H_\l,H_\l'} dr\nonumber \\
&&\3n\3n\q +\dbE\int_t^T \lan v_2(r) - D(r)\Th(r) x_2(r),\overline
\L(r) x_1(r)\ran_{H_\l,H_\l'}
dr \\
&&\3n\3n =  \mE\big\langle \overline P(t) \xi_{1},\xi_{2}
\big\rangle_{H} + \mE \int_t^T \big\langle \overline P(r)u_{1}(r),
x_{2}(r)\big\rangle_{H}dr + \mE \int_t^T \big\langle \overline
P(r)x_{1}(r),
u_{2}(r)\big\rangle_{H}dr \nonumber \\
&& \3n\3n \q  + \mE \int_t^T\big\langle \overline P(r)C(r)x_{1}(r),
v_{2}(r)\big\rangle_{H}dr + \mE
\int_t^T \big\langle \overline P(r)v_{1}(r), C(r)x_{2}(r)+v_{2}(r)\big\rangle_{H}dr\nonumber\\
&& \3n\3n \q + \mE \int_t^T \big\langle v_{1}(r), \overline
\L(r)x_2(r)\big\rangle_{H_\l,H_\l'}dr+ \mE \int_t^T \big\langle
\overline \L(r)x_1(r), v_{2}(r) \big\rangle_{
H_\l',H_\l}dr \nonumber \\
&& \3n\3n \q  -\mE \int_t^T \big\langle \Th(r) x_{1}(r),\big(B(r)^*
\overline P(r) +D(r)^* \overline P(r)C(r) +
D(r)^*\overline \L(r)\big)x_{2}(r) \big\rangle_{H}dr\nonumber  \\
&& \3n\3n \q  -\mE \int_t^T \big\langle \Th(r) x_{2}(r),\big(B(r)^*
\overline P(r) +D(r)^* \overline P(r)C(r) +
D(r)^*\overline \L(r)\big)x_{1}(r) \big\rangle_{H}dr\nonumber  \\
&& \3n\3n \q  -  \dbE \int_t^T \lan \overline P(r) D(r)\Th(r)
x_1(r),D(r)\Th(r) x_2(r) \ran _{H} dr \nonumber
\end{eqnarray}
Noting
\begin{equation*}
\begin{array}{ll}\ds
\mE \int_t^T \big\langle \overline K(r)^{-1} \overline L(r)
x_{1}(r), \overline L(r)x_{2}(r) \big\rangle_{H}dr\\
\ns\ds = - \mE \int_t^T \big\langle \big(B(r)^* \overline P(r)
+D(r)^* \overline P(r)C(r) + D(r)^*  \overline
\L(r)\big)x_{2}(r),\Th(r) x_{1}(r) \big\rangle_{H}dr
\end{array}
\end{equation*}
and
\begin{equation*}\label{3.4-eq33}
\begin{array}{ll}\ds
\dbE\int_t^T \lan \Th(r)^*R(r)\Th(r) x_1(r), x_2(r)\ran_H dr + \dbE
\int_t^T \lan \overline P(r) D(r)\Th(r) x_1(r),D(r)\Th(r) x_2(r)
\ran _{H}
dr\\
\ns\ds = \dbE\int_t^T \lan \big(R(r)+D(r)^*\overline P(r)
D(r)\big)\Th(r) x_1(r), \Th(r) x_2(r)\ran_H dr\\
\ns\ds = -\mE \int_t^T \big\langle \big(B(r)^* \overline P(r)
+D(r)^* \overline P(r)C(r) + D(r)^*\overline
\L(r)\big)x_{1}(r),\Th(r) x_{2}(r) \big\rangle_{H}dr,
\end{array}
\end{equation*}
we get from \eqref{3.4-eq32} that
\begin{eqnarray}\label{3.4-eq36}
&& \3n\3n\3n \mE\langle Gx_{1}(T),x_{2}(T)\rangle_{H} +\mE \int_t^T
\big\langle Q(r) x_{1}(r),
x_{2}(r) \big\rangle_{H}dr\nonumber\\
&& \3n\3n\3n\q - \mE  \int_t^T \big\langle \overline K(r)^{-1}
\overline L(r) x_{1}(r), \overline L(r)x_{2}(r)
\big\rangle_{H}dr\nonumber
\\
&&\3n\3n\3n = \mE\big\langle \overline P(t) \xi_{1},\xi_{2}
\big\rangle_{H} + \mE \int_t^T \big\langle \overline P(r)u_{1}(r),
x_{2}(r)\big\rangle_{H}dr + \mE \int_t^T \big\langle \overline
P(r)x_{1}(r),
u_{2}(r)\big\rangle_{H}dr \\
&& \3n\3n\3n\q  + \mE \int_t^T\big\langle \overline
P(r)C(r)x_{1}(r), v_{2}(r)\big\rangle_{H}dr + \mE
\int_t^T \big\langle \overline P(r)v_{1}(r), C(r)x_{2}(r)+v_{2}(r)\big\rangle_{H}dr\nonumber\\
&& \3n\3n\3n\q + \mE \int_t^T \big\langle v_{1}(r), \overline
\L(r)x_2(r)\big\rangle_{ H_\l,H_\l'}dr+ \mE \int_t^T \big\langle
\overline \L(r)x_1(r), v_{2}(r) \big\rangle_{
H_\l',H_\l}dr.\nonumber
\end{eqnarray}
%
%
%

Next we look at the case that $\xi_j\in
L^4_{\cF_t}(\Omega;H)$, $u_j\in
L^4_{\dbF}(\Omega;L^2(t,T;H))$ and $ v_j\in
L^4_{\dbF}(\Omega;L^2(t,T;$ $H_\l))$ ($j=1,2$).
Recall that $H_\l$ is dense in $H$, for $j=1,2$,
we can find $\{\xi_{j,n}\}_{n= 1}^\infty\subset
L^4_{\cF_t}(\Omega;H_\l)$ and $\{u_{j,n}\}_{n=
1}^\infty\subset
L^4_{\dbF}(\Omega;L^2(\Om,H_\l))$ such that
$$\ba{ll}
\ns\ds \lim_{n\rightarrow \infty}
|\xi_{j,n}-\xi_j|_{L^4_{\cF_t}(\Om;H)} =0,\\
\ns\ds \lim_{n\rightarrow \infty}
|u_{j,n}-u_j|_{
L^4_{\dbF}(\Omega;L^2(t,T;H))}
=0.
\ea
$$
Let $x_{1,n}$ (resp. $x_{2,n}$) be the
solution to (\ref{op-fsystem1}), (\resp
(\ref{op-fsystem2}))  with $\xi_{1}$
(\resp $\xi_2$) and $u_1$ (\resp $u_2$)
replaced  by $\xi_{1,n}$ (\resp
$\xi_{2,n}$) and $u_{1,n}$ (resp.
$u_{2,n}$) and $v_1$  (\resp $v_2$).
Similar to the proof of
\eqref{3.14-eq9}, for $j=1,2$, we can
get that
$$
\lim_{n\to \infty} x_{j,n}=x_j \q
\mbox{ in }L^2_{\dbF}(\Omega;C(t,T;H)).
$$
Then similar to the proof of
\eqref{3.14-eq31}, we can get that
(\ref{3.4-eq36}) holds with $\xi_j\in
L^4_{\cF_t}(\Omega;H)$, $u_j\in
L^4_{\dbF}(\Omega;L^2(t,T;H))$ and $
v_j\in
L^4_{\dbF}(\Omega;L^2(t,T;H_\l))$
($j=1,2$). Hence, $(\overline P, \overline \L)$
is a transposition solution to
\eqref{5.5-eq6}.

\vspace{0.2cm}

{\bf Step 2}. In this step, we prove
the uniqueness of the transposition
solution to \eqref{5.5-eq6}.

Assume that
$$(P_1(\cd),\Lambda_1(\cd)),
(P_2(\cd),\Lambda_{2}(\cd))\in
C_{\dbF,w}([0,T];
L^{\infty}(\Om;\cL(H))) \times
L^2_{\dbF,D,w}(0,T;\cL(H))$$ are two
transposition solutions to
\eqref{5.5-eq6}.

From \eqref{8.22-eq1}, we have that for
any $s\in[0,T)$ and $\eta\in
L^2_{\cF_s}(\Om;H)$,
\begin{equation}\label{8.22-eq2-1}
\frac{1}{2}\mE\langle P_1(s)\eta,\eta
\rangle_H=\cJ(s,\eta;\Th(\cd)x(\cd))=\frac{1}{2}\mE\langle
P_2(s)\eta,\eta \rangle_H.
\end{equation}
Thus, for any $\xi,\eta\in
L^2_{\cF_s}(\Om;H)$, we have
that\vspace{0.4mm}
$$
\mE\langle P_1(s)(\eta+\xi),\eta+\xi
\rangle_H =\mE\langle P_2(s)(\eta+\xi),
\eta+\xi \rangle_H,
$$
and
$$
\mE\langle P_1(s)(\eta -\xi),\eta -\xi
\rangle_H =\mE\langle P_2(s)(\eta
-\xi),\eta -\xi \rangle_H.\vspace{3mm}
$$
These, together with
$P_1(\cd)=P_1(\cd)^*$ and
$P_2(\cd)=P_2(\cd)^*$, imply
that\vspace{1mm}
\begin{equation}\label{8.22-eq3}
\mE\langle P_1(s)\eta,\xi \rangle_H
=\mE\langle P_2(s)\eta,\xi
\rangle_H,\qq \forall\;\xi,\eta\in
L^2_{\cF_s}(\Om;H).\vspace{2mm}
\end{equation}
Hence,
$$
P_1(s)=P_2(s)\, \mbox{ for any }  s\in
[0,T], \; \dbP\mbox{-a.s.}
$$

\vspace{0.1cm}

Let $v_2=0$ in \eqref{op-fsystem2}. By
\eqref{6.18-eq1} and noting\vspace{1mm}
$$
\big\langle K(\cd)^{-1} L(\cd)
x_{1}(\cd), L(\cd)x_{2}(\cd)
\big\rangle_{H}=-\big\langle \Th(\cd)
x_{1}(\cd), K(\cd)\Th(\cd)x_{2}(\cd)
\big\rangle_{H},
$$
we see that for any $\xi_1\in
L^2_{\cF_0}(\Om;H)$, $u_1\in
L^4_\dbF(\Om;L^2(0,T;H))$ and $v_1 \in
L^4_\dbF(\Om;L^2(0,T;H_\l))$,
\begin{equation*}\label{8.22-eq9-1}
0=\mE \int_0^T \big\langle v_1(t),
\big(\Lambda_1(t) -
\Lambda_2(t)\big)x_2(t)\big\rangle_{H_\l,H_\l'}dt.
\end{equation*}
Consequently,
\begin{equation*}\label{9.7-eq13-1}
\big(\Lambda_1 - \Lambda_2 \big)x_2 = 0
\;\mbox{ in }\;
L^{\frac{4}{3}}_\dbF(\Om;L^2(0,T;H_\l')).
\end{equation*}
This, together with Lemma \ref{lm11},
implies that
$$\Lambda_1=\Lambda_2 \mbox{ in
}L^2_\dbF(0,T;\cL_2(H;H_\l')).$$ Hence,
the desired uniqueness follows.

\ms

Next, we prove the ``if" part. The proof is
similar to  {\bf  Step 3} in the proof for Lemma
\ref{th1} (see also \cite[Theorem
2.1]{Lu-Zhang-arxiv-2019}). For readers'
convenience, we give a sketch.

For any $u\in L^2_{\dbF}(0,T;U)$, since
$\wt U$  is a dense subspace of $U$,
 there exists a sequence
$\{u_n\}_{n=1}^\infty\subset
L^2_{\dbF}(0,T;\wt U)$ such that
\begin{equation}\label{10.25-eq2}
\lim_{n\to\infty}u_n= u \q\mbox{ in
}L^2_{\dbF}(0,T;U).
\end{equation}
Let $y_n$ be the solution to
\eqref{Used-in-later-sufficiency}
corresponding with $u_n$ and a
$\xi\in L^2_{\cF_t}(\Om;H)$. Then,
similar to the proof of
(\ref{3.14-eq9}),  we can obtain that
\begin{equation}\label{10.25-eq3}
\lim_{n\to\infty}y_n= y \q\mbox{ in
}L^2_{\dbF}(\Omega;C([t,T];H)).
\end{equation}
Hence
\begin{equation}\label{10.25-eq4}\ba{ll}
\ns\ds
2\cJ(t,\xi;u(\cd))\\
\ns\ds = \lim_{n\rightarrow\infty}
\dbE\lan Gy_n(T),y_n(T)\ran_{H}+
\dbE\int_t^T\big[ \lan Q(r)
y_n(r),y_n(r)\ran_H+ \langle R(r)
u_n(r),u_n(r)\ran_{U} \big]dr.
\ea
\end{equation}
By the definition of transposition
solution in Definition \ref{4.8-def2},
we have
\begin{eqnarray}\label{10.25-eq5}
&&\3n\3n\3n\dbE\lan G y_n(T),y_n(T)\ran_{H}+
\dbE\int_t^T
\big[ \lan Q(r) y_n(r),y_n(r)\ran_H +\lan R(r)u_n(r),u_n(r)\ran_{U}\big]dr \nonumber \\
&=&\3n\dbE\lan   P(t)  \xi, \xi
\ran_{H}+ 2\dbE\int_t^T \lan
u_n(r),L(r)y_n(r)\ran_{U} dr
+ \dbE\!\int_t^T\!\!
\lan K(r)u_n(r), u_n(r)\ran_{U} dr \nonumber\\
&&  +\dbE\int_t^T \lan K(r)^{-1}L(r)y_n(r),L(r)y_n(r)\ran_H dr\\
& =&\3n\dbE\lan   P(t)  \xi, \xi
\ran_{H}+ \dbE\int_t^T \lan
K(r)\big[u_n(r)+K(r)^{-1}L(r)y_n(r)\big],u_n(r)+K(r)^{-1}L(r)y_n(r)\ran_{
U} dr.\nonumber
\end{eqnarray}
From \eqref{10.25-eq2} and
\eqref{10.25-eq3}, we have that
\begin{equation}\label{10.25-eq6}\ba{ll}
\ns\ds \lim_{n\to\infty}\dbE\int_t^T \lan K(r)\big[u_n(r)+K(r)^{-1}L(r)y_n(r)\big],u_n(r)+K(r)^{-1}L(r)y_n(r)\ran_{  U} dr\\
\ns\ds =\dbE\int_t^T \lan
K(r)\big[u(r)+K(r)^{-1}L(r)y(r)\big],u(r)+K(r)^{-1}L(r)y(r)\ran_{
U} dr.
\ea
\end{equation}
Combining
\eqref{10.25-eq4}--\eqref{10.25-eq6},
we get that
$$\ba{ll}
\ns\ds 2\cJ(t,\xi;u(\cd)) =\dbE\int_t^T
\lan
K(r)\big[u(r)+K(r)^{-1}L(r)y(r)\big],u(r)+K(r)^{-1}L(r)y(r)\ran_{
U} dr.
\ea
$$
This implies that
$$\ba{ll}
\ns\ds 2\cJ(t,\xi;u(\cd))\geq 2\cJ(t,\xi;\bar u(\cd)),\ \ \forall u\in L^2_{\dbF}(t,T;U),
\ea
$$
with $\bar u\deq -K^{-1}L y$. Hence
$-K^{-1}L$ is an optimal feedback
operator.

At last, we prove the uniqueness of the
optimal feedback operator. Suppose that
$\Th'$ is another optimal feedback
operator. According to the proof
of ``only if" part of Theorem
\ref{5.7-th1}, Riccati equation
\eqref{5.5-eq6} admits a transposition
solution, denoted by
$\big(P'(\cd),\L'(\cd)\big)$. By the
uniqueness of the Riccati equation
\eqref{5.5-eq6}, we have
$\big(P'(\cd),\L'(\cd)\big)=\big(P(\cd),\L(\cd)\big)$,
which leads to $\Th=\Th'$. This
completes the proof of Theorem
\ref{5.7-th1}.
\endpf

\section{SLQs for stochastic parabolic equations}\label{sec-example}

Stochastic parabolic equations are
widely used to describe diffusion
processes under the perturbations of
random noises (e.g.,
\cite{Kotelenez-2008}). In this
section, we consider a linear quadratic
control problem for such equation.

Let $m\in\dbN$ and $\cO\subset\dbR^m$
be a bounded domain with a $C^\infty$
boundary $\pa\cO$. Consider the
following control system:
\begin{equation}\label{10.25-eq25}
\begin{cases}
dy -   \D ydt = (a_1y+ b_1u)dt + (a_2y + b_2
u)dW(t) &\mbox{ in } \cO\times (0,T),\\
\ns\ds y=0 &\mbox{ on }\pa\cO\times (0,T),\\
\ns\ds y(0)=y_0 &\mbox{ in } \cO,
\end{cases}
\end{equation}
with the following cost functional
$$
\cJ(0,y_0;u)\triangleq
\mE\int_0^T\int_\cO (q|y|^2 +
r|u|^2)dxdt + \mE\int_\cO g|y(T)|^2dx.
$$
Here $y_0\in L^2(\cO)$ and $u\in
L^2_\dbF(0,T;L^2(\cO))$ is the control
variable. The conditions on the
coefficients will be given below.

Our optimal control problem is as
follows:

\ms

\no\bf Problem (sSLQ). \rm For each
$y_0\in L^2(\cO)$, find a $\bar
u(\cd)\in L^2_\dbF(0,T;L^2(\cO))$ such
that\vspace{0.1cm}
\begin{equation}\label{10.25-eq26}
\cJ\big(0,y_0;\bar
u(\cd)\big)=\inf_{u(\cd)\in
L^2_\dbF(0,T;L^2(\cO))}\cJ\big(0,y_0;u(\cd)\big).
\end{equation}

 Problem (sSLQ) is a concrete example of Problem
(SLQ) with the following setting:

\begin{itemize}
  \item $H=U=L^2(\cO)$;
  \item The operator $A$ is
defined as follows:
\begin{equation}\label{10.25-eq26-1}
\begin{cases}
D(A)=H^2(\cO)\cap H_0^1(\cO),\\ \ns\ds
A\f =  \D\f,\q \forall\; \f\in D(A);
\end{cases}
\end{equation}
\item $A_1y=a_1y$, $Bu=b_1u$, $Cy= a_2y$ and $Du=b_2u$;
\item The operators $Q$, $R$ and $G$ are given by
$$
\begin{cases}\ds
\lan Q(t) y,y\ran_H = \int_\cO q(t)|y|^2dx, \q
\lan R(t) u,u\ran_H = \int_\cO
r(t)|u|^2dx,\q \ae t\in [0,T],\\
\ns\ds\lan G y(T),y(T)\ran_H = \int_\cO
g|y(T)|^2dx.
\end{cases}
$$
\end{itemize}

From \eqref{10.25-eq26-1}, we get ({\bf AS0}).
To guarantee that ({\bf AS1})--({\bf
AS4}) hold, we assume the coefficients
fulfill the following conditions:
\begin{equation}\label{5.27-eq4}
\begin{cases}
\ds  a_1 \in
L^1_{\dbF}(0,T;L^\infty(\Om;C^{2m}(\overline\cO))),\q
b_1\in
L^\infty_{\dbF}(\Om;L^2(0,T;C^{2m}(\overline\cO))),\\
\ns\ds a_2,q\in
L^2_{\dbF}(0,T;L^\infty(\Om;C^{2m}(\overline\cO))),\q
b_2,  r \in
L^\infty_{\dbF}(0,T;C^{2m}(\overline\cO)), \\
\ns\ds g\in
L^\infty_{\cF_T}(\Om;C^{2m}(\overline\cO)),\q
q\geq 0,\;\;r>0,\;\;g\geq 0\;\mbox{ for a.e.
}(t,x,\om)\in [0,T]\times \cO\times\Om.
\end{cases}
\end{equation}

From \eqref{5.27-eq4}, one can check
that ({\bf AS1}) holds.

Write  $\{ \mu_j\}_{j=1}^\infty$ for the
eigenvalues of $A$ and $\{e_j\}_{j=1}^\infty$
the corresponding eigenvectors satisfying
$|e_j|_{L^2(\cO)}=1$ for all $j\in\dbN$. It is
well known that $\{e_j\}_{j=1}^\infty$
constitutes an orthonormal basis of $ L^2(\cO)$.
Hence, ({\bf AS2}) holds.

Let $\l_j = |\mu_j|^{-\frac{m}{2}}$ for
$j\in\dbN$ and
$\l=\{\l_j\}_{j=1}^\infty$. Then
$H_\l'$ is the completion of the
Hilbert space $L^2(\cO)$ with respect
to the norm
$$
\begin{array}{ll}\ds
|f|_{H_\l'}  = \sqrt{\sum_{j=1}^\infty
|\mu_j|^{-m}(1+|\mu_j|)^{-2}
|f_{j}|^2},  \qq \mbox{ for } f
=\sum_{j=1}^\infty f_{j} e_j \in
L^2(\cO).
\end{array}
$$
By the asymptotic distribution of
eigenvalues of $A$  (e.g.,
\cite[Chapter 1, Theorem 1.2.1 and
Remark 1.2.2]{Safarov}), we have
$\mu_j=\cC_1 j^{\frac{2}{m}}+O(1)$ for
some constant $\cC_1>0$. Therefore,
$$
\begin{array}{ll}\ds
|I|^2_{\cL_2(H;H_\l')} =
\sum_{j=1}^\infty \lan
e_j,e_j\ran_{H_\l'} \leq
4\sum_{j=1}^\infty |\mu_j|^{-n}
<\infty.
\end{array}
$$
Hence, the embedding from $ L^2(\cO)$
to $H_\l'$ is Hilbert-Schmidt. From the
definition of $H_\l'$, it follows that
$H_\l=D(A^{\frac{m}{2}+1})$.

Let $\a\in C^{2m}(\overline\cO)$. For any $f\in H_\l$, one has
$$
\begin{array}{ll}\ds
|\a f|_{H_\l}\3n&\ds=|\a
f|_{D(A^{{\frac{m}{2}+1}})} = \sqrt{|\a
f|_{L^2(\cO)}^2+|A^{{\frac{m}{2}+1}}
(\a f)|_{L^2(\cO)}^2}\\
\ns&\ds\leq \cC|\a|_{C^{2m}(\overline\cO)}|f|_{D(A^{\frac{m}{2}+1})} =
\cC|\a|_{C^{2m}(\overline\cO)}|f|_{H_\l}.
\end{array}
$$
From this, and recalling
\eqref{5.27-eq4}, we conclude that
$A_1\in
L^1_\dbF(0,T;L^\infty(\Om;\cL(H_\l)))$,
$C\in L^2_\dbF(0,T;$
$L^\infty(\Om;\cL(H_\l)))$, $G\in
L^\infty_{\cF_T}(\Om;\cL(H_\l))$ and
$Q\in
L^\infty_\dbF(\Om;L^2(0,T;\cL(H_\l)))$.
Similarly, we can prove that $A_1 \in
L^1_\dbF(0,T;L^\infty(\Om;\cL(H_\l')))$,
$C \in L^2_\dbF(0,T;L^\infty(\Om;$
$\cL(H_\l')))$, $G \in
L^\infty_{\cF_T}(\Om;\cL(H_\l'))$ and
$Q\in
L^\infty_\dbF(\Om;L^2(0,T;\cL(H_\l')))$.
Thus, ({\bf AS3}) holds.

\ss

Let $\wt U  =D(A^{\frac{m}{2}+1})$. Then $\wt U$
is dense in $L^2(\cO)$.  Let $\a\in
C^{2m}(\overline\cO)$. For any $f\in \wt U$, one
has
$$
\begin{array}{ll}\ds
|\a f|_{H_\l}\3n&\ds=|\a
f|_{D(A^{\frac{m}{2}+1})} =
|A^{\frac{m}{2}+1}
(\a f)|_{L^2(\cO)}\\
\ns&\ds\leq |\a|_{C^{2m}(\overline\cO)}|f|_{D(A^{\frac{m}{2}+1})} =
|\a|_{C^{2m}(\overline\cO)}|f|_{\wt U}.
\end{array}
$$
From this, and using \eqref{5.27-eq4}, we find
that $B \in L^\infty_\dbF(\Om;L^2(0,T;\cL(\wt
U;H_\l)))$, $D\in L^\infty_\dbF(0,T;$ $\cL(\wt
U;H_\l))$ and $R\in L^\infty_\dbF(0,T;\cL(\wt
U))$. Similarly, we can prove that $B \in
L^\infty_\dbF(\Om;L^2(0,T;\cL(\wt U';H_\l')))$,
$D\in L^\infty_\dbF(0,T;\cL(\wt U';H_\l'))$ and
$R\in L^\infty_\dbF(0,T;\cL(\wt U'))$.
Therefore, ({\bf AS4}) holds.

\end{document}